\documentclass{amsart}
\usepackage{lmodern}
\usepackage{amsmath, amsthm, amssymb, amsfonts}
\usepackage[normalem]{ulem}
\usepackage{hyperref}
\usepackage[all,cmtip]{xy}
\usepackage{verbatim}
\usepackage{nccmath}
\usepackage{stmaryrd}
\usepackage{caption}
\usepackage{mathrsfs}
\usepackage{mathtools}
\usepackage{esvect}
\usepackage{cite}
\usepackage{bbm}
\usepackage{eucal}

\usepackage{mathbbol}
\usepackage{tabularx}
\usepackage[toc,page]{appendix}

\usepackage{tikz-cd}

\theoremstyle{plain}
\newtheorem{thm}{Theorem}[section]
\newtheorem{cor}[thm]{Corollary}

\newtheorem{lemma}[thm]{Lemma}
\newtheorem{prop}[thm]{Proposition}

\newtheorem*{thmn}{Theorem}
\newtheorem*{corn}{Corollary}

\numberwithin{equation}{section}

\renewcommand{\arraystretch}{2}

\newcommand{\dashto}{%
	\mathrel{\tikz[baseline=-0.7ex]
		\draw[dashed,->] (0,0) -- (1.2em,0);}
}

\theoremstyle{definition}
\newtheorem{defn}[thm]{Definition}

\theoremstyle{remark}
\newtheorem{rmk}[thm]{Remark}

\newcommand{\BA}{{\mathbb{A}}}

\newcommand{\BC}{{\mathbb{C}}}

\newcommand{\BE}{{\mathbb{E}}}

\newcommand{\BL}{{\mathbb{L}}}
\newcommand{\BM}{{\mathbb{M}}}
\newcommand{\BN}{{\mathbb{N}}}

\newcommand{\BP}{{\mathbb{P}}}
\newcommand{\BQ}{{\mathbb{Q}}}
\newcommand{\BR}{{\mathbb{R}}}

\newcommand{\BZ}{{\mathbb{Z}}}

\newcommand{\CC}{{\mathcal C}}

\newcommand{\CE}{{\mathcal E}}

\newcommand{\CI}{{\mathcal I}}

\newcommand{\CL}{{\mathcal L}}

\newcommand{\CN}{{\mathcal N}}
\newcommand{\CO}{{\mathcal O}}
\newcommand{\CP}{{\mathcal P}}

\newcommand{\CU}{{\mathcal U}}

\newcommand{\CW}{{\mathcal W}}
\newcommand{\CX}{{\mathcal X}}

\newcommand{\FM}{{\mathfrak{M}}}

\newcommand{\Fz}{{\mathfrak{z}}}

\newcommand{\FW}{{\mathfrak{W}}}

\newcommand{\etalt}{{\mathrm{mult}}}

\newcommand{\di}{\mathrm{dim}}

\DeclareFontFamily{OT1}{rsfs}{}
\DeclareFontShape{OT1}{rsfs}{n}{it}{<-> rsfs10}{}
\DeclareMathAlphabet{\curly}{OT1}{rsfs}{n}{it}

\newcommand{\p}{\mathbb{P}}

\newcommand\Spec{\operatorname{Spec}}

\newcommand{\Mbar}{{\overline M}}

\newcommand{\Sym}{{\mathrm{Sym}}}

\newcommand{\ev}{{\mathrm{ev}}}

\newcommand{\Aut}{\mathrm{Aut}}

\newcommand{\FC}{\mathfrak{C}}

\begin{document}
	
	\title[Unramified Gromov--Witten and Gopakumar--Vafa invariants]
	{Unramified Gromov--Witten  and Gopakumar--Vafa invariants}

\author{Denis Nesterov}
\address{ETH Z\"urich, Departement Mathematik}
\email{denis.nesterov@math.ethz.ch}

\date{}

	\maketitle

	\begin{abstract}  Kim, Kresch, and Oh defined unramified Gromov–Witten invariants. For threefolds, Pandharipande conjectured that these invariants are equal to the Gopakumar–Vafa (BPS) invariants for Fano and primitive Calabi–Yau curve classes. We prove the conjecture using a wall-crossing method developed here. As a result, we obtain an algebro-geometric construction of Gopakumar–Vafa invariants in these cases.
	\end{abstract}

		\setcounter{tocdepth}{1}
	\tableofcontents
	\section{Introduction}
\subsection{Overview} In \cite{KKO}, Kim, Kresch, and Oh defined unramified Gromov--Witten theory, a higher-dimensional generalization of the theory of admissible covers developed by Harris and Mumford \cite{HM}. In this work, we introduce a stability condition, called $\epsilon$-\textit{unramification} (Definition	\ref{weightp}), that interpolates between standard Gromov--Witten theory and unramified Gromov--Witten theory. Using Zhou's theory of entangled tails \cite{YZ} in the setting of Fulton--MacPherson degenerations \cite{FM}, we establish a wall-crossing formula (Theorem \ref{wallcrossing}), relating Gromov--Witten invariants to unramified Gromov--Witten invariants, with the wall-crossing contributions expressed as Hodge integrals over moduli spaces of stable curves \cite{PandDeg, FP}. In the case of threefolds, we compute this wall-crossing formula explicitly, obtaining the following theorem.
	\begin{thmn}[Theorem \ref{pand}] \label{mainthm} Let $X$ be a smooth projective threefold.  There is a wall-crossing formula relating  standard Gromov--Witten (GW) invariants and unramified Gromov--Witten (uGW) invariants  of $X$ in a curve class $\beta \in H_2(X)$, 
		\[ \sum_g	\langle \gamma_1 \cdot \! \cdot \! \cdot\gamma_n  \rangle^{GW}_{g,\beta} u^{2g-2}=  \sum_g	\langle \gamma_1  \cdot \! \cdot \! \cdot\gamma_n  \rangle^{uGW}_{g,\beta} u^{2g-2} \left(\frac{\mathrm{sin}(u/2)}{u/2}\right)^{2g-2+\beta\cdot \mathrm{c}_1(X)}, \]
	where $g$ is the genus of curves, and  $\gamma_i \in H^*(X)$ are  cohomology classes. 
	\end{thmn}

On the one hand, the formula above is a three-dimensional analogue of the Gromov--Witten/Hurwitz correspondence of Okounkov and Pandharipande \cite{OPcom} for primary insertions.  On the other, it exactly matches the formula  relating Gromov--Witten and Gopakumar--Vafa invariants \cite[Section 3]{Pand3} if the class $\beta$ satisfies one of the following conditions:
\begin{itemize}
	\item  $\beta\cdot \mathrm{c}_1(X)>0$, or
\item   $\beta$ is primitive\footnote{A curve class which is not a multiple of another curve class.} and $\beta\cdot \mathrm{c}_1(X)=0$.
\end{itemize}
We call the former a \textit{Fano class}, while the latter a primitive  \textit{Calabi--Yau class}. 
We thereby conclude that unramified Gromov--Witten invariants  and Gopakumar--Vafa invariants are equal in these cases, as was conjectured by Pandharipande \cite[Conjecture 5.2.1]{KKO}. 

\begin{corn} \label{pandha}
		Let $X$ be a smooth projective threefold. If  $\beta$ is a Fano class or a primitive Calabi--Yau class, then unramified Gromov--Witten invariants and Gopakumar--Vafa invariants of $X$ in the class $\beta$ are equal. 
	
	\end{corn}

In particular, this provides an algebro-geometric construction of Gopakumar--Vafa invariants in terms of unramified Gromov--Witten invariants for such curve classes.\footnote{If $\beta\cdot \mathrm{c}_1(X)<0$, the corresponding invariants vanish by the dimension constraint.} 

The wall-crossing method developed here applies in arbitrary dimension. We refer to \cite[Section 4]{10A} for an application of Theorem \ref{wallcrossing} in the case of abelian surfaces, and to \cite{NSc} for the case of target curves.

\subsection{Gopakumar--Vafa invariants} From the perspective of algebraic geometry and curve-counting theories, Gopakumar–Vafa invariants come closest to honest curve counts as the theory of virtual fundamental classes allows. The other two major curve-counting theories -- Donaldson–Thomas theory \cite{Th} and Gromov–Witten theory \cite{Kon} -- have evident moduli-theoretic defects that prevent their invariants from being interpreted as genuine curve counts, except in a few special cases. More specifically, the moduli spaces used to define these theories are \textit{overcompactified}, leading to \textit{overcounting} in the resulting invariants.

In  Donaldson--Thomas theory, one must include embedded points and generically non-reduced curves in the moduli space, while in the case of Gromov--Witten theory there are maps with contracted components and multiple covers. For example, if we want to count the number of degree $1$ genus $g$ curves in $\p^3$ passing through two generic points, denoted by $n_{g,1}$, then the intuitive answer is:
\[ 
n_{g,1}=\begin{cases}
	1 & \text{if }g=0,\\
	0 & \text{otherwise}.
	\end{cases}
\]

Gopakumar--Vafa invariants are  equal to $n_{g,1}$.  However, Donaldson--Thomas and Gromov--Witten invariants are very far from these simple numbers. The latter are given by the Theorem above. 

The problem with Gopakumar--Vafa invariants is that they are  notoriously difficult to define. Nevertheless, we know what these invariants should be, as they admit indirect definitions. For example, one can use the formula from the Theorem  as a definition of invariants on the right-hand side, or a similar formula from \cite[Section 3]{Pand3}   in the case of an arbitrary Calabi--Yau class.  However, a direct mathematical construction has been very much desired.  

One approach is given by a sheaf--theoretic perspective and is inspired by the work of Gopakumar and Vafa themselves \cite{GV1, GV2}.  It was pursued  by many \cite{HST, Katz, KiLi}, culminating in Maulik--Toda's proposal \cite{MT}, which has far-reaching connections to representation theory, topology, and mathematical physics.  However, Maulik--Toda invariants have been shown to be equal to Gopakumar--Vafa invariants only in a few cases, primarily through explicit computations, and they are restricted to Calabi--Yau threefolds.

Moduli spaces of maps provide another solution. The idea is to improve the moduli spaces by removing unwanted maps, namely those with contracted components or arising as multiple covers. This can be achieved in several ways. The aforementioned unramified Gromov–Witten theory of \cite{KKO} leads to one such approach. Another is reduced Gromov–Witten theory, initiated in \cite{VZ08, VZ07} by Vakil and Zinger, and further constructed for  complete intersections in \cite{RMMP}. We refer to \cite{RMMP} for a detailed account of the history of reduced Gromov–Witten theory and related developments.

To the best of the author’s knowledge, both unramified Gromov--Witten and reduced Gromov–Witten theories, at least in their current form, address only contracted components, but not multiple covers. However, unlike the sheaf-theoretic approach, they are not restricted to Calabi–Yau threefolds.

If $\beta$ is a Fano or primitive Calabi–Yau class, then multiple covers do not contribute to the primary invariants. Hence, only contracted components occur. By construction, unramified Gromov--Witten theory removes them by introducing degenerations of the target, thereby yielding a more efficient compactification of the moduli space of maps, whose intersection theory gives rise to Gopakumar--Vafa invariants by the Corollary above. The wall-crossing, provided by $\epsilon$-unramification, expresses the difference between Gromov--Witten and unramified Gromov--Witten theories via the local contributions of the contracted components, which take the form of Hodge integrals \cite{PandDeg, FP}.

For one-dimensional targets, an analogous statement holds in its strongest form: the moduli spaces of admissible covers \cite{HM} are smooth after normalization, and every map lies in the closure of the locus of maps from smooth curves.

The relation between reduced Gromov--Witten and Gopakumar--Vafa invariants remains an open problem except in a few cases; see \cite{RMMP}.  On the other hand, a direct connection between unramified and reduced Gromov–Witten theories arises from the projection from the moduli space of unramified maps to the moduli space of stable maps. In light of the Corollary above and \cite[Conjecture 1.0.3]{RMMP}, it is natural to expect that the pushforward of the virtual fundamental class of the moduli space of unramified maps is closely related to the virtual fundamental class of the main component\footnote{This is the component (or components) given by the closure of the locus of maps from smooth curves; see \cite{RMMP} for further details.} of the moduli space of stable maps. Moreover, for suitable targets and curve classes, this projection is expected to describe the moduli space of unramified maps as a sequence of blow-ups of the main component. 
 
 
\subsection{Integrality and finiteness} The distinguishing features of Gopakumar--Vafa invariants are  \textit{integrality} and \textit{finiteness}. Integrality refers to the property that invariants are integers, while finiteness is the vanishing of invariants for a sufficiently large genus. Integrality and finiteness  for Fano classes were proved in \cite{Zing, DW}, using Zinger's signed counts of embedded $J$-holomorphic curves.   For Calabi--Yau classes, they were proved in \cite{IP,DIW}, using Ionel--Parker's notion of clusters. These results,  together with the Theorem, imply that unramified Gromov--Witten invariants are integers which vanish for a sufficiently large genus if the class $\beta$ satisfies the assumptions of the Theorem. 

Whether the integrality and finiteness of Gopakumar--Vafa invariants can be reproved directly using unramified Gromov--Witten invariants remains an open question. In particular, integrality is somewhat surprising, since the moduli spaces of unramified maps are Deligne--Mumford stacks. By contrast, finiteness appears more natural. In dimension one, the moduli spaces of admissible covers with prescribed ramification are empty for genera that do not satisfy the Riemann--Hurwitz formula. In higher dimensions, however, moduli spaces of unramified maps need not be empty for sufficiently large genus, since the unramified condition becomes less restrictive.\footnote{For example, consider hyperelliptic curves mapping with degree $2$ to a line in $\p^2$. Such maps have simple ramification points. By attaching rational components mapping to smooth conics  that meet the divisor at infinity with multiplicity $2$ at these ramification points, one obtains an unramified map. Since hyperelliptic curves exist in every genus, this construction produces unramified maps of degree $2$ to $\p^2$ of arbitrarily large genus.}
\subsection{Ciocan-Fontanine--Kim--Zhou formula }The formula from the Theorem is a special instance of a much more general formula from Theorem \ref{wallcrossing}, which relates both relative and absolute descendant invariants of two theories in any dimension. This more general formula is in turn  a member of a vast family of wall-crossing formulas that hold in many different contexts: quasimap theory, gauged linear sigma models,  Donaldson--Thomas theory, Hurwitz theory, etc.  All of those wall-crossing formulas have exactly the same symbolic form and differ only semantically.  The first wall-crossing formula of such kind was conjectured by Ciocan--Fontanine and Kim in the context of GIT quasimaps \cite{ CFK14, CFKmi}, then Zhou proved their conjecture in full generality \cite{YZ}. We believe that the ideas of Zhou can be used to prove all such wall-crossing formulas after some contextual modifications. 

What unites all of these formulas? In many situations, there are two ways to compactify a moduli space: we either introduce more degenerate \textit{objects} into our moduli problem, or we allow the \textit{ambient space} over which these objects are defined to become more degenerate. For example, if we want to compactify the space of maps from smooth curves to a smooth one-dimensional target, we can either follow the Gromov--Witten approach by considering arbitrarily degenerate maps from nodal curves, or the Hurwitz approach -- that is, one-dimensional unramified Gromov--Witten theory -- by considering admissible covers of isotrivial semistable degenerations (bubblings) of the target. In this case, the ambient space is the one-dimensional target, while the objects over it are maps from curves. On the other hand, in the theory of quasimaps, we trade degenerations of the source curve (rather than the target) for degenerations of the maps. Here, the source curve plays the role of the ambient space, since the maps are regarded as objects defined over it.

 Whenever there is a choice between making the objects or the ambient space more degenerate, there should exist a wall-crossing formula, which we call  the \textit{Ciocan–Fontanine–Kim–Zhou formula}, relating the associated invariants. It should be proved using hybrid moduli spaces that contain both degenerate objects and degenerate ambient spaces, together with Zhou's master spaces arising from the theory of entangled degenerations (or entangled tails, in Zhou's original terminology). In the notation of Gromov–Witten theory, the Ciocan–Fontanine–Kim–Zhou formula is expressed as follows,\footnote{In some extremal cases, the formula takes a different form and involves the meromorphic part of $I$-functions; see, for example, \cite[Theorem 1.12.1]{YZ}.} 
\begin{multline*}
	\langle \gamma_1 \psi^{k_1}_1 \cdot \!\cdot \!\cdot \gamma_n \psi^{k_n}_n \mid \gamma'_1 \Psi^{k'_1}_1\cdot \!\cdot \!\cdot \gamma'_m \Psi^{k'_m}_m   \rangle^{- }_{B}	-	\langle \gamma_1 \psi^{k_1}_1\cdot \!\cdot \!\cdot \gamma_n \psi^{k_n}_n \mid \gamma'_1 \Psi^{k'_1}_1\cdot \!\cdot \!\cdot \gamma'_m \Psi^{k'_m}_m   \rangle^{+ }_{B}\\
=\sum_{\vec{B}} \Big\langle \prod_{j\in B_0} \gamma_j \psi^{k_j}_j \ \Big| \ \gamma'_1 \Psi^{k'_1}_1\cdot \!\cdot \!\cdot \gamma'_m \Psi^{k'_m}_m \cdot  \prod^{k}_{i=1} \mu_{B_i}\Big( -\Psi_{m+i}, \prod_{j\in B_i} \gamma_j \psi^{k_j}_j \Big)  \Big\rangle^{+}_{B_0} \big/ k!,
\end{multline*}
such that:
\begin{itemize} 
	\item $\langle \hdots \mid \hdots \rangle_B^-$ and $\langle  \hdots \mid \hdots  \rangle_B^+$ denote integrals over moduli spaces with degenerate objects and degenerate ambient spaces, respectively;
	\item  the insertions to the left of the vertical bar are absolute insertions, while those to the right are relative insertions;
	\item the distinction between absolute and relative insertions is determined by the choice of ambient space: relative insertions are associated with the ambient space;
	\item  the subscript $B$ denotes the  discrete topological  data of the theory, such as the genera of curves, degrees of maps,  Chern classes, etc.,  together with the number of absolute insertions; $j \in \{1,\dots ,n\}$ is the index of absolute insertions;  and the sum is taken over all  effective partitions $\vec{B}=(B_0,B_1, \dots, B_k)$ of $B$;
	\item $\mu_{B_i}(z , \prod_{j\in B_i} \gamma_j \psi^{k_j}_j  )$ is a polynomial in the variable $z$ with coefficients in some appropriate cohomology group, constructed as a truncation of an $I$-function (also known as a Vertex function) with insertions specified by $\prod_{j\in B_i} \gamma_j \psi^{k_j}_j$; it represents the class of a moduli space parametrizing degeneracies of the objects. 
	\end{itemize}

We refer to Theorem~\ref{wallcrossing} for an explanation of the above notation in the setting of unramified Gromov--Witten theory. The Ciocan-Fontanine--Kim--Zhou formula appears in numerous contexts in enumerative geometry:
\begin{itemize}
	\item it coincides with the mirror transformation of Yukawa couplings in the context of GIT quasimaps to the quintic threefold \cite{CFK14,CFKmi};
	\item it allows to prove the Landau--Ginzburg/Calabi--Yau correspondence  \cite{RR,CJR,ZLG,YZ};
	\item it establishes a correspondence between the Donaldson--Thomas theory of $S\times C$ and the Gromov--Witten theory of moduli spaces of sheaves on a surface $S$ \cite{N,NK3, NHiggs1,NHiggs2};
	\item it establishes a correspondence between the Gromov--Witten theory of $X\times C$ and the Gromov--Witten theory of symmetric products of a variety $X$ \cite{N22};
	\item it recovers the Gromov--Witten/Hurwitz correspondence \cite{OPcom}, as shown in \cite{NSc};
	\item it relates standard Gromov--Witten and cuspidal Gromov--Witten invariants, recovering the Li--Zinger formula \cite{LiZi}, as shown in \cite{NB};
	\item it computes zero-dimensional Donaldson--Thomas invariants on threefolds \cite{NHilb}.
\end{itemize}

\subsection{Organisation of the paper} In Section \ref{secun}, we set the stage by introducing Fulton--MacPherson degenerations \cite{FM} and unramified Gromov--Witten theory \cite{KKO}. Next, in Section \ref{sectionepsilon}, we define $\epsilon$-unramification, and prove basic properties of moduli spaces of $\epsilon$-unramified maps. In Section \ref{sectionrelative}, moduli spaces of $\epsilon$-unramified maps relative to moving exceptional divisors are constructed. They generalise moduli spaces of  $\epsilon$-unramified maps, and appear  in the localisation of master spaces. 

In Section \ref{sectioninvariants}, the obstruction theory of the moduli spaces of relative $\epsilon$-unramified maps is constructed using logarithmic geometry, and the associated invariants are defined. In Section \ref{masterspaces},  we introduce master spaces which compare invariants associated to values of $\epsilon$ in adjacent chambers. The corresponding wall-crossing formulas are then proved in Section \ref{sectionwallcrossing}.  The main theoretical input is Zhou's entanglement for moduli spaces of weighted Fulton--MacPherson degenerations; the latter are discussed in Section \ref{sectionweights}.  Finally, in Section \ref{sectionthreefolds}, we evaluate the wall-crossing formulas for threefolds, using divisor and dilaton equations in unramified Gromov--Witten theory that are established in Section \ref{sectiondil}. 
\subsection{Acknowledgments}  I am grateful to  Dan Abramovich, Luca Battistella, Siao Chi Mok, Samouil Molcho, Georg Oberdieck, Rahul Pandharipande, Dhruv Ranganathan, David Rydh, and Maximilian Schimpf
 for useful discussions on related topics.  I also thank the anonymous referee for carefully reading the manuscript and for providing numerous helpful comments, which contributed significantly to improving the exposition.

This work was a part of a project that has received funding from the European Research Council (ERC) under
the European Union’s Horizon 2020 research and innovation programme (grant agreement No. 101001159). It was also supported by a Hermann-Weyl-instructorship at the Forschungsinstitut f\"ur Mathematik at ETH Z\"urich. 

\subsection{Conventions} We work over the field of complex numbers $\BC$. All curves are assumed to be projective unless stated otherwise. Cohomology and homology groups are taken with rational coefficients.

	\section{Moduli spaces of unramified maps} \label{secun}
	
	\subsection{Fulton--MacPherson degenerations} \label{FMdeg} Let $X$ be a smooth projective variety over the field of complex numbers $\BC$.  Consider a Fulton--MacPherson compactification of the configuration space of $n$ points on $X$, denoted by $X[n]$, and constructed in \cite{FM}. If $\dim(X)=1$, then 
	\[
	X[n]=\Mbar_{X,n},
	\]
	 where $\Mbar_{X,n}$ is the moduli space of stable marked curves with the fixed root component $X$. 
	 
	  Let 
	\[\pi_{X[n]}\colon \FW[n] \rightarrow X[n] \quad \text{and} \quad p_{X[n]} \colon \FW[n] \rightarrow X\]
	  be  the universal family and  the universal contraction, respectively, also constructed in \cite{FM}. A pair consisting of a variety and a map, 
	   \[ (W, p \colon W \rightarrow X),\] is a \textit{Fulton--MacPherson (FM) degeneration}  of $X$ if it is isomorphic to a  fiber of $\pi_{X[n]}$ over a closed point in $X[n]$ for some $n$.  Alternatively, 
	$W$ can be constructed as a central fiber of an iterated  blow-up of $X\times \BC$, such that the centres of blow-ups are points in the regular locus of the central fiber of the previous blow-up.
	In short, $W$ is an iterated degeneration to the normal cone of a regular point.
	
To each FM degeneration one can associate a tree graph whose vertices correspond to its irreducible components, and an edge connects two vertices if the corresponding irreducible components touch. An example  of a tree associated to a FM degeneration is depicted in Figure \ref{tree1}.  The \textit{root vertex} is the vertex corresponding to a blow-up of $X$. We will also refer to the blow-up of $X$ as a \textit{root component}. An \textit{end component} of $ W$ is a component which corresponds to a vertex with one edge in the associated tree graph, excluding the root vertex; it is isomorphic to $\p(T_xX \oplus \BC)$ up to automorphisms fixing the divisor at infinity $\p(T_xX)$. A \textit{ruled component} of $ W$ is a component which corresponds to a vertex with two edges, excluding the root vertex; it is  isomorphic to $\mathrm{Bl}_0(\p(T_xX \oplus \BC))$ up to scaling automorphisms  which fix the fibration structure given by the projection away from the origin  $0 \in \p(T_xX \oplus \BC)$.  More generally, all non-root components are blow-ups of $\p(T_xX \oplus \BC)$ at finitely many distinct points. If $\dim(X)=1$, then a FM degeneration is a more familiar \textit{bubbling} of a curve; end components are \textit{rational tails}, while ruled components are \textit{rational bridges}. Moreover, if $\dim(X)=1$, the moduli spaces of unramified maps, defined in Section \ref{unr}, coincide with moduli spaces of admissible covers of Harris and Mumford \cite{HM}. 
	
We call a union of components of $W$ whose intersection graph is a connected graph with all vertices of valency less or equal to 2 a \textit{chain} of components. More generally, a union of components whose intersection graph is a connected tree will be referred to as a \textit{tree} of components.

	We will refer to fibers of the fibration associated to $\mathrm{Bl}_0(\p(T_xX \oplus \BC))$  as \textit{fiber lines}, and to proper transforms of lines not passing through $0$ in $\p(T_xX \oplus \BC)$ as \textit{non-fiber lines}.

	The \textit{group of automorphisms} of a FM degeneration $(W,p \colon W \rightarrow X)$ is the group of automorphisms of the variety $ W$  fixing the projection $p$, which we will denote simply by $\Aut( W)$. For brevity, we will also refer to automorphisms of a FM degeneration  as automorphisms of $W$.

	\begin{figure}[!ht]
		\centering
		\begin{tikzpicture}		
			\node (0) at (0.25,1) {};
			\node (1) at (-2,2) {};
			\node (2) at (-0.5,2) {};
			\node (3) at (1,2) {};
			\node (4) at (2.5,2) {};
			\node (6) at (0.25,3) {};
			\node (7) at (1.75,3) {};
			
			\draw[black] (0.25,3)--(1.75,3);
			\draw[black] (0.25,1)--(-0.5,2);
			\draw[black] (-2,2)--(-0.5,2);
			\draw[black] (-0.5,2)--(0.25,3);
			\draw[black] (-0.5,2)--(1,2);
			\draw[black] (1,2)--(2.5,2);


			\filldraw[thick, fill = black] (0.25,1) circle (.1cm) node at (0) {};
			\filldraw[thick, fill = black] (-2,2) circle (.1cm) node at (1) {};
			\filldraw[thick, fill = black] (-0.5,2) circle (.1cm) node at (2) {};
			\filldraw[thick, fill = black] (1,2) circle (.1cm) node at (3) {};
			\filldraw[thick, fill = black] (2.5,2) circle (.1cm) node at (4) {};
			\filldraw[thick, fill = black] (0.25,3) circle (.1cm) node at (6) {};
			\filldraw[thick, fill = black] (1.75,3) circle (.1cm) node at (7) {};

			\node at (2, 3.4) {$\mathrm{end}$};
			\node at (0.65, 3.4) {$\mathrm{ruled}$};
			\node at (-1.7, 2.4) {$\mathrm{root}$};
			\node at (1.4, 2.4) {$\mathrm{ruled}$};
			\node at (2.75, 2.4) {$\mathrm{end}$};
			\node at (0.5, 0.7) {$\mathrm{end}$};

		\end{tikzpicture}
		
		\caption{Tree.} \label{tree1}
	\end{figure}

\subsection{Moduli spaces of Fulton--MacPherson degenerations}In \cite{FM}, Fulton and MacPherson constructed the space $X[n]$ as an iterated blow-up of the $n$-fold product $X^n$, but not as a moduli space. In fact, in the literature there is no construction of $X[n]$ as a moduli space of stable marked FM degenerations. Nevertheless, it is  possible to construct the algebraic stack of FM degenerations  as a groupoid, using the blow-up description of $X[n]$ and the associated forgetful morphisms. The resulting stack is sufficient for all the purposes pursued in this work. A minor drawback is that the notion of families of FM degenerations requires an \'etale-local identification with $\FW[n]$ for some $n$; see Remark~\ref{remarkFM} for further details.

\begin{defn} \label{family} A family of FM degenerations over a base scheme $B$ is a pair, consisting of an algebraic space $\CW$ and a projection $p_B$, 
\[ ( \pi_B \colon  \CW \rightarrow B, \quad p_{B} \colon  \CW \rightarrow X).\]
We require that there exists an \'etale surjection  $T \rightarrow B$ and a map $T \rightarrow X[n]$ for some $n$, such that $\pi_{B}|_{T}$ and $p_{B}|_{T}$  are pullbacks of $\pi_{X[n]}$ and $ p_{X[n]}$ via $T \rightarrow X[n]$, respectively. We will denote a family of FM degenerations simply by $\CW$. 
\end{defn}

An algebraic moduli stack of families of FM degenerations in the sense above is constructed in \cite[Section 2.8]{KKO} as a groupoid
	\[ 
	\FM_{X}:=[X[n:n] \rightrightarrows X[n]],
	\]
	where $X[n:n]$ is the Fulton--MacPherson compactification with two sets of markings, such that FM degenerations are stable with respect to both sets; the morphisms are given by forgetting either of the two sets. Geometrically, by its definition, this groupoid identifies stable marked FM degenerations with different markings but the same underling FM degeneration, i.e., it is exactly the stack of FM degenerations, such that the integer $n$ corresponds to the maximal number of irreducible components; we choose $n$ to be sufficiently large without including it in the notation of $\FM_{X}$. By the construction, $X[n]$ is a smooth atlas of $\FM_{X}$. 
\begin{rmk} \label{remarkFM}
	Note that if $X[n]$ were shown to be the moduli space of stable marked FM degenerations, then any flat morphism
	\[
	\pi_B \colon \CW \to B,
	\]
	whose geometric fibers are FM degenerations would, \'etale locally on the base $B$, admit a map to $X[n]$ such that the pullback of $\FW[n]$ is isomorphic to $\CW$. Indeed, \'etale locally we can always choose sections of $\pi_B$ that stabilize the fibers of $\CW$ and therefore induce a map from $B$ to $X[n]$ for some $n$.
	
	Hence, if $X[n]$ were the moduli space of stable marked FM degenerations, the conditions of Definition~\ref{family} would automatically be satisfied.  It should not be difficult to show that $X[n]$ is indeed a moduli space; however, we will not pursue this here and will instead follow the approach of \cite{KKO}.
	\end{rmk}

\subsection{Blow-ups and contractions} \label{blowupscont}
In what follows, we will extensively use \textit{blow-ups} and \textit{contractions} of families of FM degenerations over spectrums of discrete valuation rings. In general, this is a subject of semistable reduction theorems, first established for curves by Deligne and Mumford \cite{DM}. In the case of FM degenerations of a fixed smooth variety $X$, it takes a particularly simple form: we can use  the forgetful morphisms, 
\[
X[n] \rightarrow X[n-1],
\]
constructed in \cite{FM} and revisited in  \cite[Section 2.6]{KKO}, to blow up and contract irreducible components of FM degenerations. Since  Definition \ref{family} requires an \'etale-local identification with $\FW[n]$ for some $n$, we have to make sure to retain this property after blow-ups and contractions; forgetful morphisms guarantee that this property is indeed preserved. 

\subsubsection{Contractions} Let $R$ be a discrete valuation ring with the fraction field $K$, and let 
\[
\CW \rightarrow \Spec(R)
\] be a family of FM degenerations in the sense of Definition \ref{family}. Let $\CW^*$ and $W$ be its generic and central fibers, respectively. By Definition \ref{family}, possibly after passing to a finite \'etale cover,\footnote{This is insignificant because we will only need blow-ups and contractions to prove the properness of moduli spaces, in which case one is allowed to pass to finite extensions. Alternatively, since our moduli spaces are locally Noetherian and defined over $\BC$, it suffices to consider $R = \BC[\![t]\!]$, which does not have any nontrivial \'etale covers.}  there exists  a map 
\[\Spec(R)\rightarrow X[n],
\]
 such that $\CW$ is the pullback of the universal family $\FW[n]$. There also exist  auxiliary markings \[\sigma_i \colon \Spec(R) \rightarrow  \CW, \quad i=1,\dots ,n,\]
 induced by the markings of the universal family $\FW[n]$. Assume that
 \[
 P \subset  W\subset \CW
 \] 
 is an end component which is not a limit of an end component in the generic fiber, that is, there exists no end component $\CP^* \subset \CW^*$ at a generic fiber, such that the central fiber of the closure of $\CP^*$ is $P$. Then, we can contract $P$ by inductively applying the forgetful morphisms $X[n] \rightarrow X[n-1]$ to all markings $\sigma_i$ that stabilize\footnote{By definition, markings stabilise a component if the automorphism group of the component with the markings is finite, and it is infinite without them; the root component is assumed to be stable with or without markings. } $P$. Let  
 \[ \Spec(R) \rightarrow X[n'] \] 
 be the  map given by forgetting all markings that stabilize $P$, and let  $\mathrm{Cont}_P(\CW)$ be the pull-back of the universal family $\FW[n']$. Using the forgetful morphism between the universal families $\FW[n] \rightarrow \FW[n']$, we obtain the  contraction of $P$,
 \[ \CW \rightarrow \mathrm{Cont}_P(\CW).\]
  Then, $\mathrm{Cont}_P(\CW)$ is a family of FM degenerations in the sense of Definition \ref{family}.   The same procedure can be  applied to contract a ruled component $P_{\mathrm{rul}}\subset W\subset \CW$.
  
     \subsubsection{Blow-ups}
  To blow up a point 
  \[
  x \in W^{\mathrm{sm}}\subset  \CW
  \]
   in the regular locus of the central fiber, we introduce two additional auxiliary markings $(\sigma'_1, \sigma'_2)$  on the family $\CW$ intersecting at $x$ with distinct tangential directions.  Since the family $\CW$ is smooth around $x$, this is always possible. We may assume that they do not intersect markings $\sigma_i$ by taking $\sigma'_j$ to be one of $\sigma_i$ if this is the case, and working with total of $n+1$ markings.    By using the properness of the forgetful morphisms $X[n+2] \rightarrow X[n]$, we can lift $\Spec(R)$ to $X[n+2]$ and pull back the universal family $\FW[n+2]$, to obtain a blow-up \[\mathrm{Bl}_x  (\CW) \rightarrow  \CW,\] 
  such that $\mathrm{Bl}_x  (\CW)$ is a family of  FM degenerations in the sense of Definition \ref{family}. 
  
  The same procedure applies to blow-ups of $\CW$ at a singular locus of the central fiber, i.e., a simple normal crossing singularity at the intersection of two components. We outline how to blow up the isolated simple normal crossing singularity 
  \[
  D \subset W \subset \CW,
  \] leaving the case of a non-isolated singularity to the reader.  For simplicity, assume $R=\BC[\![t]\!]$, then by \cite[Section 2.3]{KKO}, in the formal neighbourhood around points of the singular locus $D$,  the family  $\CW$ is isomorphic to
  \begin{equation} \label{an0}
  	\mathrm{Spf} \left(\frac{\BC[\![z_1,\dots ,z_{\dim+1}, t ]\!]}{(z_1z_2-t^{b})} \right),
  \end{equation} 
which is the $A_{b-1}$-singularity.
If $b=1$, we first take a base change with respect to $t \rightarrow t^2$ and then introduce an auxiliary marking $\sigma$ which intersects $D$ orthogonally. The base change is necessary because there exists no single-valued section otherwise, which is reflected in the fact that the blow-up at $D$ for $b=1$ produces a non-reduced fiber at $t=0$. For $b=2$, such section is given by $t \mapsto (t,t, 0, \hdots, 0)$. Then, as in the case of a regular point, we can blow up the singular locus,
 using the forgetful morphisms,
 \[
 \mathrm{Bl}_D  (\CW) \rightarrow  \CW,
 \] 
 this  preserves the notion of a family from Definition \ref{family}, and results in the creation of one ruled component, resolving the $A_1$-singularity. More generally, for $b\ge 2$, by introducing auxiliary markings, we can inductively blow up singular loci $D$, thereby (partially) resolving the $A_{b-1}$-singularity and producing chains of at most $b-1$ ruled components; in this case, the base change is not necessary.

 \subsubsection{Dominating family} Finally, let   \[ \CW_{i} \rightarrow \Spec(R), \quad i=1,2,\]
  be two families of FM degenerations with isomorphic generic fibers 
  \[
   \CW^*_1\cong  \CW^*_{2}.
   \] We add both auxiliary collections of markings $\sigma_{1,i}$ and $\sigma_{2,i}$ to each family $ \CW_{i}$. By removing repeating markings, we may assume that all of them are distinct. By the properness of the forgetful morphisms $X[n_1+n_2] \rightarrow X[n_i]$, we obtain a family $\CW \rightarrow \Spec(R)$ of FM degenerations in the sense of Definition \ref{family} which dominates both families
  \[ \widetilde{ \CW}\rightarrow \CW_{i}, \quad i=1,2, \]
  and is given by an iterated blow-up.
\subsection{Unramified maps}\label{unr}	\begin{defn}[\hspace{-0.01cm}{\cite{KKO}}] \label{stab}Let $(C,\underline{p}):=(C,p_1,\dots,p_n)$ be a possibly disconnected marked nodal curve, and let $W$ be a FM degeneration of $X$. By $C^{\mathrm{sing}}$ and $C^{\mathrm{sm}}$ we denote the singular and smooth loci of $C$, respectively; the same notation applies to $W$.  For a map
\[
f\colon (C, \underline{p}) \rightarrow  W,
\]
 let $\Omega_{f}$ denote its relative cotangent sheaf.  The map is said to be  \textit{unramified} if 
		\begin{itemize} 
			\item[(1)] $\Omega_{f}|_{C^{\mathrm{sm}}}=0$,
			\item[(2)] $f^{-1}( W^{\mathrm{sing}})=C^{\mathrm{sing}}$,
			\item[(3)] images of marked points $f(p_i)$ are pairwise distinct, 
			\item[(4)] the map $f$ is \textit{admissible} over $ W^{\mathrm{sing}}$, that is,\ in a formal neighbourhood of a node $p \in C$, the map $f$ is of the following form, 
			\begin{align*} f\colon \mathrm{Spf}(\BC[\![x_1,x_2]\!]/(x_1x_2)) &\rightarrow   \mathrm{Spf}(\BC[\![z_1,\dots, z_{\dim(X)+1}]\!]/(z_1z_2)), \\
		f^*z_1=x_1^\ell, &\quad  f^*z_2=x_2^\ell,
				\end{align*}
				\item[(5)] every ruled component $P_{\mathrm{rul}} \subset  W$ contains an image of a marking or a non-fiber-line image  of an irreducible component of $C$,
				\item[(6)] every end component $P \subset  W$ contains images of two distinct markings  or a non-line image of an irreducible component of $C$.

			\end{itemize}
We say $f$ is of degree $\beta \in H_2(X,\BZ)$ if $(p\circ f)_*[C]=\beta$.  The genus of a  disconnected curve 
\[
C=\cup^\ell_{i=1} C_i
\] is defined to be its arithmetic genus, 
\[g(C):= \sum^\ell_{i=1} g(C_i)-\ell+1=1-\chi(\CO_C),\] 
where $C_i$ are connected components of $C$. 	
		\end{defn}	

  An automorphism of a map $f\colon (C,\underline{p}) \rightarrow  W$ is an automorphism of $C$ together with an automorphism of $ W$ which fix $f$, 
 \[ \Aut_{C, W}(f):=\{(g_C, g_ W)\in \Aut(C,\underline{p})\times \Aut( W) \mid f\circ g_C=g_ W\circ f\}.\]
 Later, we will also consider the automorphism group of $f$ with respect to automorphisms of the source,  
 \[
 \Aut_{C}(f):=\{g_C\in \Aut(C,\underline{p}) \mid f\circ g_C= f\}.
 \] 
 
A \textit{family} of unramified maps over a base scheme $B$ is given by a map over $B$, 
 
\[f \colon (\CC, \underline{p} ) \rightarrow  \CW,\]
such that $(\CC,\underline{p})$ and $\CW$ are families of marked nodal curves and FM degenerations over $B$,  geometric fibers of $f$ over $B$ are unramified maps, and $f$ is admissible in formal neighbourhoods of  geometric points of $B$; see \cite[Definition 3.2.1]{KKO} for the last condition.  We will denote a family of unramified maps as 
\[
 (\CC,  \CW, \underline{p}, f). 
\]
   Isomorphisms of families are considered with respect to both the target and the source, that is, two families are isomorphic 
\[
(\CC,  \CW, \underline{p}, f)\cong  (\CC',  \CW', \underline{p}', f')
\]
if 
\[ 
f \circ g_\CW=g_\CC \circ f,
\]
where $g_\CC \in \mathrm{Iso}_B((\CC,\underline{p}),(\CC', \underline{p}'))$ and $g_\CW \in \mathrm{Iso}_B(\CW,\CW')$ are $B$-relative isomorphisms of families of curves and FM degenerations. 
\begin{defn}Let 
	\[\Mbar^+_{g,n}(X,\beta) \colon (Sch/\BC)^\circ \rightarrow Grpd\] 
	be the moduli space of unramified maps of degree $\beta$ from possibly disconnected curves of genus $g$ with $n$ marked points, where $Grpd$ denotes the 2-category of groupoids. 
	\end{defn}
	\begin{thm}[\hspace{-0.01cm}\cite{KKO}] The moduli space  $\Mbar^+_{g,n}(X,\beta)$ is a proper Deligne--Mumford stack with a perfect obstruction theory. 
		\end{thm}
	\section{Moduli spaces of $\epsilon$-unramified maps} \label{sectionepsilon}
	\subsection{Weights} \label{weightsdef} In this section, we will  define a stability condition that interpolates between unramified maps and stable maps. 
	For a map $f\colon (C,\underline{p}) \rightarrow  W$,  we will assign a weight $w(x) \in \frac{1}{2}\BZ$ to each point $x \in  W^{\mathrm{sm}}$, based on how ramified $f$ is over $x$. To each end component $P \subset  W$, we assign a weight $w(P) \in \frac{1}{2}\BZ$ based on how big the curve $C$ is over $P$.  This is done as follows.

			\begin{defn} \label{weightp}  Given a map  $f\colon (C,\underline{p}) \rightarrow  W$ and a point $x\in  W^{\mathrm{sm}}$, we introduce the following notation: 
				\begin{itemize}
			\item $C_\bullet \subset  C$ is the maximal subcurve of $C$ contracted by $f$ to $x$,
			\item $C_\circ \subset C$ is the closure of the complement of $C_\bullet$,
			\item  $f_{\circ}$ is the restriction of $f$ to $C_\circ$,
			\item  $N \subset C$ is the set of nodes separating $C_\bullet$ from $C_\circ$,
			\item  $\Omega_{f_\circ}$ is the relative cotangent sheaf of $f_\circ$,
			\item $\etalt_p(\hdots )$ is the multiplicity of a $0$-dimensional sheaf at a point $p\in C$,
			\item $n_x$ is the number of marked points in the fiber of $f$ over $x$.
				\end{itemize}
				 We then define the weight of a point $x$, 
			\[w(x):=-2\chi(\CO_{C_\bullet})+ \sum_{p \in f_\circ^{-1}(x)} \etalt_p(\Omega_{f_\circ}) +2|N| + n_x/2, \]
		which is a half-integer $w(x) \in \frac{1}{2}\BZ$. A point $x \in  W^{\mathrm{sm}}$ is called a \textit{base point} if $w(x)>0$. 				
	\end{defn}

\begin{defn} \label{weightend} Given a map  $f\colon (C,\underline{p}) \rightarrow  W$ and an end component $P \subset  W$, we introduce the following notation:
	\begin{itemize}
		\item $f|_{P} \colon (C,\underline{p})|_{P} \rightarrow P$ is the restriction of $f$ to $P$,
		\item   $D:=\p(T_xX)$, $D\subset P=\p(T_xX \oplus \BC)$ is the divisor at infinity,
		\item  $N_0 \subset C$ is the set of nodes mapping to $D$,
		\item  $n_P$ is the number of marked points on  $(C,\underline{p})|_{P}$.
	\end{itemize}
		 We then define the weight of an end component $P$,
		\[w(P):=-2\chi(\CO_{C}|_{P})+\deg_{D}(f)+|N_0| +n_P/2, \] 
		which is a half-integer $w(P) \in \frac{1}{2}\BZ$.  
		\end{defn}

\subsection{Examples of weights} Let us give a few examples of weights associated to points and end components. First, consider a genus $1$ curve $E$ mapping with degree 2 onto a projective line $X=\p^1$, 
\[
f \colon E \rightarrow X.
\]
By the Riemann--Hurwitz formula, $f$ must have 4 ramification points of multiplicity $1$, that is,  $\Omega_f$ is of multiplicity $1$ at these points. Since there are no contracted components and markings, we obtain that 
\[
w(x)=0+1+0+0=1
\]
if $x\in X$ is a branching point (i.e., an image of a ramification point), and $w(x)=0$ otherwise. 

Consider now a nodal genus $1$ curve  $E\cup \p^1$ mapping to $X=\p^1$,
\[
f\colon E\cup \p^1 \rightarrow X,  
\]
 such that the source $\p^1$ is mapped onto the target $X$ with degree 2, and $E$ is contracted to a point in $X$. A degree $2$ self-map of $\p^1$ must have $2$ ramification points. Assume $E$ is attached away from those ramification points, then 
\[ 
w(x)=0+0+2+0=2
\]
if $x \in X$ is the image of the contraction of $E$,  $w(x)=1$ if $x$ is a branching point, and $w(x)=0$ otherwise. If $E$ is attached to a ramification point, then 
\[ 
w(x)=0+1+2+0=3,
\]
if $x \in \p^1$ is the image of the contraction of $E$. 

All maps discussed above appear in the moduli spaces of degree 2 stable maps to $\p^1$, \cite{Kon}. On the other hand, the moduli space of admissible covers  of $\p^1$, \cite{HM}, prohibits contracted components. Instead, it allows FM degenerations of $\p^1$. Consider a degree 2 admissible cover of a nodal genus 0 curve $W=P\cup X$ by a  nodal genus 1 curve $E \cup \p^1$, 
\[ 
f \colon E\cup \p^1 \rightarrow W,
\]
such that $E$ is mapped onto the target $P=\p^1 \subset W$ and the source $\p^1$ is mapped onto $X$. 
Since $E$ is attached to $\p^1$ via a single node, and $f$ is of degree $2$, it must be that this node is a ramification point of both $f|_{E}$ and $f|_{\p^1}$. We therefore obtain that 
\[ 
w(P)=0+2+1+0=3
\]
for the end component $P \subset W$. Observe the equality
\[ 
w(P)=w(x)=3
\]
for the image point $x \in X$  of the contraction of $E$ which is attached to a simple ramification, discussed in the previous paragraph. This equality is not a coincidence. In fact, the weight functions are tailored so that they satisfy such property, Lemma \ref{contraction}. Indeed, the end component $P \subset W$ and a point $x \in X$ as above are related by the contraction 
\[
\pi \colon W \rightarrow X
\] of the end component $P \subset W$ and the composition of the map $f$ with the contraction. 

\begin{rmk} If $\dim(X)=1$, our weight functions coincide with the weights given by multiplicities of  the branching divisor in \cite{N22}. The branching divisor $\mathrm{br}(f)$ is associated to an admissible map $f \colon C \rightarrow W$ between two nodal curves, and is  given by the cohomologically weighted support of the $0$-dimensional perfect complex
	\begin{equation} \label{branchingcomplex}
	Rf_*[f^*\Omega_C^{\mathrm{log}}\rightarrow \Omega_W^{\mathrm{log}}].
	\end{equation}
	Then,  the weight $w(x)$ from Definition \ref{weightp} is equal to the multiplicity of $\mathrm{br}(f)$ at the point $x$, 
	\[ 
	w(x)=\mathrm{mult}_x(\mathrm{br}(f)).
	\]
	 While $w(P)$ from Definition \ref{weightend} is equal to the degree of $\mathrm{br}(f)|_{P}$, 
	 \[ 
	 w(P)=\deg(\mathrm{br}(f)|_{P}). 
	 \]
	If $C$ has markings, the same conclusion holds up to the prefactor $1/2$. 
	
	In higher dimensions, the relative cotangent complex~\eqref{branchingcomplex} is not $0$-dimensional. We therefore define the weight functions so that they specialize to the branching-divisor weights when $\dim(X)=1$ while satisfying Lemma~\ref{weightP} and Lemma~\ref{contraction}. The latter lemma is essential for proving the separatedness of the moduli spaces in Theorem~\ref{ft}.
\end{rmk}
	\begin{rmk} Markings are weighted by $1/2$ to distinguish one marking from ramifications of order $1$.  Recall that in the unramified Gromov--Witten theory, at most one marking is allowed in a fiber of a map.  Any other positive rational number smaller than 1 also works. 
		
		 If the weight of markings is equal to 1, then the stability forbids marked points and ramifications of order 1 at the same time. If the weight of markings is less than 1, then the stability firstly forbids  ramifications of order 1, and then markings.  We refer to Section \ref{secret} for more details about the chamber where marked points are forbidden. 
	\end{rmk}

\subsection{Properties of weights} In this section, we establish the basic properties of weights. First, we  need to know that $w(P)$ behaves well in families; see also Section \ref{sectionweights}. Lemma \ref{weightP} provides an alternative construction of $w(P)$, which  also extends the weights to non-end components.

 Let us introduce the following notation associated to an admissible map 
\[
f\colon (C,\underline{p}) \rightarrow  W,
\]
and a non-root irreducible component  $W'$ of $W$: 
\begin{itemize}
	\item $\{D_1, \dots, D_k\} $ are the exceptional divisors of $W'$, and $D_0$ is the divisor at infinity, 
	\item   $N_j\subset C$ is the set of nodes of $C$ over $D_j$, and $N_0$ is the set of nodes over $D_0$,
	\item $\Omega^{\mathrm{log}}_{C}$ and $\Omega^{\mathrm{log}}_{ W}$ are the log cotangent bundles of $C$ (without the markings) and $ W$, respectively,
	\item $\dim(X)$ is abbreviated to $\dim$. 
	\end{itemize}
Then the following result expresses $w(P)$ in terms of degrees of logarithmic cotangent bundles. 
\begin{lemma} \label{weightP} Given an admissible map  $f\colon (C,\underline{p}) \rightarrow  W$, and  a non-root irreducible component $ W' \subset  W$, then
\begin{multline*} 
	\deg(\Omega^{\mathrm{log}}_{C}|_{W'})- \deg(f^*\Omega^{\mathrm{log}}_{ W}|_{W'})/\dim\\
	=-2\chi(\CO_{C}|_{W'})+ \sum^{k}_{j=0}|N_j|+\deg_{D_0}(f) -\sum^{k}_{j=1}\deg_{D_j}(f).
	\end{multline*}
In particular, for an end component $P \subset  W$, we have 
\[w(P)=\deg(\Omega^{\mathrm{log}}_{C}|_{P})- \deg(f^*\Omega^{\mathrm{log}}_{ W}|_{P})/\dim +n_P/2.\]
\end{lemma}
\textit{Proof.} First, 
\begin{equation*}
\Omega_{ W}^{\mathrm{log}}|_{W'} \cong \Omega_{ W'}\Big(\log\Big(\sum^{k}_{j=0} D_j\Big)\Big),\quad
\Omega^{\mathrm{log}}_{C}|_{W'}\cong \Omega^{\mathrm{log}}_{C|_{ W'}}\Big(\sum^{k}_{j=0}N_j\Big) \Big),
\end{equation*}
and there exists a short exact sequence relating the log cotangent bundle $\Omega_{ W'}(\log D)$ associated to a smooth divisor $D$ and the standard cotangent bundle, 
\begin{equation} \label{sequence}
0 \rightarrow \Omega_{ W'} \rightarrow \Omega_{ W'}(\log D) \rightarrow \CO_{D} \rightarrow 0,
\end{equation}
which can be seen directly through the definition of $\Omega_{ W'}(\log D)$ as the sheaf of $1$-forms with poles of order at most 1 along $D$. 
Moreover, if $D$ is the exceptional divisor of a blow-up at a point,
 \[\pi\colon  W' \rightarrow \widetilde{W}',\]
 then 
 \begin{equation} \label{blow}
 \Omega_{ W'}(\log D)\cong \pi^*\Omega_{\widetilde  W'}(D),
 \end{equation}
which is induced by the canonical morphism between cotangent sheaves, and follows from a calculation in the local charts of the blow-up.

In our case, $ W'$ is a blow-up of a projective space $\p^{\mathrm{dim}}$ in $k$ distinct points away from the divisor at infinity. Recall that 
\[
\mathrm{c}_1(\Omega_{\p^{\mathrm{dim}}})=-(\dim+1)[D_0].
\]  Applying the sequence (\ref{sequence}) to $D_0$ and the identification (\ref{blow}) to exceptional divisors $D_j$, we obtain the claim. 
  \qed  
\\

Next, we need to know the compatibility of weights associated to points and end components, which we establish in  Lemmas \ref{less} and \ref{contraction}. The key to these results is the following lemma. 

\begin{lemma} \label{mult} Consider a map $f\colon C \rightarrow \p^\di$ of degree $\deg_D(f)$. Let $x \in \p^\di$ be a point, such that the fiber $f^{-1}(x)$ contains only regular points of $C$, then 
	\[ \sum_{p \in f^{-1}(x)} (\etalt_p (\Omega_f)+1)\leq \deg_D(f). \]
In particular, we have 	$|f^{-1}(x)| \leq \deg_D(f).$ 
	\end{lemma}
\textit{Proof.} Let $E \subset \mathrm{Bl}_x (\p^\di)$ be the exceptional divisor of the blow-up of $\p^\di$ at $x$,
\[ \pi \colon \mathrm{Bl}_x (\p^\di) \rightarrow \p^\di,\]
 and let $\tilde{f} \colon C \rightarrow \mathrm{Bl}_x (\p^\di)$ be the unique lift of $f$ to $\mathrm{Bl}_x (\p^\di)$, which exists by our regularity assumption. The divisor $\pi^*D-E$ is nef on $\mathrm{Bl}_x (\p^\di)$, hence 
 \[\deg_E(\tilde{f}) \leq \deg_D(f). \]
On the other hand, we have  
\begin{equation} \label{multexc}
	\sum_{p \in f^{-1}(x)} (\etalt_p (\Omega_f)+1) =\deg_E(\tilde{f}), 
\end{equation}
this proves the claim. 
\qed 
\\

With Lemma \ref{mult}, we can now proceed to proving the compatibility properties of weights. 
\begin{lemma} \label{less} Consider an admissible map $f\colon (C,\underline{p}) \rightarrow  W$. Let $P\subset  W$ be an end component, then for all points $x \in P\setminus D$, we have 
	\[w(x)\leq w(P), \]
	such that equality holds, if and only if $x$ is the unique base point and the image of $f|_{P}$ is a union of lines.  
\end{lemma}
\textit{Proof.}  We will use the notation from Definition \ref{weightp}. First,
\[-2\chi(\CO_{C_\circ}|_{P})+|N|+|N_0|+K\geq 0 ,\]
where $K$ is the number of connected components of $C_{\circ}|_{P}$ which are isomorphic to $\p^1$ and do not share a node with $C_\bullet$, i.e.,\ $N\cap \p^1=\emptyset$. We also have 
\begin{equation} \label{genera}	-2\chi(\CO_{C_\bullet})-2\chi(\CO_{C_\circ}|_{P})+2|N|=-2\chi(\CO_{C}|_{P}),
	\end{equation}
which, by the inequality above, implies that  
\[ -2\chi(\CO_{C_\bullet})+|N|-K\leq-2\chi(\CO_{C}|_{P})+|N_0|. \]
 By Lemma \ref{mult}, we get 
\[ \sum_{p \in f_\circ^{-1}(x)} \etalt_p(\Omega_{f_\circ}) +|N| +K \leq  \deg_{D}(f),\]
where we also use that every connected component of $C_\circ|_{P}$ intersects $D$ (i.e., the image of $\p^1$ contributing to $K$ might not contain $x$, but it always intersects $D$). 
It is clear that $n_x\leq n_P$. Summing up these inequalities, we obtain 
\begin{multline*}
 w(x)=-2\chi(\CO_{C_\bullet})+ \sum_{p \in f_\circ^{-1}(x)} \etalt_p(\Omega_{f_\circ}) +2|N| + n_x/2 \\ 
 \leq -2\chi(\CO_{C}|_{P}) +\deg_{D}(f) +|N_0|+ n_P/2\leq w(P),
 \end{multline*}
this proves the first claim. 

 Assume for simplicity that $C_{\circ}|_P$ is connected. If $x$ is the unique base point and the image of $f_{\circ }|_P$ is a line, then by the Riemann--Hurwitz formula, the non-contracted component $C_{\circ}|_P$ is isomorphic to $\p^1$, and over its image the map $f_{\circ}|_P$  is given by  $f(z)=z^{\deg_{D}(f)}$, such that the branching points are $x$ and a point at $D$. On the other hand, the curve $C_\bullet$ contains all markings of $C|_{P}$ and is attached to the unique ramification point on $C_{\circ}|_P$ over $x$. In particular, $|N|=1 $ and $|N_0|=1$.   With this analysis we obtain that
\begin{multline*}
w(x)=-2\chi(\CO_{C_\bullet})+ \sum_{p \in f_\circ^{-1}(x)} \etalt_p(\Omega_{f_\circ}) +2|N| + n_x/2\\
=-2\chi(\CO_{C}|_{P})+\deg_{D}(f)-1+2+n_P/2=w(P).
\end{multline*}

Conversely, if $w(x)=w(P)$, then by the definition of weights and (\ref{genera}), we have 
\[
\sum_{p \in f_\circ^{-1}(x)} \etalt_p(\Omega_{f_\circ}) + n_x/2= -2\chi(\CO_{C_\circ}|_{P})+\deg_{D}(f)+|N_0|+n_P/2,\]
which by Lemma \ref{mult} and  $n_x\leq n_P$ implies that 
\[ 
-2\chi(\CO_{C_\circ}|_{P}) =-2, \quad |N|=1 , \quad |N_0|=1,
\] 
and
\[ 
\etalt_p (\Omega_{f_\circ})+1=\deg_D(f)
\]
for the unique ramification point $p$ of $f_{\circ}|_{P}$ over $x$. 
In particular, we conclude that 
\[
C_{\circ}|_{P} \cong \p^1.
\] Moreover,  since $\pi^*D-E$ is zero precisely on the fiber lines of $\mathrm{Bl}_x (\BP^{\dim})$, the inequality from Lemma \ref{mult} is satisfied, if and only if  the image of  $f_{\circ}|_P$ is a line.  Altogether this also implies that $f_{\circ}|_P$  is of the form  $f(z)=z^{\deg_{D}(f)}$, i.e., $x$ is the unique base point.  This finishes the proof in the case when $C_{\circ}|_P$ is connected. 

If $C_{\circ}|_P$ is disconnected, then the same analysis applies with the difference that $C_{\circ}|_P$ might be a disjoint union of $\p^1$'s. 
\qed

\begin{lemma} \label{contraction} Consider an admissible map
 $f\colon (C,\underline{p}) \rightarrow  W$.
	Let $\pi \colon  W \rightarrow  W'$ be the contraction of an end component of $ W$ and  $f':= \pi\circ f$ be the induced map.  
	Let $x \in  W'$ be the image of the contraction of $P \subset  W$, then the weight of $x$ with respect to $f'$ is equal to the weight of $P$ with respect to $f$, 
	\[w(x)=w(P).\]
\end{lemma}
\textit{Proof.}  Let $E \subset  W$ be the exceptional divisor at which $P$ is attached. By the admissibility of $f$, the ramification profile of $f$ over $D \subset P$ is equal to the ramification profile over $E \subset  W$, and $N=N_0$. Hence, by (\ref{multexc}), we obtain that
\begin{align*}\sum_{p \in f^{-1}(x)} \etalt_p (\Omega_f)+2|N|&=  \sum_{p \in f^{-1}(x)} (\etalt_p (\Omega_f)+1)+|N|  \\
	&=\deg_E(\tilde{f})+|N| =\deg_E (f)+|N_0| = \deg_{D}(f)+|N_0|.
	\end{align*}
Adding the contributions from the genus and the number of markings, we obtain the claim. \qed 

\begin{lemma} \label{ruledpos} Consider an admissible map $f\colon (C,\underline{p}) \rightarrow  W$.  For all ruled components $P_{\mathrm{rul}} \subset  W$ with an exceptional divisor $E$ and a divisor at infinity $D$, we have 
	\[ w(P_{\mathrm{rul}}):= -2\chi(\CO_{C}|_{P_{\mathrm{rul}}})+|N_0|+|N_1|+\deg_{D}(f)- \deg_E(f)+n_{P_{\mathrm{rul}}}/2 \geq 0. \]
	\end{lemma}
\textit{Proof.}  By Lemma \ref{mult}, we have
\[\deg_E(f) \leq \deg(f|_{P_{\mathrm{rul}}})=\deg_{D}(f),\]
we conclude 
\[w(P_{\mathrm{rul}})\geq -2\chi(\CO_{C}|_{P_{\mathrm{rul}}})+|N_0|+|N_1|+n_{P_{\mathrm{rul}}}/2.\]
Hence, all connected components of $C|_{P_{\mathrm{rul}}}$ of non-zero genera contribute non-negatively to $w(P_{\mathrm{rul}})$. Consider now a genus zero connected  component $C'$ of $C|_{P_{\mathrm{rul}}}$.  The set  $N_{0}|_{C'}$ is always non-empty, as all curves have to intersect $D$. If $N_{1}|_{C'}$ is also non-empty, then the claim follows from the inequality above.  On the other hand, if $N_1$ is empty, then $\deg_E(f)=0$. In both cases, the contribution of genus zero components to $w(P_{\mathrm{rul}})$ is therefore non-negative.   \qed

\subsection{Definition of $\epsilon$-unramified maps} \label{sectionepsilonun}
	\begin{defn} \label{epsilon} Choose a real number $\epsilon \in \BR_{>0}$. A map $f\colon (C,\underline{p}) \rightarrow  W$ is said to be $\epsilon$-\textit{unramified} if 
		\begin{itemize}
		\item[(1)] for all end components $P \subset  W$, $w(P)>  1/\epsilon$,
		\item[(2)] for all $x\in  W^{\mathrm{sm}}$, $w(x) \leq 1/\epsilon$,
		\item[(3)] the map $f$ is admissible over $ W^{\mathrm{sing}}$,
		\item[(4)] $|\Aut_{C,W}(f)|< \infty$ and all connected components map non-trivially. 
		\end{itemize}

		\end{defn}

Let us consider in more detail the $\epsilon$-unramification for some extremal values of $\epsilon$. First, if
\[
\epsilon >2,
\] then there cannot be any base points, since the minimal non-zero value of $w(x)$ is $1/2$ (which corresponds to one marking). This implies that  
\[
\Omega_{f}|_{C^{\mathrm{sm}}}=0 \quad \text{and} \quad f^{-1}( W^{\mathrm{sing}})=C^{\mathrm{sing}}.
\] On the other hand, end components of arbitrary weights are allowed. Moreover, Definition \ref{epsilon} (4) implies Definition \ref{unr} (5), while Definition \ref{epsilon} (1), (2), and Lemma \ref{less} imply Definition \ref{unr} (6). We conclude that for $\epsilon>2$, an $\epsilon$-unramified map is simply an unramified map without marked points; this will be important for the results in Section \ref{secret}. 

If 
\[
1<\epsilon \leq 2,
\] 
then the same analysis applies with the difference that markings with pairwise distinct images are allowed, which is exactly Definition \ref{unr} (3). Hence, an $\epsilon$-unramified map is an unramified map with markings in this case. 

Finally, if $\epsilon \ll 1$, then an $\epsilon$-unramified map is  a stable map to $X$. It is clear that all stable maps are $\epsilon$-unramified for a small enough $\epsilon\in \BR_{>0}$. However, to show that for a fixed degree $\beta$ and genus $g$, there exists a small enough $\epsilon$ such that all $\epsilon$-unramified maps of such degree and genus are stable maps, we need to know that the possible weights of end components $w(P)$ are bounded from above. This follows from the boundedness result for maps that are subject to $\epsilon$-unramification established in Lemma \ref{alg}.

Let us introduce a shortened notation for the extremal values of $\epsilon$. 
\begin{align*}
	\epsilon&=+ \quad \text{ if }1<\epsilon \leq 2, \\
	\epsilon&=-\quad \text{ if } \epsilon\ll 1. 
\end{align*}
By the discussion above, we have
\begin{align*}
	+\text{-unramified maps } &= \text{ unramified maps}, \\
	- \text{-unramified maps } &= \text{ stable maps}. 
\end{align*}

A family of $\epsilon$-unramified maps over a base scheme $B$ is defined as for unramified maps in Section \ref{unr}, requiring that the geometric fibers over $B$ satisfy the $\epsilon$-unramification instead of the unramification.
\begin{defn}
Let \[\Mbar^\epsilon_{g,n}(X,\beta) \colon (Sch/\BC)^\circ \rightarrow Grpd\] 
be the moduli space of $\epsilon$-unramified maps of degree $\beta$ from possibly disconnected curves of genus $g$ with $n$ markings.
\end{defn} 
\subsection{Moduli spaces of $\epsilon$-unramified maps}
There are two properties that we want $\epsilon$-unramification to satisfy:
\begin{itemize}
	\item[$\bullet$] openness in the stack of all admissible maps to FM degenerations,
	\item[$\bullet$] properness.
\end{itemize}
Openness is needed for the obstruction theory and the construction of virtual fundamental classes, while properness is required for integrating cohomology classes against the virtual fundamental classes; both are necessary for the definition of invariants in Section~\ref{sectioninvariants}. We will verify that these properties are indeed satisfied.

Let 
 \[\FM_{g,n}(\FW,\beta) \colon (Sch/\BC)^\circ \rightarrow Grpd\] 
be the moduli space of maps to FM degenerations $W$, subject to the following conditions:
\begin{itemize}
	\item for all end components $P \subset W$, $w(P)>0$,
	\item the map $f$ is admissible over $ W^{\mathrm{sing}}$,
	\item $|\Aut_{C,W}(f)|< \infty$ and all connected components map non-trivially. 
\end{itemize}	
A family of such maps is defined as in Section \ref{unr}.   It is clear that this space contains all geometric points of  $\Mbar^\epsilon_{g,n}(X,\beta)$,
 \[
 \Mbar^\epsilon_{g,n}(X,\beta)(K) \subseteq \FM_{g,n}(\FW,\beta)(K),
 \] for all $\epsilon \in \BR_{>0}$. Also, the defining conditions of $\FM_{g,n}(\FW,\beta)$ are open in the stack of all admissible maps. Hence, to establish the openness of $\Mbar^\epsilon_{g,n}(X,\beta)$ inside the stack of all admissible maps, it suffices to verify this inside $\FM_{g,n}(\FW,\beta)$. 
\begin{lemma} \label{alg} The moduli space $\FM_{g,n}(\FW,\beta)$ is a quasi-separated Deligne--Mumford stack of finite type. 
	\end{lemma}
\textit{Proof.} In \cite[Section 2.8]{KKO}, the stack of FM degenerations is shown to be algebraic.  By the argument in \cite[Proposition 3.3.2]{KKO}, the stack $\FM_{g,n}(\FW,\beta)$ is naturally a locally closed substack of the stack of all maps (not necessarily admissible) to FM degenerations. We conclude that $\FM_{g,n}(\FW,\beta)$  is a quasi-separated algebraic stack from the algebraicity and quasi-separatedness of mapping stacks, e.g.,\ see \cite[Theorem 1.2]{HR}.  Moreover, by the condition on the finiteness of automorphisms, \cite[Section 0DSL]{stacks-project} implies that it is Deligne--Mumford.

Since the stack of all FM degenerations is locally of finite type,  $\FM_{g,n}(\FW,\beta)$ is also locally of finite type by \cite[Theorem 1.2]{HR}. To establish that it is of finite type, we therefore have to show that it is quasi-compact. Since the moduli spaces of stable maps to a projective target are quasi-compact for a fixed topological data, and our maps are assumed to be admissible, it is enough to show that the family of FM degenerations appearing in $\FM_{g,n}(\FW,\beta)$ is quasi-compact, i.e.,\ FM degenerations have a bounded number of components. This in turn is equivalent to showing that the associated intersection graphs have a bounded number of chains.  The argument from \cite[Lemma 3.1.2]{KKO} can be adjusted to our setting. We, however, provide another argument below.

First, by composing a map $f \in  \FM_{g,n}(\FW,\beta)$ with the contraction,
\[ 
p\colon W \rightarrow X,
\]
and applying the stabilisation, we obtain a map  $f' \in \Mbar_{g,n}(X,\beta)$.  To reconstruct the map $f$ from $f'$, we inductively re-attach back the components of the FM degeneration $W$ (together with rational components of $C$), such that at each step  a new end component is attached. Using arguments from the preceding section and the stability conditions for maps in $\FM_{g,n}(\FW,\beta)$, we will show that the number of components that can be attached to maps in $\Mbar_{g,n}(X,\beta)$ to obtain a map in $\FM_{g,n}(\FW,\beta)$ is bounded from above. 

Let us first show that the length of chains is bounded. Consider the locus of maps without contracted components $\Mbar^\circ_{g,n}(X,\beta)$. We define the maximal intersection number of maps of genus $g'\leq g$ with the exceptional divisors $E_x$ of the blow-ups $\mathrm{Bl}_x(X)$ at points $x\in X$,  
\[ 
w_{\max}:=\mathrm{max}\{ \mathrm{deg}_{E_x}(f) \mid x\in X, f\in \Mbar^\circ_{g',  n}(X,\beta)   \}. 
\]
It is finite, since  $\Mbar_{g',n}(X,\beta)$ is bounded. By the admissibility, the map over the first end component $P_1$ that we attach to $\mathrm{Bl}_x (X)$ intersects the divisor at infinity with multiplicity $w_1$, such that 
\[
w_1 \leq w_{\max}.
\]  When we attach the next component $P_2$ to $\mathrm{Bl}_{x'}(P_1)$, the map must intersect the divisor at infinity of $P_2$  with multiplicity $w_2$, such that 
\[ 
w_2\leq w_1. 
\]
Moreover, the equality holds, if and only if the image of the map over $P_1$  is a union of lines passing through the point $x'$, because the divisor $\pi^*D-E$ on $\mathrm{Bl}_{x'} (P_1)$ vanishes precisely on the fiber lines (cf.\ the proof of Lemma \ref{mult}).  If the equality holds, in
order for the map to have finitely many automorphisms, it either should contain a marked point or a contracted component, or it is a ramified cover of lines with finitely many automorphisms (with respect to both the target and the source). Hence, we cannot have an infinite sequence of components $P_k$ with $w_{k}=w_{k-1} $, as with every
new component in such a sequence we increase the number of markings or the genus
of the source curve.   We conclude that the length of chains of components must be bounded from above by a number that depends only on $w_{\max}$, $g$, and $n$. 


It remains to argue that there cannot be infinitely many chains of finite length. This follows from the fact that, by Lemma~\ref{contraction} and the defining properties of maps in $\FM_{g,n}(\FW,\beta)$, new components can be attached only at base points, or the number of markings has to increase, or the genus of the source curve has to increase. By the boundedness of stable maps, the number of base points in $X$ and $\p^{\dim}$ is bounded above by  a number that depends only on $\beta$, $g$, and $n$.

\qed
\begin{lemma} \label{open}
The moduli space $\Mbar^\epsilon_{g,n}(X,\beta)$ is an open substack of $\FM_{g,n}(\FW,\beta)$. 
	\end{lemma}
\textit{Proof.} First, it is not difficult to see that $\epsilon$-unramification is a locally constructable condition by stratifying the space of FM degenerations into loci of degenerations with a fixed intersection graph. Hence we can use the valuative criterion of openness for locally constructible sets (or the valuative criterion of closedness for its complement).\footnote{By \cite[Section 0060]{stacks-project}, a constructable set is closed, if and only if it is stable under specialisation; by \cite[Lemma 054F]{stacks-project} the stability under specialisations can be check on spectrums of discrete valuations rings.} In fact, we will show that each of the conditions of $\epsilon$-ramification is open.  Definition \ref{epsilon} (4) is standard, hence we will deal with Definition \ref{epsilon} (1) and  Definition \ref{epsilon} (2).  Let $R$ be a discrete valuation ring and $K$ be its fraction field. Let 
\[
f \colon (\CC,\underline{p}) \rightarrow  \CW
\]
be a map to a FM degeneration over $\Spec(R)$. 

Assume that there exists  a  point $x^* \in   \CW^*$ in the generic fiber, such that 
\[
w(x^*)> 1/\epsilon.
\]  Both ramifications and genera of contracted components can only increase at the central fiber. Hence, the central fiber of the extension of the point $x^*$, denoted by $x\in  W$, satisfies 
\[ w(x)\geq w(x^*)> 1/\epsilon.\]
We conclude that a family  does not satisfy Definition \ref{epsilon} (1) at the central fiber if it does not satisfy it at the generic fiber. By the valuative criterion of closedness  for locally constructible subsets, we obtain that  Definition \ref{epsilon} (1) is open.

Assume there exits an end component $\CP^* \subset  \CW$, such that  
\[
w(\CP^*)\leq 1/ \epsilon.
\] Let $\CP$ be the closure of $\CP^*$ inside $\CW$. Let $P'$ be the central fiber of $\CP$; note that it might no longer be an end component. Consider the restriction of $f$ to  $\CP$, 
\[ f|_{\CP} \colon (\CC,\underline{p})|_{\CP} \rightarrow \CP.\] 
We choose an end component of $P \subset P'$. Let 
\[
P'_{\mathrm{chain}}\subset P'  \subset \CP
\] be the chain of components which connects $P$ to the root component, including $P$ itself. Contract all components in the central fiber $P'$ but $P'_{\mathrm{chain}}$, to obtain a family $\CP'$. By construction, the central fiber of $\CP'$ is chain of ruled components with the end component $P$. Note that it is not necessarily $\epsilon$-unramified, but it is $\epsilon$-unramified over $P$, because the contraction left $P$ intact. In particular, the weight of $P$ in $\CP'$ is equal to the weight of $P$ in $\CP$.  Composing the map $f|_{\CP}$ with the contraction 
\[
\CP \rightarrow \CP',
\] we obtain 
 \[f' \colon (\CC,\underline{p})|_{\CP}\rightarrow  \CP',\]
 such that $f|_{\CP^*}=f'|_{\CP^*}$. 
 The quantity 
\[\deg(\Omega^{\mathrm{log}}_{C}) -\deg(f^*\Omega^{\mathrm{log}}_{ W})/\dim+n/2\] is constant in families. In particular, 
 by Lemma \ref{ruledpos}, we conclude that 
 \[\omega(P) \leq \omega(\CP^*).\]
Since the contraction  $\CP \rightarrow \CP'$ did not change the value of $w(P)$,  the same holds for  the  family  \[ f|_{\CP} \colon (\CC,\underline{p})|_{\CP} \rightarrow \CP.\] Hence  a family  does not satisfy Definition \ref{epsilon} (2) at the central fiber if it does not satisfy it at the generic fiber. By the valuative criterion  of closedness  for locally constructible subsets, we obtain that  Definition \ref{epsilon} (2) is open. 
\qed 

\begin{cor} \label{ft}
	 The moduli space $\Mbar^\epsilon_{g,n}(X,\beta)$ is a quasi-separated Deligne--Mumford stack of finite type. 
\end{cor}
\textit{Proof.} The result follows immediately from  Lemma \ref{alg} and Lemma \ref{open}.
\qed 
\\

\begin{thm} \label{epsilonthm}  The moduli space $\Mbar^\epsilon_{g,n}(X,\beta)$ is a proper Deligne--Mumford stack. 
	\end{thm}
\textit{Proof, Step 1.} 
Let $R$ be a discrete valuation ring with the fraction field  $K$. Consider a family of $\epsilon$-unramified maps over $\Spec(K)$, 
\[
f^*\colon (\CC^*,\underline{p}^*) \rightarrow  \CW^*.
\]    We first show the existence part of the valuative criterion for quasi-separated Deligne--Mumford stacks. 

Let us introduce the following notation associated to the family $f^*$: 
\begin{itemize}
	\item $\CC^*_{\bullet}$ is  the union of components of $\CC^*$ contracted by $f^*$, 
	\item $\CC^*_{\circ}$ is the normalisation of the closure of the complement of $\CC^*_{\bullet}$, 
	\item $f^*_{\circ}$ is the restriction of $f^*$ to $\CC^*_{\circ}$, 
	\item  $\underline{p}^*_{\bullet}$ and $\underline{p}^*_{\circ}$ are the subsets of the markings $\underline{p}^*$ which lie on $\CC^*_{\bullet}$ and on $\CC^*_{\circ}$, respectively,
	\item  $\underline{n}^*_{\bullet}$ are the markings on $\CC^*_{\bullet}$ corresponding to nodes which separate it from $\CC^*_{\circ}$, 
	\item $\underline{n}^*_{\circ}$ are  the markings of $\CC^*_{\circ}$ corresponding to nodes which separate it from  $\CC^*_{\bullet}$ and the nodes which were resolved after taking the normalisation.
\end{itemize}

The   map 
\[
f^*_{\circ} \colon (\CC_{\circ}^*,  \underline{p}_{\circ}^*,  \underline{n}_{\circ}^*)   \rightarrow  \CW^* 
\]
does not have contracted components. Hence, by Proposition \ref{propexitence}, there exists an extension 
\[ 
f_{\circ} \colon (\CC_{\circ},  \underline{p}_{\circ},  \underline{n}_{\circ})   \rightarrow  \CW,
\]
whose ramification profiles are constant in the family, and the markings are contained in the same fiber of $f_\circ$ over the central fiber of the family if they are contained in the same fiber of $f_\circ$  over the generic fiber of the family. Moreover, we may assume that $f_\circ$ has finitely many automorphisms with respect to the target and the source by contracting unstable components.    

The components $(\CC_{\bullet}^*,\underline{p}_{\bullet}^*, \underline{n}_{\bullet}^*)$, contracted at the generic fiber,  define a family of stable curves. Hence, by the properness of stable curves, again possibly after a finite base change, there exists  an extension $(\CC_{\bullet},\underline{p}_{\bullet}, \underline{n}_{\bullet})$ over $\Spec(R)$. We take a common finite base change for both extensions. 

We then glue back the two families of curves via the markings $\underline{n}_{\bullet}$ and $\underline{n}_{\circ}$. We also glue the markings from $\underline{n}_{\circ}$ to each other if they correspond to the normalised nodes. As the result, we obtain the family of maps over $\Spec(R)$, 
 \[ f \colon (\CC,  \underline{p})  \rightarrow  \CW,\]
 which remains admissible, as the gluing is done away from the singularities of $\CW$. Moreover, it is flat over $\Spec(R)$, since every irreducible component $\CC$ dominates $\Spec(R)$. The family $(\CC_{\bullet},\underline{p}_{\bullet}, \underline{n}_{\bullet})$ is contracted by the map $f_\circ$ to the images of some markings $\underline{n}_\circ$. Observe now that by the defining properties of the family $f_\circ$, the weights of all points remain constant. More precisely, for a base point $x^* \in \CW^*$, the central fiber of its extension $x\in W$ satisfies  
 \[
 w(x^*)=w(x).
 \]
Moreover, all base points of the central fiber are extensions of the base points of the generic fiber. In particular, the map $f$ satisfies Definition \ref{epsilon} (2) by the $\epsilon$-unramification of the generic fiber, as well as Definition \ref{epsilon} (3), (4). However, it might not satisfy Definition \ref{epsilon} (1). We deal with this next. 
 \\
 
 \textit{Step 2. }
Our construction of
\[
f \colon (\CC,\underline{p}) \rightarrow \CW
\]
does not exclude the creation of end components in the central fiber,
\[
P \subset W \subset \CW,
\]
that violate the $\epsilon$-unramification condition,
\[
w(P)\leq \epsilon,
\]
since $f_\circ$ was assumed to be unramified away from the ramification points arising from the generic fiber.

 Assume  that there exists such an end component $P$.   By the $\epsilon$-unramification  of the generic fiber, this end component cannot be the limit of an end component at the generic fiber, i.e., there does not exist an end component  $\CP^* \subset \CW^*$, whose closure inside $\CW$ gives $P$ at the central fiber. Hence, as, explained in Section \ref{blowupscont}, we can contract $P$. By Lemma \ref{contraction}, the image of the contraction $x$ satisfies 
\[
w(x)\leq 1/\epsilon.
\]

By repeating this procedure inductively, we can get rid of all end components which do not satisfy $\epsilon$-unramification without creating any points that violate $\epsilon$-unramification. Moreover, the results of such contractions do not violate Definition \ref{epsilon} (3), (4). The admissibility is  clearly preserved, and  after contracting a component, the image of the contraction satisfies $w(x)>0$, hence it is a base point, and if the component on which it lies had finitely many automorphisms before the contraction, then it still has finitely many automorphisms after the contraction. Overall, after inductively  contracting all such end components, we arrive at the family which satisfies all conditions of $\epsilon$-unramification, Definition \ref{epsilon}. 
\\

\textit{Step 3.} We will now deal with the uniqueness part of the valuative criterion.  Consider two families  of $\epsilon$-unramified maps over $\Spec(R)$,
 \[f_i\colon (\CC_i,\underline{p_i}) \rightarrow  \CW_i, \quad i=1,2,\] whose generic fibers are isomorphic. Possibly after a finite base change, we can find a family 
 \[\tilde{f}\colon (\tilde\CC,\tilde{\underline{p}}) \rightarrow  \widetilde{\CW},\] which dominates both.  It is given by taking a dominating family of the FM degenerations, as explained in Section \ref{blowupscont},  and extending maps (possibly after a base change, and  blow-ups on the source curves); by the separatedness of maps to a fixed target, both $f_j$ must extend to the same map $\tilde{f}$.  In particular, there exist contraction maps
 \[
  \pi_i \colon  \widetilde{\CW} \rightarrow  \CW_i, \quad i=1,2,
   \]
 whose compositions with $\tilde{f}$ give $f_i$.   By induction, we may assume that $\pi_i$ contracts either one end component or one  ruled component (and not a tree of components). By Lemma \ref{contraction}, contractions must introduce points whose weights are equal to the weights of contracted components. Hence,  by the $\epsilon$-unramification of both families, if an end component is contracted by $\pi_1$, it is also contracted by $\pi_2$, and vice versa. Similarly, if a ruled component is contracted by $\pi_1$, then it must be unstable by the admissibility of $f_1$ over the singularities of the target (i.e., the map over this component has infinitely many automorphisms); indeed, by the local description of admissibility, if there is a ruled component over a singularity of the target, such that the map extends to it, then, over this component, the map must be given by covers of the fiber lines with no ramifications over the smooth locus.  This  ruled component therefore must also be contracted by $\pi_2$, and vice versa.  We conclude that $ \CW_1$ and $\CW_2$ must be isomorphic. The claim then follows from the uniqueness of extensions of stable maps to a fixed target. 
\qed

\begin{prop} \label{propexitence} Let $R$ be a discrete valuation ring with the fraction field  $K$. Let
	\[
	f^* \colon (\CC^*,\underline{p}^*) \rightarrow \CW^*
	\]
	be a family of admissible maps over $\Spec(K)$. Assume that $f^*$ is ramified over  points $\{x^*_1, \hdots, x^*_m\}$ of the target $\CW^*$ and is unramified elsewhere, with ramification profiles given by the vector of partitions $(\eta^1,\hdots,\eta^m)$. Then, possibly after a finite base change, there exists an extension
	\[
	f \colon (\CC,\underline{p}) \rightarrow \CW
	\]
	over $\Spec(R)$ whose central fiber is an admissible map, such that it has ramification profiles $(\eta^1,\hdots,\eta^m)$ over the extensions of the points $\{x^*_1,\hdots,x^*_m\}$ and is unramified elsewhere. Moreover, the marked points $\underline{p}$ lie in the same fiber of $f$ over the central fiber of the family, if and only if they lie in the same fiber of $f$ over the generic fiber of the family.
\end{prop}

\begin{proof}
The arguments presented in \cite[Proposition 4.2.1]{KKO} can be adapted to prove this result. We, however, give a slightly different proof based on the properties of the weights established in Lemma \ref{less} and Lemma \ref{contraction}. We begin by addressing the ramifications of the extensions.
 
First, blow up the points $x_j^* \in \CW^*$ to obtain a \textit{blown-up} FM degeneration with $m$ exceptional divisors, 
\[ 
\CE_j^* \subset \CW^*(m):=\mathrm{Bl}_{x_1^*,\hdots,x_m^*}(\CW^*), \quad j=1,\hdots, m.
\]
See Section \ref{sectionrelative} for more details about this construction. 
The map $f^*$ extends to a map to $\CW^*(m)$, 
\begin{equation} \label{genadm}
f^* \colon (\CC^*,\underline{p}^*) \rightarrow \CW^*(m),
\end{equation}
such that the intersection multiplicity of $\CC^*$ with $\CE_j^*$ is equal to $\eta^j$, cf.\  (\ref{multexc}). Next, possibly after a finite base change, by the properness of marked FM degenerations, we can choose an extension of $(\CW^*,x_1^*,\hdots,x_m^*)$ over $\Spec(R)$, such the marked points remain disjoint. Blowing up these marked points then yields a blown-up FM degeneration over $\Spec(R)$,
\[
\CE_j \subset \CW(m)'.
\]
By the properness of stable maps to projective targets,  possibly after a finite base change, there exists an extension of the map 
\[
f' \colon (\CC',\underline{p}') \rightarrow \CW(m)',
\]
which, however, might not be admissible over the singularities of the central fiber of $\CW(m)'$ and over the exceptional divisors $\CE_j$ (i.e., there is a contracted component over $\CE_j$, or the ramification index jumps), and  might also have ramifications elsewhere. The admissibility is local with respect to the target; since the generic fiber (\ref{genadm}) is admissible, it can be achieved by performing expansions (i.e., blow-ups) along the singularities and the smooth divisors $\CE_j$ of the central fiber, following Li's proof of  properness of relative stable maps \cite[Section 3]{LiRel}. Here, we rely on the fact that the singularities of the FM degenerations are of the form (\ref{an0}). Hence, by \cite[Lemma 3.9]{LiRel},  possibly after a finite base change, there exists
\[
f \colon (\CC,\underline{p}) \rightarrow \CW(m),
\]
an extension whose central fiber is admissible, and the ramification profiles over $\CE_j \subset   \CW(m)$ are equal to $\eta^j$. Note that we allow the markings of $\CC$ to map to $\CE_j$ whenever they map to $\CE_j^*$. We refer to Remark \ref{KKOproof} for an alternative argument from \cite[Proposition 4.2.1]{KKO} establishing admissibility in this case.

Next, assume $f$ is ramified over a regular point $x$ of the central fiber of $\CW(m)$ away from $\CE_j$.  We distinguish two cases: (1) $f$ has a contracted component over $x$, and (2) $f$ is ramified without contracted component over $x$. 
\\
 
 	\textit{Case 1.} Assume $f$ contracts an irreducible component $C_{\bullet}$ of the central fiber of $\CC$ to a point $x\in \CW(m)$. By inductively blowing up the image point of the contraction, we can move this contracted component onto an end component of the target. More precisely, following the argument from \cite[Lemma 4.1.1]{KKO}, this is achieved as follows.  First, possibly after taking an \'etale cover of $\Spec(R)$, we equip $\CC$ with two auxiliary markings 
 	\[
 	\sigma_1, \sigma_2 \colon \Spec(R) \rightarrow \CC,  
 	\]
 	which
 	\begin{itemize}
 		\item do not intersect and lie in the regular locus of the central fiber of  $\CC$,
 		\item are  contained in $C_{\bullet}$ at the central fiber, 
 			\item have distinct images with respect to $f$ at the generic fiber. 
 			\end{itemize}
The existence of sections satisfying the first two conditions follows from the smoothness of $\CC \to \Spec(R)$ around a regular point of the central fiber. The third condition can be arranged by first choosing such a section~$\sigma_1$ and then removing the generic fiber $f^{-1}(f(\sigma_1))^*$ from~$\CC$. The resulting curve,
$
\CC \setminus f^{-1}(f(\sigma_1))^*,
$
is still smooth over~$\Spec(R)$; therefore, the existence of a section~$\sigma_2$ satisfying the third condition follows as well.
 
 Consider now the images of $\sigma_1$ and $\sigma_2$ with respect to $f$, which define sections of $\CW(m)$ intersecting at the central fiber, 
 \[
 f(\sigma_1), f(\sigma_2) \colon \Spec(R) \rightarrow \CW(m). 
 \]
 By the properness of marked FM degenerations, there exists $\CW(m)'$, isomorphic to $\CW(m)$ at the generic fiber, such that the sections $f(\sigma_1), f(\sigma_2)$ become disjoint at the central fiber. Such $\CW(m)'$ is constructed by an iterated blow-up of $\CW(m)$ followed by contraction of ruled components; the central fiber of $\CW(m)'$ has an additional end component.  Let  
 \[
 f' \colon (\CC',\underline{p}') \rightarrow \CW(m)'
 \]
be the extension of $f$ to $\CW(m)'$, which exists possibly after a finite base change; it is equal to $f$ away from the new  end component, and $\CC'$ is a blow-up of $\CC$. Moreover, (the proper transform of) $C_{\bullet}$ is no longer contracted by  construction of $\sigma_1$, $\sigma_2$, and $\CW(m)'$; indeed, the central fibers of $f'(\sigma_1), f'(\sigma_2)$ are two distinct points. However, $f'$ might not be admissible over the singularity at which the new end component is attached, for which we apply Li's expansions. Let 
 \[
 f'' \colon (\CC'',\underline{p}'') \rightarrow \CW(m)''
 \]
 be the new admissible map after expansions at the singularities. Overall, we have introduced a chain of ruled components with an end component $P \subset \CW(m)''$. Moreover, $C_{\bullet}$ was a stable contracted component, i.e., it was either a non-rational curve, or a rational curve with at least three special points.  It therefore cannot be  mapped by $f''$ onto a line with one base point over $P$. Assuming no ruled components were introduced in between $P$ and the rest of $ \CW(m)''$,  by Lemma \ref{less} and Lemma \ref{contraction}, the weights of all points on $P$ with respect to $f''$ are therefore strictly smaller than the weight of $x$ with respect to $f$. If there are ruled components, then by contracting  end components, we can make any ruled component into an end component, and again use Lemma \ref{less} and Lemma \ref{contraction}; hence, the same conclusion applies to the ruled components. We thereby replaced the point $x$ with a chain of components whose points have strictly smaller weights, preserving the admissibility of the map. 
 \\
 
 \textit{Case 2.} Assume $f$ has a ramification without contracted component over a point at the central fiber $x\in \CW(m)$. By \cite[Lemma 4.2.2]{KKO} and the unramification of the generic fiber, the projectivisation of the derivative 
 \[
 \p(f)  \colon \CC \dashto \p(TX)
 \]
 does not extend at points of $\CC$ over $x$ (since the ramification order jumps). Choose one such point. By the resolution of indeterminacies of maps from surfaces, after inductively blowing up $\CC$ at this point, one eventually reaches a stage at which the exceptional curve of the final blow-up either contains more than one indeterminacy, or $\p(f)$ extends to a non-constant map away from the indeterminacies, or it contains a marking distinct from an indeterminacy point. Let $\CC'$ denote the resulting surface, obtained by performing these blow-ups and then contracting the exceptional curves arising from the earlier blow-ups.
 Denote 
 \[
  \p(f)' \colon \CC' \dashto \p(TX)
 \]
the rational map given by $\p(f)$.  Consider now the composition 
\[
f' \colon (\CC',\underline{p}') \rightarrow \CW(m)
\]
of the blow-down $C' \rightarrow C$ with $f$. The map~$f'$ has a contracted rational component, which we denote by~$C'_\bullet$. Using the same procedure as in \textit{Case~1}, we introduce an end component into the central fiber of~$\CW(m)$, onto which the proper transform of~$C'_\bullet$ is mapped. Let
\[
f'' \colon (\CC'',\underline{p}'') \rightarrow \CW(m)'
\]
be the resulting map.

Note that
\[
\p(f)'=\p(f'')
\]
around the component~$C'_\bullet$, since both maps agree at the generic fiber. Hence, by construction of $C'_\bullet$ and \cite[Lemma~4.2.2]{KKO}, the map~$f''$ cannot send~$C'_\bullet$ onto a line with a single base point. After making the map admissible over the singularities as in \textit{Case 1}, by Lemma \ref{less} and Lemma \ref{contraction}, we therefore arrive at
\[
f''' \colon (\CC''',\underline{p}''') \rightarrow \CW(m)'',
\]
where $\CW(m)''$ contains an additional end component together with a chain of ruled components, such that the weights of all points on these component with respect to $f'''$ are strictly smaller than the weight of $x$ with respect to $f$. 
\\

After repeating the procedures in \textit{Cases 1, 2} and relabelling the maps, we must eventually arrive at a map
\[
f \colon (\CC,\underline{p}) \rightarrow \CW(m),
\]
such that the weights of all points away from~$\CE_j$ are zero, i.e., $f$ is unramified away from~$\CE_j$. Indeed, this follows because the weights of the points strictly decrease at each iteration. After blowing down~$\CE_j$, we obtain a map
\[
f \colon (\CC,\underline{p}) \rightarrow \CW
\]
whose ramification profiles are constant throughout the family.

Finally, suppose that two markings from~$\underline{p}$, say $p_1$ and~$p_2$, lie in the same fiber of~$f$ over the central fiber of the family, while their generic fibers are distinct. Then their images
\[
f(p_1),\, f(p_2) \colon \Spec(R) \rightarrow \CW
\]
are intersecting markings of the target~$\CW$. Using the properness of the moduli space of marked FM degenerations, we can separate the sections~$f(p_1)$ and~$f(p_2)$. To preserve the admissibility of the map, further expansions, as in the previous steps, might be necessary.
\end{proof}	


\begin{rmk} \label{KKOproof}
In \cite[Proposition 4.2.1]{KKO}, an admissible extension is constructed by introducing auxiliary markings to prevent contracted components over the singularities of the central fiber, as is done in \textit{Cases 1, 2} of the proof above; cf.\ \cite[Lemma 4.1.2]{KKO}. By \cite[Proposition 2.2]{LiRel}, once there are no contracted components, admissibility is automatic. Note that further expansions do not produce new contracted components lying over the singularities, but only components mapping to the smooth locus of newly created ruled components. Indeed, if there were a new contracted component of this kind, then it would already have been present before the expansion, since it would not be destabilized by the contraction of a ruled component. Hence, the process terminates after finitely many steps. The same argument can also be applied to establish admissibility over smooth divisors.
\end{rmk}	

\section{Moduli spaces of  relative $\epsilon$-unramified maps} \label{sectionrelative}
\subsection{Maps with prescribed ramifications} In the theory of admissible covers of curves \cite{HM,ACV}, we have liberty to consider covers with prescribed ramification profiles over points on the \textit{target}, which are required to be constant in families. The higher-dimensional generalisation of ramified admissible covers will be approached by blowing up points in the target and using the theory of maps relative to exceptional divisors  in the sense of Li \cite{LiRel}. These spaces will be essential for the statement of the wall-crossing formula, Theorem \ref{wallcrossing}.    

In \cite[Definition~3.1.1]{KKO}, ramified maps with prescribed ramification at points on the \emph{source} were introduced. However, they are not sufficient for our purposes, since the wall-crossing construction involves gluing components of the target, which in turn requires the ramification profiles on both sides of the gluing to agree. This is analogous to the theory of admissible covers: unlike in Gromov--Witten theory, gluing maps requires gluing both the sources and the targets while preserving admissibility at the nodes.\footnote{In other words, Hurwitz theory defines a cohomological field theory only with respect to the gluing of targets.}

\subsection{Blown-up FM degenerations}	
Let 
\[
( W, x_1,\dots,  x_m)
\]be a FM degeneration of $X$ with $m$ distinct regular points $x_j \in  W^{\mathrm{sm}}$. We blow up the variety $W$ at $x_j$,
\[ W(m):=\mathrm{Bl}_{x_1, \dots, x_m} (W),\]
and let 
\[
E_j \subset  W(m), \quad  j=1, \hdots, m, 
\] be the associated exceptional divisors.  We call $ W(m)$ a \textit{blown-up FM degeneration}, or simply a FM degeneration, when it is clear from the context whether it is blown-up or not.  

 To each blown-up FM degeneration, we can similarly associate  a tree graph with one-vertex edges corresponding to exceptional divisors of blown-up components. As before, an \textit{end component} of $ W(m)$ is a component which corresponds to a vertex with one edge in the associated tree graph, excluding the root vertex. A \textit{ruled component} of $ W(m)$ is a component which corresponds to a vertex with two edges, excluding the root vertex. Note that a ruled component can  now  be a blown-up end component, we refer to such components as \textit{marked ruled components}.  It is productive to think of $W(m)$ as just marked FM degenerations. In particular, below, a family of blown-up FM degenerations is defined as a family of FM degenerations with $m$ sections.   An example of a tree associated to a FM degeneration is depicted in  Figure \ref{tree2}.

\begin{defn}A family of blown-up FM degenerations is defined to be a family of FM degenerations in the sense of Section \ref{FMdeg} together with $m$ sections, 
	\[
	\begin{tikzcd}[row sep=normal, column sep = normal]
		\CW  \arrow[r,"\pi_B"] &  B, \arrow[l,bend left=25, "\sigma_{j}"]  &\hspace{-0.8cm} j=1,\dots, m,
	\end{tikzcd}
	\]
which we blow up. 
 Let  \[\FM_{X,m}\colon (Sch/\BC)^\circ \rightarrow Grpd\] 
  denote the moduli space of blown-up FM degenerations\footnote{By definition, this is simply a moduli space of FM degenerations with $m$ markings.  }
   of $X$. Let
  \[
  \mathfrak{W}(m) \rightarrow \FM_{X,m} 
  \]
  be the universal blown-up FM degenerations, and
  \[ 
  \mathfrak{E}_j  \subset  \mathfrak{W}(m), \quad j=1, \hdots, m,
  \]
 be the universal exceptional divisors,   which are the blow-ups of the universal sections. 

  \end{defn}

\begin{figure}[!ht]
	\centering
	\begin{tikzpicture}		
		\node (0) at (0.25,1) {};
		\node (1) at (-2,2) {};
		\node (2) at (-0.5,2) {};
		\node (3) at (1,2) {};
		\node (4) at (2.5,2) {};
		\node (6) at (0.25,3) {};
		\node (7) at (1.75,3) {};
		\node (8) at (1, 0.5) {};
		\node (9) at (1, 1.5) {};
		\node (10) at (2.7,3) {};
		\node (11) at (1, 3.6) {} ;

		\draw[black] (0.25,3)--(1.75,3);
		\draw[black] (0.25,1)--(-0.5,2);
		\draw[black] (-2,2)--(-0.5,2);
		\draw[black] (-0.5,2)--(0.25,3);
		\draw[black] (-0.5,2)--(1,2);
		\draw[black] (1,2)--(2.5,2);
		\draw[black] (0.25,1)--(0.80, 1.4);
		\draw[black] (0.25,1)--(0.80, 0.60);
		\draw[black] (1.75,3)--(2.45,3);
		\draw[black] (0.25,3)--(0.80, 3.40);

		
		(6) -- (11);

		\filldraw[thick, fill = black] (0.25,1) circle (.1cm) node at (0) {};
		\filldraw[thick, fill = black] (-2,2) circle (.1cm) node at (1) {};
		\filldraw[thick, fill = black] (-0.5,2) circle (.1cm) node at (2) {};
		\filldraw[thick, fill = black] (1,2) circle (.1cm) node at (3) {};
		\filldraw[thick, fill = black] (2.5,2) circle (.1cm) node at (4) {};
		\filldraw[thick, fill = black] (0.25,3) circle (.1cm) node at (6) {};
		\filldraw[thick, fill = black] (1.75,3) circle (.1cm) node at (7) {};

		\node at (2.75, 3.4) {$\mathrm{marked \ ruled}$};
		\node at (-1.7, 2.4) {$\mathrm{root}$};
		\node at (1.4, 2.4) {$\mathrm{ruled}$};
		\node at (2.75, 2.4) {$\mathrm{end}$};

	\end{tikzpicture}
	\caption{Tree with markings.} \label{tree2}
\end{figure}

 \subsection{Relative $\epsilon$-unramified maps}Weight functions $w(x)$ and $w(P)$ are defined in the same way for points
 \[
 x\in  W(m)^{\mathrm{sm}} \setminus \cup_j E_j
 \] and end components $P \subset W(m)$ as in Section \ref{weightsdef}.  Associated to blown-up FM degenerations, we will consider the following relative $\epsilon$-unramified maps. For 
 \[
 1<\epsilon \leq 2,
 \] they are a higher-dimensional generalisation of admissible covers with prescribed ramifications at marked points of the target \cite{HM,ACV}.

\begin{defn} \label{relative} Choose a real number $\epsilon \in \BR_{>0}$. A map $f \colon (C,\underline{p}) \rightarrow  W(m)$ is said to be \textit{relative} $\epsilon$-$unramified$ if 
\begin{itemize}
	\item[1)] for all end components $P \subset  W(m)$, $w(P)> 1/\epsilon$,
\item[2)] for all $x\in  W(m)^{\mathrm{sm}}\setminus \cup_j E_j$, $w(x) \leq 1 /\epsilon$,
\item[3)] the map $f$ is admissible over $ W(m)^{\mathrm{sing}}$ and $\cup_j E_j$,
\item[4)] $|\Aut_{C,W(m)}(f)|<\infty$ and all connected components map non-trivially. 

	\end{itemize}

	\end{defn}



 \begin{defn} \label{def} Let 
\[\Mbar^\epsilon_{g,n,m}(X,\beta, \underline{\eta}) \colon (Sch/\BC)^\circ \rightarrow Grpd \] 
be the moduli space of relative $\epsilon$-unramified maps of degrees $\beta$ from possibly disconnected curves of genus $g$ with $n$ markings and $m$ exceptional divisors on the target.  Ramifications over the exceptional divisors $E_j \subset W(m)$ are specified by a vector of partitions 
\[
\underline{\eta}=\{\eta^1,\dots, \eta^m\},
\]
where 
\[
\eta^j=(\eta^{j}_1,\dots, \eta^{j}_{\ell(\eta^j)}) \in \BZ_{>0}^{\ell(\eta^j)}.\]
 If $\eta^j=(1)$ for all $j$, then we drop ramification profiles from the notation, 
\[
\Mbar^\epsilon_{g,n,m}(X,\beta):=\Mbar^\epsilon_{g,n,m}(X,\beta, \underline{\eta}) \quad \text{ if }   \eta^j=(1) \text{ for all }j.
\]
For every exceptional divisor $E_j$, the points of contact $f^{-1}(E_j)=\{p_{j,1}, \hdots, p_{j,i_{\ell(\eta^j)}}\}$ are  ordered in agreement with their multiplicities, i.e.,\  for any two points $p_ {j,i_1}$ and $p_{j,i_2}$ in $f^{-1}(E_j)$, we have 
\[
i_1\geq i_2 \quad \text{ if }    \eta_{i_1}^j\geq \eta_{i_2}^j.
\] 
 We call such order a \textit{standard order}. There are $|\Aut(\eta^j)|$ ways of putting a standard order on points in the fiber $f^{-1}(E_j)$. 
\end{defn}
Spaces $\Mbar^\epsilon_{g,n,m}(X,\beta, \underline{\eta})$ will be the intermediate moduli spaces needed for the wall-crossing.  Note that if $\dim(X)>1$, then the class $\beta$ bounds the size of ramification profiles, but does not determine it, unlike in the case $\dim(X)=1$.  For example, if we consider degree $3$ stable maps from genus $0$ curves to $\p^2$,
\[
\Mbar_0(\p^2,3),
\]
then the sizes of possible ramifications range from $1$ to $3$, 
\[
1\leq |\eta|\leq 3,
\] 
such that $|\eta|=1$ arises  for an embedding of a smooth curve, $|\eta|=2$ is associated to the cusp of a cuspidal cubic, and $|\eta|=3$ is given by degree $3$ covers of a line. In general, the size $|\eta|$ will depend on the multiplicities of singularities of curves inside the class $\beta$ and the multiplicity of the class itself. This contrasts, of course, with
\[
\Mbar_0(\p^1,3),
\]
for which all maps have $|\eta|=3$.


\begin{thm} \label{thmrel} The moduli space $\Mbar^\epsilon_{g,n,m}(X,\beta, \underline{\eta})$ is a proper Deligne--Mumford stack. 
\end{thm}
\textit{Proof.} The proofs of Lemma \ref{alg} and Lemma \ref{open} readily extend to show that the moduli space $\Mbar^\epsilon_{g,n,m}(X,\beta, \underline{\eta})$ is a quasi-separated Deligne--Mumford stack of finite type. The argument presented in the proof of Theorem~\ref{epsilonthm} also holds almost unchanged for $\Mbar^\epsilon_{g,n,m}(X,\beta, \underline{\eta})$, as we now explain.

 Consider a family
 \[
 f \colon (\CC^*,\underline{p}^*) \rightarrow \CW(m)^*
 \]
 of relative $\epsilon$-unramified maps over $\Spec(K)$. Using the same procedure as in \textit{Step~1} of the proof of Theorem~\ref{epsilonthm}, we decompose $f$ into the map
 \[
 f^*_{\circ} \colon (\CC_{\circ}^*, \underline{p}_{\circ}^*, \underline{n}_{\circ}^*) \rightarrow \CW(m)^*,
 \]
 which has no contracted components, and the contracted part
 \[
 (\CC_{\bullet}^*,\underline{p}_{\bullet}^*, \underline{n}_{\bullet}^*).
 \]
 The exceptional divisors of~$\CW(m)^*$ can be viewed as markings of the family~$\CW^*$, over which $f^*_{\circ}$ has fixed ramification profiles. Hence, Proposition~\ref{propexitence} applies in this setting. Indeed, its proof already involved viewing points on FM degenerations as exceptional divisors. The remainder of the argument is identical to that of Theorem~\ref{epsilonthm}.
 \qed

	\section{Weighted Fulton--MacPherson degenerations} \label{sectionweights}

\subsection{Definition of weighted Fulton--MacPherson degenerations}
For the construction of master spaces in Section \ref{masterspaces}, we require a notion of weights of the end components for families of FM degenerations. In particular, these weights must be well defined with respect to smoothings of end components. To achieve this, we will need to equip all irreducible components of FM degenerations with weights, which is the subject of this section.

\begin{defn} \label{weighted}
	Given a blown-up FM degeneration $W(m)$, let 
	\[ W(m)=   W(m)_1 \cup \hdots \cup  W(m)_k\] be the decomposition  of $ W(m)$ into its irreducible components.  A \textit{weighted blown-up FM degeneration}, 
	\[( W(m), \underline d):=( W(m), (d_1, \dots, d_k)),\] is a  blown-up FM degeneration  $ W(m)$ with a weight 
	\[d_i \in \frac{1}{2\dim}\BZ\] assigned to every irreducible component $ W(m)_i$. The total weight of $( W(m),\underline{d})$ is 
	\[
	d:=\sum_i d_i. 
	\]
\end{defn}

	A family of weighted FM degeneration $( \CW(m),\underline{d})$ is a family of FM degenerations with weights assigned to irreducible components of every  fiber 
	\[
	\CW(m)_b=\CW(m)_{b,1} \cup \hdots \cup  \CW(m)_{b,k}
	\] over a geometric point $b\in B$,
	\[\underline{d}(b)=(d_1,\dots, d_k).\]
	We require that, \'etale locally on the base $B$, there exist two collections of sections of the family $\CW(m)$,
	\[\{\sigma_{-,1}, \dots, \sigma_{-,k_-}\} \quad \text{and} \quad 
	\{\sigma_{+,1}, \dots, \sigma_{+,k_+}\},\]  disjoint from the singular locus and exceptional divisors, such that 
	\[ \underline{d}(b)=n_+\sum_i \underline{\deg}(\sigma_{+,i}(b))-n_-\sum_i \underline{\deg}(\sigma_{-,i}(b))\] 
for all geometric points $b \in B$. Here, 
	\[
	\underline{\deg}(...)
	\] stands for the vector of degrees of points on each component, and $n_+$ and $n_-$ are some rational numbers. We call this condition \textit{continuity} of weights. 

\begin{defn} \label{moduliwieght}
Let
\[
\FM_{X,m,d} \colon (Sch/\BC)^\circ \rightarrow Grpd
\]
be the moduli space of weighted blown-up FM degenerations of $X$ of total weight $d$, such that:
\begin{itemize}
	\item[(1)] the number of irreducible components and their weights are bounded,
	\item[(2)] the weights of end components are at least $d_0>0$, where $d_0$ is a fixed number,
	\item[(3)] the weights of all ruled components of a FM degeneration and of all its smoothings are non-negative.
\end{itemize}
   The bounds and $d_0$ will depend on the moduli spaces of maps that we consider, hence, we suppress them from the notation. 
\end{defn}

Note that we allow negative weights on components that are neither end nor ruled. We refer the reader to Remarks~\ref{weight0} and~\ref{polar} for the reasons why we define moduli spaces of weighted FM degenerations in this way.

The last condition in Definition~\ref{moduliwieght} means that if we smooth a singularity of a weighted FM degeneration, then the resulting FM degeneration must also have ruled components of non-negative weight. This is an open condition because the locus of smoothings of a FM degeneration is open. We must impose it; otherwise, the smoothing of ruled components could be obstructed by the non-negativity requirement in the presence of negative weights. Equivalently, this condition can be expressed purely in terms of the sums of the weights of the various components; see the proof of Lemma~\ref{lemmaweights}.

  The space $\FM_{X,m,d}$ admits a projection, 
\begin{equation} \label{weightsprojection}
 \FM_{X,m,d} \rightarrow \FM_{X,m},
\end{equation}
by forgetting the weights.  Open neighbourhoods of FM degenerations are given by smoothings of the simple normal crossing singularities. The weights of smoothings of a FM degeneration are uniquely determined by the weights of the FM degeneration by the continuity of weights. Hence, the projection is \'etale.  The space $\FM_{X,m}$  is a smooth algebraic stack with a normal crossing boundary, the same therefore holds for $\FM_{X,m,d}$.

\subsection{Weights associated to maps} Consider a map 
\[
[f\colon (C,\underline{p}) \rightarrow  W(m)] \in \Mbar^\epsilon_{g,n,m}(X,\beta,\underline{\eta}).
\]
 The target $ W(m)$ can be naturally weighted by the
degrees of the log cotangent bundles associated to $W(m)$ together with exceptional divisors $E_j$, 


\begin{equation} \label{weights}
  W(m)_\ell \mapsto d_\ell \in \frac{1}{2\dim}\BZ,
 \end{equation} 
where 
 \begin{multline*} 
  d_\ell= \deg\big(\Omega^{\mathrm{log}}_{C}\big(\sum_{j,i} p_{j,i}\big)\big|_{W(m)_\ell}\big) \\
  -\deg\big(f^*\Omega^{\mathrm{log}}_{W(m)}\big(\log\big(\sum_{j} E_j \big)\big)\big|_{W(m)_\ell}\big)\big/\dim+n_{ W(m)_\ell}/2.
\end{multline*}

\noindent By Lemma \ref{weightP}, the total weight is given by 
\begin{align*}
	d=2g-2+ \beta \cdot \mathrm{c}_1(X)/\dim+n/2+\sum_j (\ell(\eta^j)-|\eta^j|). 
\end{align*}
Also by Lemma \ref{weightP}, the weights above evaluated at end components recover the weights given by the weight function $w(P)$, defined in Definition \ref{weightend}. 

\begin{lemma} \label{lemmaweights}The weights given by (\ref{weights}) define a morphism 
	\[\Mbar^\epsilon_{g,n,m}(X,\beta, \underline{\eta}) \rightarrow \FM_{X,m,d}.\] 
	\end{lemma}

\textit{Proof, Step 1.}
 Weights associated to admissible maps from curves to FM degenerations behave well in families because we are taking degrees of vector bundles  on curves. 
However, we have to comment on how to align this construction with Definition \ref{weighted}. In fact, continuity of weights in Definition \ref{weighted} is equivalent to the following two conditions which must be satisfied with respect to generalisatons of geometric points in the space of FM degenerations: 
\begin{itemize} 
	\item[(1)] weights of two adjacent components add up whenever a singularity between them is smoothed out,
	
	\item[(2)] weights are constant whenever no singularities are smoothed out or no new components are created.
	
\end{itemize}

It is clear that the above condition are implied by the continuity condition in Definition \ref{weighted}. Let us show the converse. Let 
\[
\CW(m) \rightarrow B
\] be a family of FM degenerations with weights satisfying the conditions (1) and (2) above.  Away from the singular locus of the fibers, the family is smooth. 
Hence,  \'etale locally on the base $B$, it admits sections which extend any given point away from the singular locus on any given closed fiber $ \CW(m)_{b} \subset  \CW(m)$. Let 
\[
\{\sigma_{-,1}(b), \dots, \sigma_{-,k_-}(b)\} \quad \text{and} \quad \{\sigma_{+,1}(b), \dots, \sigma_{+,k_+}(b)\}
\] be two sets of distinct points on a fiber $ \CW(m)_{b}$, away from the singular locus and the exceptional divisors, defined as follows. The set of points $\{\sigma_{+,1}(b), \dots, \sigma_{+,k_+}(b)\}$ is supported on components with positive weights, and 
\[
n_{+}\sum_i\underline{\deg}(\sigma_{+,i}(b))
\]
 is equal to these weights. While the set of points $\{\sigma_{-,1}(b), \dots, \sigma_{-,k_-}(b)\}$ is supported on components with negative weights, and   
 \[
 -n_{-}\sum_i\underline{\deg}(\sigma_{-,i}(b))
 \] is equal to these weights. We now extend these points to two collections of sections 
 \[
 \{\sigma_{-,1}, \dots, \sigma_{-,k_-}\} \quad \text{and} \quad \{\sigma_{+,1}, \dots, \sigma_{+,k_+}\} 
 \] on some \'etale neighbourhood around $b \in B$, such that this neighbourhood contains at most only smoothings of singularities on $ \CW(m)_b$, but no further degenerations of $ \CW(m)_b$ (this is an open condition in the moduli stack of FM degenerations). Degrees of these sections must satisfy the two requirements above. Hence by the choice of the neighbourhood, their degrees recover the weights of the family 
 \[
 \CW(m)\rightarrow B
 \] in this neighbourhood. Indeed, in this case, they can be determined solely in terms of the conditions (1) and (2). 
 
 Since the weights associated to admissible maps from curves (\ref{weights}) satisfy the  conditions (1) and (2), being given by degrees of vector bundles, we conclude that they are continuous in the sense of Definition \ref{weighted}.
\\

\textit{Step 2}. It remains to show that the weights associated to an admissible map 
\[
f\colon (C,\underline{p}) \rightarrow  W(m)
\]
 satisfy the non-negativity requirements from Definition \ref{moduliwieght}. End components have positive weights that are bounded  by Lemma \ref{weightP} and the defining properties of $\epsilon$-unramification,  Definition \ref{unr}. The fact that ruled components of $W(m)$ have non-negative weights follows from Lemma \ref{ruledpos}. We have to show that the ruled components of all smoothings of $W(m)$ also have non-negative weights. To show  that weights (\ref{weights}) indeed satisfy this property, we choose a non-ruled component 
\[
W(m)' \subset W(m),
\] 
which is also non-end and non-root.  We then choose a chain of components, 
\[
P_{\mathrm{chain}} \subset W(m)
\]
which contains $W(m)'$, an end component, and the root component. We  contract all components but those contained in $P_{\mathrm{chain}}$, 
\[ 
p\colon W(m) \rightarrow P'_{\mathrm{chain}},
\]
where $P'_{\mathrm{chain}}$ is the result of the blow-down of exceptional divisors in $P_{\mathrm{chain}}$. 
 The image of the component $W(m)'$  becomes ruled in  $P'_{\mathrm{chain}}$ after the contraction. The result of such contraction is equivalent to the smoothing of contracted components. Indeed, if a component is smoothed over the spectrum of a discrete valuation ring $\Spec(R)$, then it can also be contracted in the total space of the family. All ruled components of smoothings of $W(m)$ can be  obtained in this way.
 
 We now consider the map from the curve given by the composition with the contraction,
 \[
 p \circ f\colon  (C,\underline{p}) \rightarrow P'_{\mathrm{chain}},
 \]
  which is admissible. By Lemma~\ref{weightP}, the weight of $W(m)'$ after the contraction is equal to the sum of the original weight of $W(m)'$ and the weights of the contracted components. Hence, it is equal to the weight associated with the smoothing of the contracted components by \textit{Step~1}, which is non-negative by Lemma~\ref{ruledpos}. This completes the proof.
\qed


\begin{rmk}\label{weight0} Unlike in the set-up of \cite{YZ}, we allow ruled components of weight~$0$, and non-ruled and non-end  components of negative weight if $\dim(X)>1$. This is because the weights defined in~(\ref{weights}) have these properties. A ruled component $P_{\mathrm{rul}}$ has weight~$0$ if the restriction of the curve to this component, $C|_{P_{\mathrm{rul}}}$, is a union of $\p^1$’s without marked points whose image is a union of fiber lines and at least one non-fiber line, and such that there are no ramifications away from the divisors $D$ and~$E$ (see Section~\ref{FMdeg} for the definition of fiber and non-fiber lines). These components are stable with respect to the stability of unramified maps introduced in~\cite{KKO}, and more generally with respect to $\epsilon$-unramification. See Figure~\ref{fig:unstable} for a pictorial representation of such components.
	
	Note that if $\dim(X)=1$, then the notion of a non-fiber line does not apply. Hence, if  $\p^1$ maps to a rational bridge with no ramifications away from $0$ and~$\infty$, it is unstable.
	
\end{rmk}

\begin{rmk} \label{polar}
	 It is unclear how to construct a universal polarisation for FM degenerations in moduli spaces of $\epsilon$-unramified maps. If $\dim(X)=1$, such polarisation is provided by the branching divisor \cite{N22}. While if $\epsilon=-$, then the target FM degenerations are fixed to be $X$, hence any polarisation of $X$ works. However, it seems that there does not exist a construction that works for an arbitrary dimension and an arbitrary value of $\epsilon$. 
	\end{rmk}


\section{Obstruction theory and invariants} \label{sectioninvariants}

\subsection{Obstruction theory and admissibility}
There are several ways to equip moduli spaces of maps to normal crossing targets with an obstruction theory. The main technical difficulty lies in the fact that the admissibility condition is not open in the space of all maps. Hence, the standard obstruction theory for maps to projective varieties must be modified. The first solution to this problem was proposed by Li \cite{Lideg} in the context of maps to a target with a fixed smooth divisor. Li directly identifies the part of the standard obstruction theory of maps that preserves the admissibility condition.

Other approaches modify the moduli problem in the way that the admissibility becomes an open condition. This is achieved by equipping the targets and sources with auxiliary structures. It can  either be auxiliary logarithmic structures of Kim \cite{KiL} or orbifold structures of Abramovich--Fantechi \cite{AF}. We will follow the approach of Kim; see \cite{AMW} for an overview of different approaches to maps with singular targets and their relations. 

It is also possible to use the logarithmic techniques without introducing auxiliary logarithmic structures, as it was done in \cite{KKO}.  However, the drawback is that it produces only a relative obstruction theory over a stack which is not necessarily smooth but only pure-dimensional, which, for example, complicates the analysis of torus localisation formulas. By the main result of \cite{LM}, the virtual fundamental classes of \cite{KiL} and \cite{KKO} agree after pushing forward the former class with respect to the morphism which forgets the auxiliary logarithmic structures.  

\subsection{Relative $\epsilon$-unramified log maps} \label{auxi}
We define 
\[ \Mbar^\epsilon_{g,n,m}(X, \beta,\underline{\eta})^{\log} \colon (Sch/\BC)^\circ \rightarrow Grpd \]
to be the moduli space of relative $\epsilon$-unramified log maps in the sense of \cite[Section 5.2]{KiL}. More precisely, a family over a base scheme $B$ in  $\Mbar^\epsilon_{g,n,m}(X, \beta,\underline{\eta})^{\log}$ consists of the following data:
\begin{itemize}
	\item A locally-free log structure $N_B$ on $B$.
	\item A minimal log prestable curve 
	\[((\CC,\underline{p}), M_\CC) \rightarrow (B,N_B)\] 
	and an extended log twisted  FM degeneration 
	\[(\CW(m), M_{\CW(m)}) \rightarrow (B,N_B),
	 \] 
	satisfying the corank condition \cite[Section 5.2.2 (3)]{KiL};   for both the source and the target, the logarithmic structures do not include markings or exceptional divisors.
	\item A log map over $(B,N_B)$,
	  \[
	\begin{tikzcd} [row sep=small, column sep = small]
	((\CC,\underline{p}), M_{\CC} )  \arrow{dr} \arrow{rr}{f} &&  (\CW(m), M_{\CW(m)}) \arrow{dl} \\
		& (B,N_B), 
	\end{tikzcd}
	\]
	 whose underlying classical map is relative $\epsilon$-unramified in the sense of Definition \ref{relative}. 
 \end{itemize}
 
Associated to $\Mbar^\epsilon_{g,n,m}(X, \beta,\underline{\eta})^{\log}$, we have the universal curve, the universal FM degeneration, and the universal map, 
\begin{align*}
	&\tau \colon (\FC,\underline{p}) \rightarrow \Mbar^\epsilon_{g,n,m}(X, \beta,\underline{\eta})^{\log}, \\
	& \pi \colon \mathfrak{W}(m) \rightarrow \Mbar^\epsilon_{g,n,m}(X, \beta,\underline{\eta})^{\log}, \\
	& F \colon \FC\rightarrow \mathfrak{W}(m).
\end{align*} 
Let 
\[ 
T^{\mathrm{log}}_{\pi, \mathfrak{E}} =T^{\mathrm{log}}_{ \pi} \Big(-\log\Big(\sum^m_{j=1} \mathfrak{E}_j \Big)\Big) 
\] 
be the $\pi$-relative logarithmic tangent bundle of $\mathfrak{W}(m)$ associated to the universal logarithmic structure on $\mathfrak{W}(m)$ and  $\Mbar^\epsilon_{g,n,m}(X, \beta,\underline{\eta})^{\log}$, and the  exceptional divisors of $\mathfrak{W}(m)$,
\[ 
\mathfrak{E}_j  \subset  \mathfrak{W}(m), \quad j=1, \hdots, m. 
\]

 Let 
 \[
 \FM^{\log}_{X,m}\quad \text{and} \quad \FM^{\log}_{g,n+n'} 
 \]
  be the moduli space of extended log twisted blown-up FM degenerations and   the moduli space of minimal log prestable curves, respectively. The extra markings in  $\FM^{\log}_{g,n+n'}$ will correspond to the points of contact of the curve with the exceptional divisors.  The fibred product 
  \[
 \mathfrak{S}= \FM^{\log}_{X,m} \times_{\log} \FM^{\log}_{g,n+n'}
  \] parametrises a pair of log twisted  FM degenerations and minimal log prestable curves over  the same logarithmic base
    \[
  \begin{tikzcd} [row sep=small, column sep = small]
  	((\CC,\underline{p},\underline{p}_1,\hdots,\underline{p}_m ), M_{\CC} )  \arrow{dr} &&  (\CW(m), M_{\CW(m)}) \arrow{dl} \\
  	& (B,N_B). 
  \end{tikzcd}
  \]
  The stack $\mathfrak{S}$ is formally smooth by the discussion in \cite[Section 6.2.3]{KiL}.  
  
  We also have a forgetful morphism between stacks, 
  \[ 
  \Mbar^\epsilon_{g,n,m}(X, \beta,\underline{\eta})^{\log} \rightarrow   \mathfrak{S},
  \] 
   such that the extra markings
  \[
  \underline{p}_j:=\{p_{j,1},\hdots, p_{j,\ell(\mu^j)}\}=f^{-1}(E_j) 
  \]
  are the points of contact with the exceptional divisors $E_j$.

 
\begin{thm} \label{obst} The space $\Mbar^\epsilon_{g,n,m}(X, \beta, \underline{\eta})^{\log}$ is a proper Deligne--Mumford stack, and there exists a morphism 
	\[
	\phi \colon (R\pi_*F^*T^{\mathrm{log}}_{\pi, \mathfrak{E}})^\vee  \rightarrow \BL_{\Mbar^\epsilon_{g,n,m}(X,\beta,\underline{\eta})^{\log} /\mathfrak{S}},
	\]  which defines a perfect obstruction theory on $\Mbar^\epsilon_{g,n,m}(X,\beta,\underline{\eta})^{\log}$ relative to $\mathfrak{S}$. 
\end{thm}
\textit{Proof.} By \cite[Theorem 6.3.1]{KiL}, the fact that  $\Mbar^\epsilon_{g,n,m}(X, \beta, \underline{\eta})$ is a proper  Deligne--Mumford stack implies that  $\Mbar^\epsilon_{g,n,m}(X, \beta,\underline{\eta})^{\log}$ is also a proper  Deligne--Mumford stack. The former is established in Theorem \ref{thmrel}. 

By \cite[Section 7]{KiL}, the complex $(R\pi_*F^*T^{\mathrm{log}}_{\pi, \mathfrak{E}})^\vee$ indeed defines a perfect obstruction theory. The obstruction-theory morphism $\phi$ constructed in the same way as in for moduli spaces of stable maps to a fixed projective target \cite{BF},  using the properties of the logarithmic tangent bundles in relation to the deformation theory of schemes with logarithmic structures \cite{Ol}.  See also the proof of  \cite[Proposition 5.1.1]{KKO} for more details.  \qed

\subsection{Virtual fundamental class} \label{Virtual}
By Theorem \ref{obst} and \cite{BF}, we have a virtual fundamental class associated to the obstruction theory $\phi$, 
\[
[\Mbar^\epsilon_{g,n,m}(X, \beta,\underline{\eta})^{\log}]^{\mathrm{vir}} \in H_*(\Mbar^\epsilon_{g,n,m}(X, \beta,\underline{\eta})^{\log}).
\]
Using the natural  morphism
\[
p\colon \Mbar^\epsilon_{g,n,m}(X, \beta,\underline{\eta})^{\log} \rightarrow \Mbar^\epsilon_{g,n,m}(X, \beta,\underline{\eta}),
\] 
which forgets the auxiliary logarithmic structures, we define the virtual fundamental class of $\Mbar^\epsilon_{g,n,m}(X, \beta,\underline{\eta})$, 
\[
[
\Mbar^\epsilon_{g,n,m}(X, \beta,\underline{\eta})]^{\mathrm{vir}}:=p_*[\Mbar^\epsilon_{g,n,m}(X, \beta,\underline{\eta})^{\log}]^{\mathrm{vir}}.
\]
The morphism $p$ is a virtual normalisation in the sense of \cite{LM}. In the context of admissible covers, such a normalisation was considered in \cite{ACV} using orbifold curves; in this case, $p$ is a normalisation in the usual sense of the term. In \cite{AF}, orbifold structures were also used in the case of higher-dimensional targets in place of logarithmic structures. Moreover, by the main result of \cite{LM}, the virtual fundamental class $
[\Mbar^\epsilon_{g,n,m}(X,\beta,\underline{\eta})]^{\mathrm{vir}}$
defined above agrees with the one constructed using the arguments of \cite[Section 5]{KKO}.

\subsection{Relative tangent line bundles} \label{reltangentbundle} For each exceptional divisor $E_j$, a blown-up FM degeneration $ W(m)$ admits a contraction to a blow-up of $X$ at a point, 
\[  W(m) \rightarrow \mathrm{Bl}_{x_j}(X), \quad j=1,\dots, m, \]
where $x_j$ is the image of the exceptional divisor  $E_j \subset  W(m)$ after the contraction of all components except the root component $X$. The exceptional divisor $E_j$ is mapped isomorphically onto the exceptional divisor of $\mathrm{Bl}_{x_j}(X)$.   Let
\[\pi \colon \mathfrak{W}(m)\rightarrow \Mbar^\epsilon_{g,n,m}(X,\beta,\underline{\eta})\]
be the universal  blown-up FM degeneration, and consider the blow-up of the diagonal, 
\[q \colon  \mathrm{Bl}_\Delta (X\times X) \rightarrow X, \]
which we treat as the universal blow-up up of $X$ at a point, since for a point $x \in X$, we have $q^{-1}(x)=\mathrm{Bl}_x(X)$. 
By  \cite[Section 9]{NHilb},  the contraction  
\[ W(m) \rightarrow  \mathrm{Bl}_{x_i}(X) \hookrightarrow \mathrm{Bl}_\Delta (X\times X)\] can be defined in families, giving rise to the universal contraction
\[ \tau_j  \colon \mathfrak{W}(m) \rightarrow \mathrm{Bl}_\Delta (X\times X). \]
Consider now the universal inclusion of the exceptional divisor, 
\[\iota_j \colon  \mathfrak{E}_j\hookrightarrow \mathfrak{W}(m),\]
with the associated relative normal bundle $\CN_{\mathfrak{E}_i}$. 
We define the associated tangent line bundle as follows,
\begin{equation} \label{relativepsi}
\widetilde{\BL}(E_j)^\vee:= \pi_{*}(\CN_{\mathfrak{E}_i} \otimes (\tau_j \circ \iota_j)^* \CO_{\p (TX)}(1)).
\end{equation}
Note that $\widetilde{\BL}(E_j)^\vee$ is a line bundle, because $\CN_{E_j}$ is isomorphic to $\CO_{\p (T_x X)}(-1)$ over a closed point in  $\Mbar^\epsilon_{g,n,m}(X,\beta,\underline{\eta})$. However, we will be primarily interested in a certain modification of the line bundle $\widetilde{\BL}(E_j)$. Let 
\[\mathfrak{Z}_{0,j} \subset \FM_{X,m,d} \] 
be the reduced closed substack parametrizing FM degenerations, such that the $j$-th exceptional divisor lies on a marked ruled component (i.e.,\ a ruled component with a one-vertex edge), and this component is of weight 0 with respect to the weights from Lemma \ref{lemmaweights}. See Remark \ref{weight0} for how such components arise.   The substack $\mathfrak{Z}_{0,j}$ is a divisor. We then twist $\widetilde{\BL}(E_j)$ by $\CO(-\mathfrak{Z}_{0,j})$, 
\begin{equation} \label{twist}
\BL(E_j):=\widetilde{\BL}(E_j)\otimes \CO(-\mathfrak{Z}_{0,j}). 
\end{equation}
 
The reason we need these twisted line bundles is because they appear in the analysis of the localisation formula, Proposition \ref{isomophism}; see the end of Section \ref{sectionend} for an explanation. We also suggest to consult \cite[Section 9]{NHilb} for more about cotangent line bundles in the context of FM degenerations. 
\subsection{Invariants} \label{Sectioninv}
There exist the usual structures needed to define Gromov--Witten type invariants:
\begin{itemize}
	\item Virtual fundamental classes,
	\[[\Mbar^\epsilon_{g,n,m}(X, \beta, \underline{\eta})]^{\mathrm{vir}}\in H_{*}(\Mbar^\epsilon_{g,n,m}(X, \beta, \underline{\eta})), \]
	which are constructed in Section \ref{Virtual}.
	\item Evaluation morphisms at markings,
	\[\ev_i: \Mbar^\epsilon_{g,n,m}(X, \beta, \underline{\eta}) \rightarrow X, \quad f \mapsto (p\circ f)(p_i),\]
	where $p\colon W \rightarrow X$ is the contraction map associated to a FM degeneration. 
	 Exceptionally for $\epsilon=+$, we obtain  evaluation morphisms to $\p(TX)$ by the construction from \cite[Lemma 3.2.4]{KKO},
	 \[ \p (\ev_i): \Mbar^+_{g,n,m}(X, \beta,\underline{\eta}) \rightarrow \p(TX). \]
	 \item Classes associated to the cotangent line bundles, defined by $L_{i}|_{(C,\underline{p})}=T^*_{p_i}C$,
	 \[ \psi_i:=\mathrm{c}_1(L_i). \]
\end{itemize}	 
 
There also exist relative analogues of these structures. First, given a partition $\eta=(\eta_1,\dots ,\eta_{\ell(\eta)})$ of arbitrary size $|\eta|=\sum_i \eta_i$, we  define 
\vspace{0.1cm}
	\begin{align}  \label{symm}
		\begin{split}
	\CI_\eta\Sym_X\p(TX)&:= 	 \underbrace {\p(TX)  \times_X \dots \times_X \p(TX)}_{\ell(\eta) \text{ times}}, \\
\CI\Sym_X\p(TX)&:= \coprod_{\eta}  \CI_\eta\Sym_X\p(TX), 
\end{split}
\end{align}
where the disjoint union is taken over all  partitions $\eta$.  We then have the following additional structures: 
\begin{itemize}
	\item   Relative evaluation morphisms associated to the exceptional divisors $E_j$,
	\[\widetilde{\ev}_j: \Mbar^\epsilon_{g,n,m}(X, \beta,\underline{\eta}) \rightarrow \ \CI_{\eta^j}\Sym_X\p(TX), \quad f \mapsto f(C)\cap E_j.\]
	\item Relative $\psi$-classes associated to the cotangent line bundles from (\ref{relativepsi}), 
	\[ \Psi_j=\mathrm{c}_1(\BL(E_j)). \]	
	\end{itemize}
 
 Using the structures above, we are now ready to defined the invariants associated to the spaces  $\Mbar^\epsilon_{g,n,m}(X, \beta, \underline{\eta})$. 
  
\begin{defn} \label{inv} Given classes $\gamma_i \in H^*(X)$ and $\gamma'_j \in H^*(  \CI\Sym_X\p(TX))$, we define relative $\epsilon$-unramified Gromov--Witten invariants of $X$, 
	\begin{multline*} \langle \gamma_1 \psi_1^{k_1} \cdot \! \cdot \! \cdot  \gamma_n \psi_n^{k_n} \mid  \gamma'_1  \Psi_1^{k'_1} \cdot \! \cdot \! \cdot  \gamma'_m \Psi_m^{k'_m}  \rangle^\epsilon_{g,\beta}:= \\
		 \frac{1}{\prod^m_{j=1} |\Aut(\eta^j)|}\int_{[\Mbar^\epsilon_{g,n,m}(X, \beta, \underline{\eta})]^{\mathrm{vir}}}\prod^n_{i=1} \ev_i^*(\gamma_i)  \psi_i^{k_i} \cdot \prod^m_{j=1}   \widetilde{\ev}_{j}^*(\gamma'_j)  \Psi_j^{k'_j},
		\end{multline*}
	such that $\Aut(\eta):= \prod_i S_{m_i}$, where $m_i$ are the multiplicities of the parts of the partition $\eta$. 
	
Exceptionally for $\epsilon=+$,  we allow to pull back classes $\gamma_i \in H^*(\p (TX))$, using the evaluation morphisms $\p(\ev_i)$. The ramification profiles $\eta^j$ are implicitly specified by the connected components of  $\CI\Sym_X\p(TX)$  on which the classes $\gamma'_j$ are supported.  Observe that these invariants specialise to Gromov--Witten  and unramified Gromov--Witten invariants for the associated values of $\epsilon$ if $m=0$. 
	\end{defn}

\section{Master spaces} \label{masterspaces}

\subsection{Overview of master spaces}Having defined the spaces $\Mbar^\epsilon_{g,n,m}(X, \beta,\underline{\eta})$, we are now ready to address the wall-crossing phenomenon associated with the variation of the parameter $\epsilon \in \BR_{>0}$. In this section, we construct master spaces using the theory of entangled end components introduced by Zhou in the one-dimensional case~\cite{YZ}. The wall-crossing formula, Theorem~\ref{wallcrossing}, will be established with the help of torus localisation on the master spaces. It provides an explicit comparison of $\epsilon$-unramified invariants for values of $\epsilon$ separated by a single wall. By applying Theorem~\ref{wallcrossing} inductively, we obtain a formula comparing invariants for extremal values of $\epsilon$, corresponding to the moduli spaces of stable and unramified maps, which will be studied in the case of threefolds in Section~\ref{sectionthreefolds}.

\subsection{Walls} For fixed $\beta$, $g$, $n$, and $\underline{\eta}$,  the space of $\epsilon$-unramified stabilities, introduced in Definition \ref{unr},  is the real half-line $\BR_{>0}$. By the boundedness of the space of admissible maps to FM degenerations with end components of positive weight, Lemma \ref{alg}, there exist finitely many values\footnote{Indeed, by the boundedness of this moduli space, weights of points and end components can take only finitely many values, since they are half-integers. } 
\[
\epsilon_0 \in \BR_{>0},
\] such that the moduli space  $\Mbar^\epsilon_{g,n,m}(X, \beta, \underline{\eta})$ changes abruptly, as $\epsilon$ crosses $\epsilon_0$.  We call such values of $\epsilon$ \textit{walls}. By the definition of $\epsilon$-unramification, the reciprocals of walls are half-integers, 
\[
d_0:=1/\epsilon_0 \in \frac{1}{2}\BZ.
\] 
We call open intervals in $\BR_{>0}$ between walls \textit{chambers}.  For a fixed wall $\epsilon_0 \in  \BR_{>0}$, let 
\[
\epsilon_-, \epsilon_+ \in \BR_{>0}
\]
 be some real numbers, such that
 \[
\frac{\epsilon_0}{1+\epsilon_0/2}<  \epsilon_- < \epsilon_0< \epsilon_+ < \frac{\epsilon_0}{1-\epsilon_0/2}.
 \]
 In other words, $\epsilon_-$ and $\epsilon_+$ are values in two adjacent chambers separated by the wall $\epsilon_0$. 
 

\subsection{Entanglement of Zhou} \label{ent} We start with the definition of the space of FM degenerations with entangled end components. It is a higher-dimensional analogue of the space of curves with entangled rational tails defined in \cite[Section 2.2]{YZ}.  Zhou's entanglement is given by an inductive blow-up of the boundary strata of the moduli space of weighted curves, which aims to remove stabilisers of curves with rational tails.   Various results from \cite[Section 2]{YZ} extend to our case, since the proofs rely just on the fact that we have a pair 
\[
\mathfrak{Z}\subset \FM,
\] such that $ \mathfrak{Z}$ is a  normal crossing divisor and $\FM$ is a smooth algebraic stack.

 Consider the moduli space of weighted blown-up FM degenerations $\FM_{X,m,d}$ from Definition \ref{moduliwieght}, such that $d_0=1/\epsilon_0$ for a wall $\epsilon_0 \in \BR_{>0}$. Let \[\mathfrak{Z}_k \subset \FM_{X,m,d}\] be the reduced closed substack which parametrizes weighted FM degenerations with at least $k$ end components of weight $d_0$. By the \'etale projection (\ref{weightsprojection}), the space 
 \[
 \FM_{X,m,d}
 \] 
 admits a smooth atlas, whose charts are given by Zariski open subsets of Fulton--MacPherson compactifications  
 \[ 
 U \subset  X[n'+m]. 
 \]
In such a chart, the locus of singular marked FM degenerations is a normal crossing divisor by \cite[Theorem 3]{FM}. Hence, $\mathfrak{Z}_1\subset \FM_{X,m,d}$ is also a normal crossing divisor. Moreover, in formal neighbourhoods,  $\mathfrak{Z}_k$ is the stratum of $\mathfrak{Z}_1$ given by the intersection  of $k$ local hyperplane branches of $\mathfrak{Z}_1$. 
 
  Let $h$ be the maximum number of weight $d_0$ end components of FM degenerations in $\FM_{X,m,d}$, then $\mathfrak{Z}_h$ is the deepest stratum of $\mathfrak{Z}_1$, and therefore is smooth.   Set 
  \[\mathfrak{U}_h:=\FM_{X,m,d}.\] Let 
\[\mathfrak{U}_{h-1} \rightarrow \mathfrak{U}_h\]
be the blow-up of $\mathfrak{U}_h$ at $\mathfrak{Z}_h$. For $i\in \{1,\dots, h-1\}$, we define $\mathfrak{U}_{i-1}$ to be the blow-up of $\mathfrak{U}_{i}$ at the proper transform of $\mathfrak{Z}_i$ in $\mathfrak{U}_{i}$,  denoted by $\mathfrak{Z}_{(i)}$, which is smooth.  We then set 
\[\widetilde{\FM}_{d,m}^X:=\mathfrak{U}_0,\]
this is the space of FM degenerations with \textit{entangled} end components, its purpose is to remove the $\BC^*$-scaling automorphisms of end components. 
\begin{defn} 
Given a FM degeneration $ [W(m)] \in \FM_{X,m,d}(\BC)$, an \textit{entanglement} of $ W(m)$ is a point in the fiber of 
\[
\widetilde{\FM}_{X,m,d} \rightarrow \FM_{X,m,d}
\] over $[W(m)]$. We denote an entanglement of $ W(m)$ by $(W(m), e)$, where $e$ is a closed fiber of the projection above. 
\end{defn}

 Let 
 \[
 \{P_1, \dots, P_\ell\}
 \] be the set of the end components of $ W(m)$ of weight $d_0$. At 
 \[
 [W(m)] \in  \FM_{X,m,d}(\BC),
 \]  the locus $\mathfrak{Z}_1$ has local branches $\mathfrak{H}_j$,  where an end component $P_j$ remains intact  (i.e., it does not smooth out). Define
\[k:= \mathrm{min}\{ i \mid \text{the image of } ( W(m),e) \text{ in } \mathfrak{U}_i \text{ lies in } \mathfrak{Z}_{(i)}   \}.\]
 The image of $( W(m),e)$ in $\CU_k$ lies on the intersection of $k$ branches $\{\mathfrak{H}_{j_1}, \dots, \mathfrak{H}_{j_k}\}$. The corresponding set of end  components $\{P_{j_1}, \dots, P_{j_k}\}$ are called \textit{entangled.}

\begin{defn} We define the divisor
	\[\mathfrak{Y}_i \subset \widetilde{\FM}_{X,m,d}\]
	to be the closure of the locus of FM degenerations with exactly $i+1$ entangled end components. Equivalently, $\mathfrak{Y}_i$ is the proper transform of the exceptional divisor of the blow-up $\mathfrak{U}_{i} \rightarrow \mathfrak{U}_{i-1}$. 
\end{defn}

\begin{defn} \label{defncalibration} We define the \textit{calibration bundle} on $\widetilde{\FM}_{X,m,d}$, 
	\[ \BM_{\widetilde{\FM}_{X,m,d}}:=\CO(-\mathfrak{Z}_{1}).  \]
Let 
\[M\widetilde{\FM}_{X,m,d}:=\p(\BM_{\widetilde{\FM}_{X,m,d}} \oplus \CO_{\widetilde{\FM}_{X,m,d}})\]
be the moduli space of FM degenerations with \textit{calibrated} end components.
	\end{defn}
A $B$-valued point of $M\widetilde{\FM}_{X,m,d}$ is given by 
\[ ( \CW(m),e,\CL,v_1,v_2),\]
such that:
\begin{itemize}
\item $( \CW(m),e) \in \widetilde{\FM}_{X,m,d}(B)$ is a FM degeneration with entanglement,
\item  $\CL$ is a line bundle on $B$, such that $v_1 \in H^0(B, \BM_B\otimes \CL)$ and $v_2 \in H^0(B,  \CL)$ are sections with no common zeros.
\end{itemize} 
We refer to the data above as a \textit{calibration} of $\CW(m)$. 

\subsection{Example of entanglement} \label{exampleenta}
 In a smooth chart, a formal local description of the pair $(\FM_{X,m,d}, \mathfrak{Z}_1)$  for $h=2$ is provided by 
\[
(\BA^2, Z_1),
\] 
where $Z_1=\{xy=0\}$. The blow-up procedure described above tells us to blow up $\BA^2$ at 
\[
Z_2=\{0\} \subset \BA^2.
\] The proper transform $\widetilde{Z}_1$ of $Z_1$ is given by two disjoint lines intersecting the exceptional divisor $E=\p^1$ at two distinct points. The entanglement of a point $w \in Z_1 \subset \BA^2$ is a point $\widetilde{w} \in \mathrm{Bl}_0 (\BA^2)$ lying over it. It has two entangled end components if 
\[\widetilde{w} \in \BA^1 \setminus \{0\}=E \setminus (\widetilde{Z}_1\cap E);
\] otherwise, it has one entangled end component or no end components at all. The divisor $Y_1$ is simply the exceptional divisor $E$ of $\mathrm{Bl}_0 (\BA^2)$. 

Taking the stack quotient of $\BA^2$ by $(\BC^*)^2$ corresponds to introducing $\BC^*$-scaling automorphisms of end components in  $(\FM_{X,m,d}, \mathfrak{Z}_1)$. The locus of two entangled end components is then given by 
\[ 
[(\BA^1 \setminus \{0\})/(\BC^*)^2] \subset [\p^1/(\BC^*)^2],
\]
which is isomorphic to $B\BC^*$. In particular, the result of the entaglement is the reduction of the stabiliser group of the locus, parametrising FM degenerations  with two end components, from $(\BC^*)^2$ to  $\BC^*$. The remaining stabiliser can be dealt with by means  of introducing non-zero cotangent vectors at the points of $Z_1$, or, in other words, a section of the pullback of $\CO(-Z_1)$,  which Zhou calls \textit{calibration}.

\subsection{Master space}

By allowing the strict inequality of Definition \ref{epsilon} (1) to become a non-strict one, we can construct a non-separated space containing both 
  $\Mbar^{\epsilon_-} _{g,n,m}(X, \beta,\underline{\eta})$ and $\Mbar^{\epsilon_+} _{g,n,m}(X, \beta,\underline{\eta})$, which is the space of $\epsilon_0$-semi-unramified maps $\FM^{\epsilon_0}_{g,n,m}(X,\beta,\underline{\eta})$, defined below. If this space was separated, we could use its intersection theory to derive a relation between $\epsilon$-unramified invariants associated to $\epsilon_-$ and $\epsilon_+$. However, it is not.   To make this space into a proper Deligne--Mumford stack, called master space, Zhou's entaglement and calibration is required. This leads to the notation of $\epsilon_0$-unramified maps introduced in Definition \ref{epsilon0unr}.   Let us start with the definition of $\epsilon_0$-semi-unramified maps.
  
 \begin{defn} \label{relative2} Let $\epsilon_0 \in  \BR_{>0}$ be a wall. A map  $f \colon (C,\underline{p}) \rightarrow  W(m)$ is said to be $\epsilon_0$-$semi$-$unramified$ if 
	\begin{itemize} 
		\item[(1)] for all end components $P \subset  W(m)$, $w(P)\geq 1/\epsilon_0$,
		\item[(2)] for all $x\in  W(m)^{\mathrm{sm}}\setminus \cup_j E_j$, $w(x) \leq 1 /\epsilon_0$,
		\item[(3)] the map $f$ is admissible over $ W(m)^{\mathrm{sing}}$ and $\cup_j E_j$,
		\item[(4)] $|\Aut_C(f)| < \infty$, and $|\Aut_{C,P_{\mathrm{rul}}}(f|_{P_{\mathrm{rul}}})|<\infty$ for all ruled components $P_{\mathrm{rul}} \subset  W(m)$; all connected components of the source map non-trivially. 
		
	\end{itemize}

\end{defn}

\begin{defn} Given an admissible map $f \colon (C,\underline{p}) \rightarrow  W(m)$.  An end component $P \subset  W(m)$ is \textit{constant} if there is only one base point on $P$ and the image of $f|_{P}$ is a union of lines. 
	
	 An end component $P \subset  W(m)$ is \textit{strictly constant} if there is only one base point on $P$ and the image of $f|_{P}$ is a union of lines passing through the base point, or, equivalently, $\Aut_{C,P}(f|_{P})=\BC^*$. 
\end{defn}

Observe that the difference between Definition~\ref{epsilon} and Definition~\ref{relative2} is that the condition on the weights for end components is given by a non-strict inequality, and that we require finiteness of automorphisms only for ruled components. In particular, end components might have infinitely many automorphisms; that is, they might be strictly constant.

If $\dim(X) \geq 2$, we must distinguish between strictly constant and constant end components. Although maps over constant end components do not necessarily have infinitely many automorphisms, they are nevertheless excluded from the moduli space of $\epsilon$-unramified maps by Definition~\ref{epsilon} (1) and Lemma~\ref{less}; see also Definition~\ref{unr} (6).
  Constant end components will be essential for the properness of the master space. 

\begin{defn}Let  
\[ \FM^{\epsilon_0}_{g,n,m}(X,\beta,\underline{\eta}) \colon (Sch/\BC)^\circ \rightarrow Grpd\] 
 be the moduli space of $\epsilon_0$-semi-unramified maps of degree $\beta$ from possibly disconnected curves of genus $g$ with $n$ markings and with $m$ exceptional divisors on the target.
 \end{defn}

 The moduli space $\FM^{\epsilon_0}_{g,n,m}(X,\beta,\underline{\eta})$ is a quasi-separated algebraic stack of finite type by the arguments from Lemma \ref{alg} and Lemma \ref{open}. Next, we will introduce entanglement and calibration by taking the fibred product, 
\[M\FM^{\epsilon_0}_{g,n,m}(X,\beta,\underline{\eta}):=\FM^{\epsilon_0}_{g,n,m}(X,\beta,\underline{\eta})\times_{\FM_{X,m,d}} M\widetilde{\FM}_{X,m,d}. \]
A $B$-valued point of this moduli space is
\[\xi=(\theta,\lambda)=(\CC,  \CW(m), \underline{p},f, e, \CL, v_{1}, v_{2}),\]
where $\theta$ is the underlying $\epsilon_0$-semi-unramified map, and $\lambda$ is a calibration of $\CW(m)$, 
\[\theta=(\CC, \CW(m), \underline{p},f) \ \text{and} \ \lambda=(e, \CL, v_1,v_2),\]
 the meaning of gadgets in $\lambda$ is the same as in discussion after Definition \ref{defncalibration}.

We are now ready to define the master space. In short,  we put $\epsilon_-$- and $\epsilon_+$-unramified maps at different ends of the calibration bundle, while the stabilisers of maps with constant end components are removed by means of the entanglement.  Both aspects of the definition are needed for the properness of the master space; see the beginning of Section \ref{sectionprop} for more details about the properness. 

\begin{defn} \label{epsilon0unr} A family of $\epsilon_0$-semi-unramified maps with calibrated end components over a base scheme $B$,
\[(\CC,  \CW(m), \underline{p},f, e, \CL, v_{1}, v_{2}),\]
is $\epsilon_0$-unramified if 
\begin{itemize}
	\item[1)] all constant end component of geometric fibers over $B$ are entangled,
	\item[2)] if a geometric fiber $ \CW(m)_{b}$ has end components  of weight $d_0$, then end components of $\CW(m)_b$ contain all points of weight $d_0$,
	\item[3)] if $v_1(b)=0$, then $(\CC,  \CW(m),\underline{p},f)_b$ is $\epsilon_+$-unramified,
	\item[4)] if $v_2(b)=0$, then $(\CC,  \CW(m), \underline{p},f)_b$ is $\epsilon_-$-unramified. 
\end{itemize} 
\end{defn}
\begin{defn} \label{masterspacedefn} The \textit{master space},
\[M \Mbar^{\epsilon_0}_{g,n,m}(X,\beta,\underline{\eta}) \subseteq M\FM^{\epsilon_0}_{g,n,m}(X,\beta,\underline{\eta}),\]
is  the moduli space of  $\epsilon_0$-unramified maps with calibrated end components of degree $\beta$ from possibly disconnected curves of genus $g$ with $n$ markings and with $m$ exceptional divisors on the target.
\end{defn}
\begin{prop} The moduli space $M \Mbar^{\epsilon_0}_{g,n,m}(X,\beta,\underline{\eta})$ is a quasi-separated Deligne--Mumford stack of finite type with a perfect obstruction theory.
	\end{prop}  
\textit{Proof.}   $M \Mbar^{\epsilon_0}_{g,n,m}(X,\beta,\underline{\eta})$  is an open substack of  $M\FM^{\epsilon_0}_{g,n,m}(X,\beta,\underline{\eta})$ by \cite[Lemma 4.1.4]{YZ}.  Since the latter stack is algebraic, quasi-separated, and of finite type, so is the space $M\Mbar^{\epsilon_0}_{g,n,m}(X,\beta,\underline{\eta})$.  Automorphism groups of objects in  the space $M \Mbar^{\epsilon_0}_{g,n,m}(X,\beta,\underline{\eta})$ are finite by \cite[Lemma 4.1.10]{YZ}. Indeed, the only potential source of infinite automorphisms is the strictly constant end components; their automorphisms are removed by the entanglement and calibration.   Hence, the stack $M \Mbar^{\epsilon_0}_{g,n,m}(X,\beta,\underline{\eta})$ is  Deligne--Mumford by \cite[Section 0DSL]{stacks-project}. \qed


\subsection{Properness of  master spaces} \label{sectionprop}
We now deal with the properness of the master space $M \Mbar^{\epsilon_0}_{g,n,m}(X,\beta,\underline{\eta})$. We will follow the strategy of \cite[Section 5]{YZ}.   Given a discrete valuation ring $R$ with the fraction field $K$, let 
\[\xi^*=(\CC^*,  \CW(m)^*, \underline{p}^*,f^*, e^*, \CL^*, v_{1}^*, v_{2}^*) \in M \Mbar^{\epsilon_0}_{g,n,m}(X,\beta,\underline{\eta})(K)  \] be a family of $\epsilon_0$-unramified map with calibrated tails over $\Spec(K)$.  Let 
\[\theta^*=(\CC^*,  \CW(m)^*, \underline{p}^*,f^*) \ \text{and} \ \lambda=(e^*, \CL^*, v_1^*,v_2^*)\]
be the underlying $\epsilon_0$-semi-unramified map and the calibration of $\CW(m)^*$, respectively. 

Let us start by giving a brief overview of the proof. First, by the properness of the projection
\begin{equation} \label{projectionm}
	M\widetilde{\FM}_{X,m,d} \rightarrow \FM_{X,m,d},
\end{equation}
an extension of the calibration $\lambda^*$ to $\Spec(R)$ is determined by the extension of the family $\CW(m)^*$. The latter also determines the extension of the rest of $\theta^*$, by the properness of the moduli space of stable maps to a projective target. Hence, the proof will consist of two related but distinct tasks:
\begin{itemize}
	\item[(1)] classifying all extensions $\CW(m)$ of $\CW(m)^*$, such that the resulting extension of $\theta^*$ is $\epsilon_0$-semi-unramified,
	\item[(2)] selecting an extension $\theta$, such that $(\theta,\lambda)$ is $\epsilon_0$-unramified. 
\end{itemize}

For the task (1),  there is not much we can do to the central fiber of $\CW(m)$ for $\theta$ to be $\epsilon_0$-semi-unramified. As we show in Lemma \ref{neq01d} and Lemma \ref{eqd0}, we can change only the singularity type of weight $d_0$ end components in the total space of the family by blowing them up and down (or move a weight $d_0$ point onto an end component).  Hence, the classification of such extensions will depend on an array of rational numbers that record the singularity types of end components and the degree of the base change of $\Spec(R)$. We differentiate two cases: (i) the generic fiber $\theta^*$ does not have end components of weight $d_0$, and (ii) the generic fiber $\theta^*$ has end components of weight $d_0$. They are treated in Section \ref{Firstcase} and Section \ref{Secondcase}, respectively. Note that some extensions can have constant end components, and they are essential for the properness of the master space; in fact, they give rise to the wall-crossing components in Section \ref{secfixed}. 

For the task (2), we will use  \cite[Lemma 2.9.1]{YZ}, \cite[Lemma 2.9.2]{YZ}, and \cite[Lemma 2.10.1]{YZ}, which describe  the effect of modifying the central fiber of $\CW(m)$ on the calibration $\lambda$. This is explained in the proof of Theorem \ref{masterprop}. 

In the case of the local example considered in Section \ref{exampleenta}, modifying the central fiber of $\CW(m)$ corresponds to passing from 
\[
(w_1(t),w_2(t))\colon \mathrm{Spf}(\BC[\![t]\!]) \rightarrow  [\BA^2/(\BC^*)^2]
\] 
to 
\[
(t^{k_1}\cdot w_1(t), t^{k_2}\cdot w_2(t)).
\] The result of this modification on the entanglement and the calibration is clear from this perspective. For example, if we want to modify the family in the way that it becomes entangled, we need to find $k_1$ and $k_2$, such that the order of vanishing of the resulting functions above are equal; this is the content of \cite[Lemma 2.9.1]{YZ} and \cite[Lemma 2.9.2]{YZ}.  On the other hand, such a modification also affects the pullback of $\mathcal{O}(-Z_1)$ -- it acquires a twist by the $(k_1+k_2)$-multiple of the closed point of a discrete valuation ring; this is the content of \cite[Lemma 2.10.1]{YZ}.

The difference between the set-up of \cite{YZ} and ours is the existence of ruled components of weight 0, whose origin is explained in Remark \ref{weight0}.  In practice, this means we can apply modifications to whole chains of ruled components connected to end components, while in \cite{YZ} only end components (rational tails) appear. In both cases, the local description of modifications is the same. 
\subsubsection{No weight $d_0$ end components at the generic fiber} \label{Firstcase} Assume $\theta^*$ does not have end components of weight $d_0$.  By the arguments from  Theorem \ref{epsilonthm}, we can always extend the family $\theta^*$ to a family over $\Spec (R)$, possibly after a finite base change.   Let 
\[\theta_+:=(\CC_+,  \CW(m)_+, \underline{p}_+,f_+) \in  \FM^{\epsilon_0}_{g,n,m}(X,\beta,\underline{\eta})(R)\]
be an $\epsilon_0$-semi-unramified extension of $\theta^*$ to $\Spec(R)$, such that:
\begin{itemize}
\item the central fiber is $\epsilon_+$-unramified at all end components, i.e., Definition \ref{epsilon} (1) is satisifed for all end components for $\epsilon=\epsilon_+$,
\item the central fiber is $\epsilon_+$-unramified at all points, i.e., Definition \ref{epsilon} (2) is satisifed for all points for $\epsilon=\epsilon_+$, unless a point is the limit of a weight $d_0$ base point from the generic fiber. 
\end{itemize}	
  More explicitly, the family $\theta_+$ is constructed by firstly applying \textit{Step 1} of Theorem \ref{epsilonthm}, which produces a family whose weights at points of the central fiber are equal to weights of points at the generic fiber. We then apply \textit{Step 2} of Theorem \ref{epsilonthm} with respect to $\epsilon_+$-unramification for all end components. The family $\theta_+$ is  closest to being an $\epsilon_+$-unramified extension of $\theta^*$, except that it might have weight $d_0$ points which come from the generic fiber. 
  
  The calibration $\lambda_+$ is given by the unique
extension of $\lambda^*$ to  by the properness of (\ref{projectionm}). Overall, we obtain 
\[\xi_+=(\theta_+, \lambda_+)\in  M\FM^{\epsilon_0}_{g,n,m}(X,\beta,\underline{\eta})(R).\]

We similarly define 
\[\xi_-=(\theta_-, \lambda_-) \in M\FM^{\epsilon_0}_{g,n,m}(X,\beta,\underline{\eta})(R)\]
to be the extension of $\xi^*$ to $\Spec(R)$, such that the family $\theta_-$ is the unique  $\epsilon_-$-unramified extension of $\theta^*$, and is given by the contraction of all weight $d_0$ end components at the central fiber of $\theta_+$. Again, the calibration is given by the unique extension of $\lambda^*$ by the properness of (\ref{projectionm}). 

  Between $\xi_+$ and $\xi_-$, there exist intermediate extensions of $\xi^*$ that we will classify. By the uniqueness of extensions of calibrations, it is sufficient to classify extensions of $\theta^*$.
  
   \begin{defn} \label{modifiation}
	Let $R$ be a discrete valuation ring, and $(\CC,  \CW(m), \underline{p},f)$ be a family of $\epsilon_0$-semi-unramified maps
	 over $\Spec(R)$. A \textit{modification} of $(\CC,  \CW(m), \underline{p},f)$ is a family of $\epsilon_0$-semi-unramified maps
	$(\CC',  \CW(m)', \underline{p}', f')$  over $\Spec (R')$, such that
	\[(\CC',  \CW(m)', \underline{p}', f') |_{\Spec(K')} \cong (\CC,  \CW(m), \underline{p},f)|_{\Spec(K')},\]
	where $R'$ is a finite extension of $R$ with the fraction field $K'$. 
	
\end{defn}

 Consider  the ordered set consisting of weight $d_0$ end components  and weight $d_0$ base points,
 \[\{P_1 ,\dots, P_h, x_{h+1}, \dots, x_\ell\},\]
  of the central fiber 
 \[
 W(m)_{+} \subset \CW(m)_{+}
 \] of the FM degeneration in $\theta_+$.
We define a vector of $\ell$ numbers associated to the family $\theta_+$,
\[\underline{b}=(b_{1},\dots,b_{\ell})\in (\BQ\cup \{\infty\})^\ell, \]
as follows. Set $b_{i}$
to be $\infty$ for a degree $d_0$ base point $x_i$. The definition of $b_{i}$ for an end component $P_i$ is more involved. First,  an end component $P_i$ might be followed by a chain of ruled components of weight 0,
\[ P_{i,1} \cup \hdots \cup P_{i,k_i},\]
such that $P_i$ is attached to $P_{i,k_i}$, and we set $P_{i,k_i+1}:=P_i$.  See Remark \ref{weight0} for an explanation of how 
 weight 0 ruled components arise.\footnote{Such components do not exist in dimension one.} Let 
 \[
 D_{i,j} \subset P_{i,j}, \quad j=1,\dots, k_i+1,
 \] be the divisors at infinity of the ruled and end components.  By the local description of the universal family of FM degenerations \cite[Section 2.3]{KKO}, in the formal neighbourhood around points of $D_{i,j}$, the family  $\CW(m)_{+}$ looks like 
\begin{equation} \label{an}
	\mathrm{Spf} \left(\frac{\BC[\![z_1,\dots ,z_{\dim+1}, t ]\!]}{(z_1z_2-t^{b_{i,j}})} \right).
	\end{equation}
We will refer to singular locus of (\ref{an}) as the $A_{b_{i,j}-1}$-singularity.     For an end component $P_{i}$, we define 
\begin{equation} \label{bi}
	b_{i}= \sum^{k_i+1}_{j=1}b_{i,j}.  
\end{equation}
The quantity $b_i$ is the order of contact of $\Spec(R)$ with the local branch $\mathfrak{H}_i$ of $\mathfrak{Z}_1$ associated to the component $P_i$. 

Let 
\[
\eta^{i,j}=(\eta^{i,j}_{1}, \hdots, \eta^{i,j}_{\ell(\eta^{i,j})})
\] be the ramification profile of $f_+$ over the divisor  $D_{i,j}$, and 
\[
m_{i,j}:=\mathrm{lcm}(\eta^{i,j}_{1}, \hdots, \eta^{i,j}_{\ell(\eta^{i,j})}).
\] By the admissibility over the singular locus, we have 
\[
m_{i,j} \mid b_{i,j}. 
\]
Indeed, the total space of the curve around a node corresponding to the ramification of order $\eta^{i,j}_s$ over $D_{i,j}$ looks like $\mathrm{Spf}(\BC[\![ x_1,x_2,t]\!]/(x_1x_2-t^{b_{ij}/\eta^{i,j}_s }$)).

 All  modifications of $\theta_+$ (and therefore of $\theta_-$) are given by two types of operations: finite base changes, and blow-ups and contractions around $D_{i,j}$ and $x_i$. The result of these modifications will be a change of singularity type of $\theta_+$ around $D_{i,j}$ and $x_i$. Hence, the classification will depend on an array of $\ell$ rational numbers,
\[\underline{a}=(a_1,\dots, a_\ell)\in \BQ_{\geq 0}^\ell ,\] such that the numerator of $a_i$   keeps track of the singularity type around $D_{i,j}$ and $x_i$ (depending on where the index $i$ corresponds to an end component $P_i$ or a point $x_i$), while the denominator is responsible for the degree of an extension of $R$. 

For example, if $a_i=a_i'/r$,  and $i$ correspond to a point $x_i$, we firstly take the degree $r$ base change, and then by inductively blowing $x_i$ up and contracting all exceptional divisors but the last one, we introduce an end component with the  $A_{a_i'-1}$-singularity.

 If $i$ correspond to an end component $P_i$, after taking a degree $r$ base change, we  inductively contract ruled components $P_{i,j}$ together with the end component $P_i$ (following the index $j$ in the reverse order), and (partially) resolve the $A_{b_{i,j}-1}$-singularity of the non-contracted  component with the largest index $j$; if $a_i'$ is strictly smaller than  $r\sum^{j'}_{j=1}b_{i,j}$ for some $j'$, then all components with the index $j \geq j'$ are contracted.  Note that such operations are independent of the fractional presentation of $a_i$ up to finite base changes. All modifications of $\epsilon_0$-unramified maps are therefore encoded by rational numbers. The precise statement is the following lemma. 
\begin{lemma} \label{neq01d} Assume $\theta^*$ does not have weight $d_0$ end components. 
	For each $\underline{a}=(a_1,\dots, a_\ell)\in \BQ_{\geq 0}^\ell$, such that 
	\[\underline{a}\leq \underline{b},\]
	 there exists a  modification $\theta_{\underline{a}}$ of $\theta_+$ 
	with the following properties:
	\begin{itemize}
		\item $\theta_{\underline{a}}$ is determined by $\underline{a}$,
		\[\theta_{\underline{a}}\cong \theta_{\underline{a}'} \iff \underline{a}=\underline{a}',\]
		up to a finite base change. 
		\item  given a  modification $\theta$ of $\theta_+$, then there exists $\underline{a}\in \BQ_{
		\geq 0}^\ell$, such that  \[\theta \cong \theta_{\underline{a}}\] 
		up to a finite base change. 
		\item   $\theta_{\underline{a}}$ is $\epsilon_+$-unramified, if and only if $\underline{a}=\underline{b}$, and $b_{i}< \infty$ for all $i$,  
		\item $\theta_{\underline{a}}$ is $\epsilon_-$-unramified, if and only if $\underline{a}=\underline{0}$. 
	\end{itemize}
\end{lemma}

\textit{Proof.} We start with the construction of $\theta_{\underline{a}}$. 
Let us choose a fractional presentation of $(a_1,\dots, a_\ell)$ with a common denominator
\[(a_1,\dots, a_\ell)=\left(\frac{a'_1}{r},\dots, \frac{a'_\ell}{r}\right).\]
Let 
\[
\Spec (R') \rightarrow \Spec (R)
\] be a degree $r$ totally ramified cover. We take a base change of $\theta_+$ with respect  to  $\Spec (R') \rightarrow \Spec (R)$. Abusing the notation, we will denote the base change also by $\theta_+$.   We will construct 
\[\theta_{a_i}=(\CC_{a_i},  \CW(m)_{a_i}, \underline{p}_{a_i},f_{a_i})\]
 by blowing up and down $\CW(m)_+$ at the divisors $D_{i,j}$ if the index $i$ corresponds to an component, or at a point $x_i$ if $i$ corresponds to a point.
\\

\textit{Case 1.} If the index $i$ corresponds to an end component $P_i$ and $a'_i=rb_i$, we leave the family $\theta_+$ intact at $P_i$.   Assuming
\[ r\sum^{j'-1}_{j=1}b_{i,j} \leq a'_i < r\sum^{j'}_{j=1}b_{i,j} \quad \text{for some }j' \in \{1,\dots,k_i+1\},\]
such that $\sum^{0}_{j=1}b_{i,j}:=0$, we proceed as follows. 
  We take the resolution of singularities of $ \CW(m)_+$ at $D_{i,j'}$ for all $j\geq j'$; this is a consecutive blow-up at the singular loci. Then, we  consider the divisor  given by the exceptional divisors of the resolution together with proper transforms of ruled and end components  $P_{i,j}$,  
\begin{equation} \label{exc}
P_{i,j'-1} \cup \CE_{j'} \cup P_{i,j'} \cup \CE_{j'+1} \cup P_{i,j'+1} \cup  \hdots \cup  \CE_{k_i} \cup P_{i,k_i} \cup   \CE_{k_i+1} \cup P_{i,k_i+1}, 
\end{equation}
where we set $P_{i,j'-1}=\emptyset$ if $j'=1$, and $P_{i,k_i+1}:=P_i$. 
Around $D_{i,j}$, the family $\CW(m)_+$ has the $A_{b_{i,j}-1}$-singularity, and we took a degree $r$ base change;
the exceptional divisor $\CE_{j}$ of the resolution at $D_{i,j}$ is therefore a chain of $rb_{i,j}-1$ ruled components. We order these ruled components together with $P_{i,j}$'s, such that $P_{i,j'-1}$ is at place 0, the ruled component attached to $P_{i,j'-1}$ is at place 1, etc. We  contract all ruled components and the end component in the divisor (\ref{exc}), except the component $P_{i,j'-1}$ and the component in   $\CE_{j'}$ which is at place 
\[a''_i= a'_i- r\sum^{j'-1}_{j=1}b_{i,j}< rb_{i,j'}.   \]
We denote the resulting family by $\CW(m)_{a_i}$. It has the $A_{a''_i-1}$-singularity at $D_{i,j'}$ (or, more precisely, at the divisor on $P_{i,j'-1}$ corresponding to $D_{i,j'}$ after the birational modifications), to which an end component is attached; this is the unique non-contracted component in $\CE_{j'}$. We may assume that $a''_i$ is divisible by
\[
m_{i,j'}=\mathrm{lcm}\bigl(\eta^{i,j'}_{1}, \hdots, \eta^{i,j'}_{\ell(\eta^{i,j'})}\bigr)
\]
after a finite base change.\footnote{A base change does not affect the ratio $a'_i/r$, since both the numerator and the denominator are multiplied by the same positive integer.}

The map $f_{+}$ then defines a rational map to $\CW(m)_{a_i}$.
We set
\[f_{a_i} \colon \CC_{a_i} \rightarrow \CW(m)_{a_i}\] 
to be the minimal admissible resolution of indeterminacies of this rational map, obtained by modifying $\CC_+$ in exactly the same way as  $\CW(m)_+$. At this stage, it is necessary for $a''_i$ to be a multiple of $m_{i,j'}$, otherwise, the map does not extend.  A local model of such an extension is provided by 
\begin{align*}
f\colon \mathrm{Spf}(\BC[\![ x_1,x_2,t]\!]/(x_1x_2-t^{b/\eta })) &\rightarrow \mathrm{Spf}(\BC[\![ z_1,\hdots,z_{\dim+1},t]\!]/(z_1 z_2-t^{b})),  \\
x_1 \mapsto z^{\eta}_1,& \quad  x_2 \mapsto z^{\eta}_2.
\end{align*}
In this local model, the exceptional divisors of the resolutions of the target and the source consist of chains of $b-1$ and $b/\eta-1$ ruled components, respectively. On the target, we contract all but one exceptional component, denoted by $P$, namely the $a$-th exceptional component when counted from the side of $\{z_1=0\}$. We choose $a$ to be divisible by $\eta$. The resulting singularity is of type $A_{a-1}$.

We proceed similarly for the source, contracting all but the $(a/\eta)$-th exceptional component, denoted by $B$; the resulting singularity is of type $A_{a/\eta-1}$. Around each singularity of the central fiber, the extension of $f$ is locally given by
\[
x'_1 \mapsto (z'_1)^{\eta}, \qquad x'_2 \mapsto (z'_2)^{\eta},
\]
where $x'_i$ and $z'_i$ are local coordinates in Zariski neighborhoods of the singular points on the exceptional components $B$ and $P$, respectively. In particular, $B$ maps onto a line in $P$ by a degree $\eta$ map with exactly two ramification points, each of ramification order $\eta$.

Finally,  the markings $\underline{p}_+$ extend to a markings $\underline{p}_{a_i}$ of $\CC_{a_i}$. We thereby obtain a family of maps  $\theta_{a_i}$, which is also admissible by the local description above. The result  of this modification is a constant end component attached to the divisor $D_{i,j'}$. 
\\

\textit{Case 2.} If the index $i$ corresponds to a base point $x_i$,  we  inductively blow up  $a'_i$ times the family $ \CW(m)_+$, starting with a blow-up at $x_i$ and then continuing with a blow-up at a point  of the exceptional divisor of the previous blow-up away from the divisor at infinity. We then contract all ruled components in the exceptional divisor, leaving the end component intact.  The resulting family is $ \CW(m)_{a_i}$, it has the $A_{a'_i-1}$-singularity at the divisor at infinity $D_i$ lying over $x_i$. As in \textit{Case 1}, we may assume that $a'_i$ is a multiple of the least common multiple of ramification orders (cf.\ the proof of Theorem \ref{epsilonthm}).  The map $f_{C_+}$ defines a rational map to $\CW(m)_{a_i}$. We set 
\[f_{a_i} \colon \CC_{a_i} \rightarrow \CW(m)_{a_i}\]
to be its minimal admissible resolution of indeterminacies, as in \textit{Case 1}. The markings $\underline{p}_+$ extend to the marking $\underline{p}_{a_i}$ of the family $\CC_{a_i}$. 
\\

The family 
\[\theta_{\underline{a}}=(\CC_{\underline{a}} ,  \CW(m)_{\underline{a}} , \underline{p}_{\underline{a}},f_{\underline{a}})\] is constructed by applying all of the above modifications to $D_{i,j}$ and $x_i$. It is not difficult to verify that the central fiber of $\theta_{\underline{a}}$ is indeed $\epsilon_0$-semi-unramified, since we either changed the singularity type of end or ruled components (however, end components might become constant in the process), or moved a weight $d_0$ point on an end component.  A further base change does not change the ratio of the singularity type and the degree of the base change. Up to a finite base change, the resulting family is uniquely determined by $\underline{a}=(a_1,\dots, a_\ell)$ among other families given by the same construction, because of the singularity types around $D_{i,j}$ and $x_i$.  

Given now an arbitrary $\epsilon_0$-semi-unramified modification 
\[\theta=(\CC, \CW(m),\underline{p},f)\] of $\theta_+$. Possibly after a finite base change, there exists  $ \widetilde{ \CW}(m)$ that dominates both $ \CW(m)$ and $ \CW(m)_+$. By taking an admissible resolution of the map  from $\CC$  to  $ \widetilde{ \CW}(m)$, we obtain a family
\[\tilde \theta=(\widetilde{C}, \widetilde{ \CW}(m), \tilde{\underline{p}},\tilde{f}),\]  that dominates both $\theta$ and $\theta_+$.  We take a minimal family with this property. The family $\widetilde{ \CW}(m)$ is a blow-up of $ \CW(m)_+$ and of $ \CW(m)$ at the central fibers. In other words, to obtain $ \CW(m)$ from $ \CW(m)_+$, we firstly blow up $ \CW(m)_+$ and then contract end and ruled components.    By the assumption of minimality, $\epsilon_0$-semi-unramification of $\theta$ and $\theta_+$, and Lemma \ref{contraction}, these must be blow-ups at $D_{i,j}$ and $x_i$ on $ \CW(m)_+$, while contracted components must be components in the resulting exceptional ruled and end components. These are the modifications described in \textit{Steps 1, 2} of the proof. Moreover, the modifications in  \textit{Steps 1, 2} are the only blow-ups followed by contractions that produce $\epsilon_0$-semi-unramified maps.   We conclude that 
\[  \CW(m) \cong  \CW(m)_{\underline{a}}\]
for some $\underline{a}$, after a finite base change.  The uniqueness of maps follows from seperatedness of the moduli space of maps to a fixed target. Hence we obtain that, possibly after a finite base change,
\[\theta \cong \theta_{\underline{a}}\]
for some $\underline{a}\in\BQ^\ell_{\geq 0}$, where $\underline{a}$ is determined by the singularity types of $\CW(m)$ around $D_{j,i}$ and $x_i$.
\qed
\\

\subsubsection{Weight $d_0$ end components at the generic fiber} \label{Secondcase} We will now deal with case when $\theta^*$ has weight $d_0$ end components. Let 
\[\CP^*_1, \dots, \CP^*_\ell \subset \CW(m)^*\] be the set of entangled weight $d_0$ end components in the generic fiber. Using arguments of Theorem \ref{epsilonthm}, we can construct an $\epsilon_0$-semi-unramified extension of $\theta^*$, possibly after a finite base change. As before, we choose a particular extension, 
\[\theta_+:=(\CC_+,  \CW(m)_+, \underline{p}_+,f_+) \in  \FM^{\epsilon_0}_{g,n,m}(X,\beta,\underline{\eta})(R),\]
 which is $\epsilon_+$-unramified over the complement of the central fibers of $\CP_{+,i} \subset \CW(m)_+$ and over central fibers of $\CP_{+,i}$ for which there are no weight $d_0$ base points in the generic fiber, while over central fibers of other $\CP_{+,i}$ we take any extension.  The calibration extends uniquely, hence we obtain a family of $\epsilon_+$-semi-unramified with calibrated components,  
\[\xi_+=(\theta_+, \lambda_+)\in  M\FM^{\epsilon_0}_{g,n,m}(X,\beta,\underline{\eta})(R).\]

  We define an array of $\ell$  numbers associated to the set of entangled weight $d_0$ end components above,
\[\underline{b}=(b_1,\dots,b_\ell)\in (\BQ\cup \{-\infty\})^\ell, \]
by setting $b_i=0$ if $\CP^*$ does not have a weight $d_0$ base point, and $b_i=-\infty$ otherwise.  Modifications of $\theta_+$ can also be classified by an array of $\ell$ rational numbers 
\[\underline{a}=(a_1,\dots, a_\ell)\in \BQ^\ell, \]
 as is shown in the next lemma. 
 \begin{lemma} \label{eqd0} Assume $\theta^*$ has weight $d_0$ end components.	For each $\underline{a}=(a_1,\dots, a_\ell)\in \BQ^\ell$, such that 
 	\[\underline{a}\geq \underline{b},\] there exists a  modification $\theta_{\underline{a}}$ of $\theta_+$ with the following properties: 
 	\begin{itemize}
 		\item $\theta_{\underline{a}}$ is determined by $\underline{a}$,
 		\[\theta_{\underline{a}}\cong \theta_{\underline{a}'} \iff \underline{a}=\underline{a}',\]
 		up to a finite base change. 
 		\item  given a  modification $\theta$ of $\theta_+$, there exist $\underline{a}\in \BQ^\ell$, such that 
 		\[\theta \cong \theta_{\underline{a}}\]
 		up to a finite base change. 
 		\item   $\theta_{\underline{a}}$ is $\epsilon_+$-unramified, if and only if $\underline{a}=\underline{b}$, and $b_i >-\infty$ for all $i$. 
 	\end{itemize}
 	\end{lemma}
 
 \textit{Proof.} We will construct families $\theta_{\underline{a}}$, the properties listed above follow from the same arguments as in Proposition \ref{neq01d}. As before, we choose a fractional presentation of $(a_1,\dots, a_\ell)$ with a common denominator
 \[(a_1,\dots, a_\ell)=\left(\frac{a'_1}{r},\dots, \frac{a'_\ell}{r}\right).\]
  We then take a base change of $\theta_+$ with respect  to totally ramified degree $r$ cover  $\Spec (R') \rightarrow \Spec (R)$. Abusing the notation, we will denote the base change also by $\theta_+$. 
  \\
  
  \textit{Case 1.} Assume $\CP^*_i$ does not have a weight $d_0$ base point. Let  $P_{+,i} \subset \CP_{+,i}$ be the central fiber of $\CP^*_i$, it is an end component of weight $d_0$ possibly with a chain of ruled components of weight 0 attached to it. If $a_i=0$, we leave the family $\theta_+$ intact at $P_{+,i}$. If $a_i\leq b'_i$, where $b'_i$ is the total singularity type of the ruled components of $P_{+,i}$ defined in (\ref{bi}), then like in the first case of the proof of Lemma \ref{neq01d}, we resolve $A_{b'_{ij}-1}$-singularities and then consecutively contract ruled components. The number of ruled components we contract is determined by $a_i$. If $a_i>b'_i$, we firstly contract all ruled components except the last one, which becomes an end component. Then we blow up at its divisor at infinity $D_{+,i}$, and then inductively continue with the blow-up of the divisor at infinity of the exceptional ruled component of the previous blow-up. Afterwards,
  we contract all exceptional ruled components together with the end component except the last exceptional ruled component. The number of times we do this is determined by $a_i$.  We call the resulting family $\CW(m)_{a_i}$. The map $f_{C_+}$ defines a rational map to $\CW(m)_{a_i}$. We set 
  \[f_{a_i} \colon \CC_{a_i} \rightarrow \CW(m)_{a_i}\]
  to be its minimal admissible resolution of indeterminacies.  The marking $\underline{p}_+$ extends to the marking $\underline{p}_{a_i}$ of $\CC_{a_i}$. 
  \\

 \textit{Case 2.} Assume now $\CP^*_i$ has a weight $d_0$ base point. If $a_i=0$, we leave the family $\theta_+$ intact at $P_{+,i}$. If $a_i>0$, we repeat the same consecutive contraction of ruled components followed by blow-ups at the divisor at infinity of $P_{+,i}$ as in the previous case. If $a_i<0$, we consecutively blow up $P_{+,i}$ at the base point, contracting all exceptional ruled components except the last exceptional end component.   
  \qed
  \\
  
For the next theorem, we will use \cite[Lemma 2.9.1]{YZ}, \cite[Lemma 2.9.2]{YZ} and \cite[Lemma 2.10.1]{YZ}, which were proved in the context of a one-dimensional target. The proofs extend to an arbitrary pair $(\FM, \mathfrak{Z})$, such that $\FM$ is a smooth algebraic stack and $\mathfrak{Z}\subset \FM$ is a normal crossing divisor. Our pair ($\FM_{X,m,d}, \mathfrak{Z}_1$) satisfies these requirements.  The only difference is that in our case there might be weight $0$ ruled components. The definition of $b_i$ in (\ref{bi})  accounts for the presence of such components. 

\begin{thm} \label{masterprop}  The moduli space $M\Mbar^{\epsilon_0}_{g,n,m}(X,\beta,\underline{\eta})$ is proper. 
\end{thm}

\textit{Proof.} Let $\xi^*$ be a family of $\epsilon_0$-unramified map with calibrated end components over $\Spec(K)$. 
\\

\textit{Case 1.} Assume $\xi^*$ does not have an end component of weight $d_0$. Consider the extension $\xi_-$. Let 
\[\delta:=\mathrm{ord}(v_{-,1})-\mathrm{ord}(v_{-,2}),\]
where $\mathrm{ord}(...)$ is the order of vanishing of a  section. 
Let $\xi_{\underline{a}}$ be the degree $r$ modification of $\xi_-$, associated to the array of rational numbers,
\[
(a_1,\dots, a_\ell)=\left(\frac{a'_1}{r},\dots, \frac{a'_\ell}{r}\right),
\]
 which was constructed in Lemma \ref{neq01d}. We define 
\[|\underline{a}|:=\sum_i a_i.\]
By \cite[Lemma 2.10.1]{YZ}, we have 
\[\mathrm{ord}(v_{\underline{a},1})-\mathrm{ord}(v_{\underline{a},2})=r (\delta - |\underline{a}|).\] 
Assume $\delta>0$, i.e.,\ $v_{-,1}=0$ at the closed point of $\Spec(R)$, and  $v_{-,2}\neq0$. 
By   \cite[Lemma 2.9.1]{YZ} and the end of  \cite[Lemma 5.2.5]{YZ}, the modification $\xi_{\underline{a}}$ is $\epsilon_0$-unramified, if and only if
\begin{itemize}
	\item $|\underline{a}| \leq \delta,$ 
	\item $0<a_i \leq b_i \text{ for all } i=1, \dots,\ell,$
	\item $\text{if $|\underline{a}|<\delta$, then $a_i=b_i$ for all $i=1,\dots, \ell$, }$
	\item $\text{if $a_i<b_i$, then $a_i$ is maximal among $a_1, \dots, a_\ell$ for all $i=1,\dots, \ell$.}$
	\end{itemize}
There is a unique solution to the system above.  The  modification $\xi_{\underline{a}}$ associated to the this solution is the unique $\epsilon_0$-unramified extension with calibrated end components by Lemma \ref{neq01d}. 

If $\delta\leq 0$, i.e.,\ $v_{-,1}\neq 0$ at the closed point, then $\xi_-=\xi_{\underline{0}}$ is $\epsilon_0$-unramified. It is the unique extension, such that $\theta_-$ is $\epsilon_-$-unramified, by the seperatedness of  $\Mbar^{\epsilon_-}_{g,n,m}(X,\beta,\underline{\eta})$. Any other extensions $\xi_{\underline{a}}$ have a weight-$d_0$ end component and satisfy $\mathrm{ord}(v_{\underline{a},1}) - \mathrm{ord}(v_{\underline{a},2}) < 0$; hence, they are not $\epsilon_0$-unramified.
\\

\textit{Case 2.} Assume now $\xi^*$ has an end component of weight $d_0$. Consider the extension $\xi_+$. Let 
\[\delta:=\mathrm{ord}(v_{+,1})-\mathrm{ord}(v_{+,2}).\]
For a degree $r$ modification $\xi_{\underline{a}}$ of $\xi_+$, constructed in Lemma \ref{neq01d}, let
\[|\underline{a}|:=\sum_i a_i,\]
and let $a_{+,i}$ be the order of vanishing of sections defined by $\CP_{+,i}$ via the construction used in the proof of \cite[Lemma 2.9.2]{YZ}. 
 We have 
\[\mathrm{ord}(v_{\underline{a},1})-\mathrm{ord}(v_{\underline{a},1})=r (\delta + |\underline{a}|).\] 
 By \cite[Lemma 2.9.2]{YZ} and the end of \cite[Lemma 5.4.1]{YZ}, a modification $\xi_{\underline{a}}$ is $\epsilon_0$-unramified, if and only if 
\begin{itemize}
	\item $\delta+|\underline{a}| \geq 0,$
	\item  $a_i \geq b_i \text{ for all } i=1, \dots,\ell,$
	\item $\text{if $\delta+|\underline{a}|>0$, then $a_i=b_i$ for all $i=1,\dots, \ell$, }$
	\item if $a_i>b_i$, then $a_{+,i}-a_i$ is maximal among $a_{+,1}-a_1, \dots, a_{+,\ell}-a_{\ell}$ for all $i=1,\dots, \ell$.
\end{itemize}
There is a unique solution to the system above. The modification $\xi_{\underline{a}}$ associated to this solution is the unique $\epsilon_0$-unramified extension by Lemma \ref{eqd0}.  \qed

\subsection{Obstruction theory of master spaces} \label{obstmasterspace} To endow the master spaces with an obstruction theory, we introduce auxiliary logarithmic structures, as in Section \ref{auxi}. Let 
\[ M\Mbar^{\epsilon_0}_{g,n,m}(X,\beta,\underline{\eta})^{\log} \colon (Sch/\BC)^\circ \rightarrow Grpd \]
 be the moduli space of relative $\epsilon$-unramified log maps in the sense of \cite[Section 5.2]{KiL}. A family over a base scheme $B$ in  $ M\Mbar^{\epsilon_0}_{g,n,m}(X,\beta,\underline{\eta})^{\log}$ consists of the following data:
\begin{itemize}
	\item A locally-free log structure $N_B$ on $B$. 
	\item A minimal log prestable curve 
	\[((\CC,\underline{p}), M_\CC) \rightarrow (B,N_B)\] 
	and an extended log twisted FM degeneration 
	\[(\CW(m), M_{\CW(m)}) \rightarrow (B,N_B), \] 
	satisfying the corank condition \cite[Section 5.2.2 (3)]{KiL}; for both the source and the target, the logarithmic structures do not include markings or exceptional divisors.
	\item A calibration $(e, \CL, v_1,v_2)$ of $\CW(m)$ over $B$. 
	\item A log stable map
	\[
	\begin{tikzcd} [row sep=small, column sep = small]
		((\CC,\underline{p}), M_{\CC} )  \arrow{dr} \arrow{rr}{f} &&  (\CW(m), M_{\CW(m)}) \arrow{dl} \\
		& (B,N_B), 
	\end{tikzcd}
	\]
	whose underlying classical map is  $\epsilon_0$-unramified in the sense of Definition \ref{epsilon0unr}. 
\end{itemize}

Observe that the calibration of $\CW(m)$  and the auxillary logarithmic structures are  independent from each other. More precisely, this means that the moduli space of log twisted FM degenerations with calibrated end components is given by the fibered product
\[
	\begin{tikzcd} [row sep=small, column sep = small]
		  M \widetilde{\FM}^{\log}_{X,m,d} \arrow{d} \arrow{r} &  	  M \widetilde{\FM}_{X,m,d}\arrow{d} \\
	\FM^{\log}_{X,m,d} \arrow{r}  & 	\FM_{X,m,d}, 
\end{tikzcd}
\]
where $\FM^{\log}_{X,m,d}$ is the moduli space of log twisted weighted FM degenerations. As in Section \ref{auxi}, we also have the fibered product
\[
 M\mathfrak{S}:=  M \widetilde{\FM}^{\log}_{X,m,d} \times_{\log}, \FM^{\log}_{g,n+n'}
 \] 
 parametrising log twisted FM degenerations with calibrated end components and minimal log prestable curve over the same logarithmic base. This is a formally smooth stack, since $\mathfrak{S}$ is. Using the same notation for the universal objects as in Section \ref{auxi}, we have the following result.

 	\begin{thm} \label{obstmaster} The space $M\Mbar^{\epsilon_0}_{g,n,m}(X, \beta, \underline{\eta})^{\log}$ is a proper Deligne--Mumford stack, and there exists a morphism 
 		\[
 		\phi \colon (R\pi_*F^*T^{\mathrm{log}}_{\pi, \mathfrak{E}})^\vee  \rightarrow \BL_{M\Mbar^{\epsilon_0}_{g,n,m}(X,\beta,\underline{\eta})^{\log} /M\mathfrak{S}},
 		\]  which defines a perfect obstruction theory on $M\Mbar^{\epsilon_0}_{g,n,m}(X,\beta,\underline{\eta})^{\log}$ relative to $M\mathfrak{S}$. 
 	\end{thm}
 	
 	\textit{Proof} As in the proof of Theorem \ref{obst},  \cite[Theorem 6.3.1]{KiL} and Theorem \ref{masterprop} imply that $M\Mbar^{\epsilon_0}_{g,n,m}(X, \beta, \underline{\eta})^{\log}$ is a proper Deligne--Mumford stack. The existence of the obstruction theory follows from \cite[Section 7]{KiL}. 
 		\qed 
\\

As in Section \ref{Virtual}, we define the virtual fundamental class of $M\Mbar^{\epsilon_0}_{g,n,m}(X, \beta, \underline{\eta})$ by pushing forward the virtual fundamental class from $M\Mbar^{\epsilon_0}_{g,n,m}(X, \beta, \underline{\eta})^{\log}$. The associated invariants are then defined as in Section \ref{Sectioninv}. 
  
\section{Wall-crossing} \label{sectionwallcrossing}

\subsection{Overview of wall-crossing} In this section, we will analyse the torus localisation formula, Theorem \ref{wallcrossing}, associated to the master spaces, constructed in Section \ref{masterspaces}. In Proposition \ref{isomophism}, the  torus fixed components responsible for the wall-crossing terms are described in terms of Vertex spaces.   They arise due to the constant end components of FM degenerations in the master spaces.  

\subsection{Vertex spaces} \label{sectionvertex} Let $\BC^*$ act on  $\p(TX\oplus \BC)$ by scaling the first summand 
\begin{equation} \label{vertexaction}
t \cdot (TX\oplus \BC):=(t\cdot TX)\oplus \BC, \quad t \in \BC^*.
\end{equation}
Consider the divisor at infinity 
\[
\p(TX) \subset  \p(TX\oplus \BC),
\]
 defined by requiring the coordinate  corresponding to the second summand $\BC$ to be zero. We also set
\[\mathbf{z}:=\BC_{\mathrm{std}}, \quad e_{\BC^*}(\mathbf{z})=z,\]
where $\BC_{\mathrm{std}}$ is the weight 1 representation of $\BC^*$ on the vector space $\BC$.

\begin{defn}  We define 
	\[V_{g,n,\eta}(X) \subset \Mbar_{g,n} (\p(TX\oplus \BC)/ X)\] to be the $X$-relative moduli space of stable maps from possibly disconnected curves of genus $g$ with $n$ markings to 
	\[ \p(TX\oplus \BC) \rightarrow X,\]
	which are admissible over 
	\[ 
	\p(TX) \subset \p(TX\oplus \BC),
	\] and each connected component is mapped non-trivially. The ramification profile over $\p(TX)$ is fixed to be $\eta=(\eta_1, \dots, \eta_{\ell(\eta)})$, which also fixes the degree of maps to be $|\eta|=\sum_j\eta_j$. 	We put a standard order on points in $f^{-1}(\p(TX))$ and  endow $V_{g,n,\eta}(X)$ with the obstruction theory of stable maps relative to the divisor $\p(TX)$ in the sense of Li \cite{LiRel}.\footnote{The obstruction-theory complex can be defined as in Section \ref{auxi}; note that expanded degenerations at $\p(TX)$ are not present in the moduli space, hence no auxiliary logarithmic structures associated to the target singularities are needed.}

\end{defn}

%

The space $V_{g,n,\eta}(X)$ inherits a $\BC^*$-action from $\p(TX\oplus \BC)$. Observe that $V_{g,n,\eta}(X)$ is not proper, but its $\BC^*$-fixed locus $V_{g,n,\eta}(X)^{\BC^*}$ is.  
 Since maps belonging to $V_{g,n,\eta}(X)$ are admissible over $\p(TX) \subset \p(TX \oplus \BC)$, we also have  a natural evaluation morphism
\[\widetilde{\ev} \colon V_{g,n,\eta}(X) \rightarrow  \CI_\eta \Sym_X\p(TX) \subset \CI \Sym_X\p(TX),\]
defined as in  Section \ref{sectioninvariants}. 
\begin{defn} \label{Idefinition}  We define the  classes, called the $I$-\textit{function} and its truncation, in the cohomology groups $H^*(\CI_{\eta} \Sym_X \p(TX))[z^{\pm}]$ and  $H^*(\CI_{\eta} \Sym_X \p(TX))[z]$, respectively, 
	\begin{align*}
		I_{g,n,\eta}(z,\prod^{n}_{i=1}\gamma_i \psi_i^{k_i})&:=\prod_j \eta_{j} \cdot\Big(\sum_j \mathrm{c}_{\dim-j}(X)z^{j}\Big) \\
	& \hspace{3cm}\cdot   \widetilde{\ev}_*\left( \frac{ [V_{g,n,\eta}(X)^{\BC^*}]^{\mathrm{vir}}  \cap  \prod^{n}_{i=1} \ev_i^*(\gamma_i)\psi_i^{k_i}}{e_{\BC^*}(N^{\mathrm{vir}}_{V_{g,n,\eta}(X)^{\BC^*}/ V_{g,n,\eta}(X)})} \right), \\
		\mu_{g,n,\eta}(z, \prod^{n}_{i=1}\gamma_i \psi_i^{k_i})&:=[ I_{g,n,\eta}(z, \prod^{n}_{i=1}\gamma_i \psi_i^{k_i})]_{z^{\geq 0}},
	\end{align*}
where $[\hdots]_{z^{\geq0}}$ is the truncation given by taking non-negative powers of $z$. 

\end{defn}

	If $\eta=(1)$, we will  drop $\eta$ from the notation, i.e.,\ we will write 
\[\mu_{g,n}(z, \prod^{n}_{i=1}\gamma_i \psi_i^{k_i}):=\mu_{g,n,(1)}(z, \prod^{n}_{i=1}\gamma_i \psi_i^{k_i}).  \]
The case $\eta=(1)$ will feature  in the wall-crossing for threefolds, Section \ref{sectionthreefolds}.

\begin{rmk} The prefactor 
	\[
	\sum_j \mathrm{c}_{\dim-j}(X)z^{j}
	\]
	is simply the equivariant Euler class  of the tangent bundle $TX\otimes \mathbf{z}$ with the trivial $\BC^*$-action. The tangent bundle $TX$ freely acts on  $V_{g,n,\eta}(X)$ by translations. It will arise in Proposition \ref{isomophism} due to the translation automorphisms of the end components.  
	
The prefactor $\prod_j \eta_j$ in the definition of the $I$-function together with the standard order on the points of contact with the divisor at infinity produce an overall factor 
	\[\Fz(\eta):=|\Aut(\eta) |\cdot \prod_j \eta_j. \]
	Such quantities  appear in degeneration and localisation formulas, e.g.,\ see \cite[Section 3.8]{GV} or \cite[Theorem 2]{MR}.  For the same reasons, they will appear in the analysis of the wall-crossing components in Proposition \ref{isomophism}. The factor $\Aut(\eta)$ arises from the ordering of the nodes of curves mapping to targets with several components, while $\prod_j \eta_j$ arises from the ramification of maps over the singularities of the target. The latter also admits an interpretation in terms of gluing logarithmic structures across the nodes, as explained in Proposition~\ref{isomophism1}.
\end{rmk}

\subsection{Moduli space of end components} \label{sectionend} Let 
\[ G_X\subset \Aut_X(\p (TX \oplus \BC))\]
be the group scheme of $X$-relative automorphisms of $\p (TX \oplus \BC)$ which fix the relative divisor at infinity $\p (T X) \subset \p (TX \oplus \BC)$. We define the moduli stack of end components, 
\[\mathfrak{B}:=BG_X = [X/G_X],\]
which carries a universal family $\mathfrak{P}:=[\p (TX \oplus \BC)/G_X]$  with the universal divisor $\mathfrak{D}:=[\p (T X)/G_X]$ at infinity, 
\[\mathfrak{D}\xhookrightarrow{\iota} \mathfrak{P} \xrightarrow{\pi} \mathfrak{B}. \]
The group $G_X$ fixes $\p(TX)$ pointwise, hence there exists a map 
\[\mathfrak{D} \rightarrow \p(TX).  \]
We denote the $\mathfrak{B}$-relative normal bundle of $\mathfrak{D}$ inside $\mathfrak{P}$ by $\CN_{\mathfrak{D}}$,  and the relative tautological bundle of  $\p(TX)$ by $\CO_{\p(TX)}(1)$.  We define the cotangent line bundle associated to ${\mathfrak{D}}$ and the associated  $\psi$-class as follows, 
\[ \BL(D)^\vee:= \pi_*(\CN_{\mathfrak{D}} \otimes  \iota^*\CO_{\p(TX)}(-1)),  \quad \Psi(D):=\mathrm{c}_1(\BL(D)).\]

Given a weighted FM degeneration $ W(m)$, let
 \[
\{ P_1, \dots, P_k\}
 \] 
 be the set of its end components of weight $d_0$. By \cite[Proposition 2.3]{Frie}, the space 
 \[
 \Theta_i\cong \BC
 \] of the first-order smoothing of  an end component $P_i \subset  W(m)$ admits a natural identification
\[ \Theta_i \cong H^0(D_{i}, \CN_{D_{i}}\otimes \CN_{E_i})\cong (\BL(D_i)|_{W(m)}\otimes \widetilde{\BL}(E_i)|_{W(m)})^\vee,\]
where  $\CN_{D_{i}}$ and  $\CN_{E_i}$ are  the normal bundles of the divisor at infinity and the exceptional divisor associated to $P_i$ (the latter is defined in (\ref{relativepsi})), while $\BL(D_i)|_{W(m)}$ and $\widetilde{\BL}(E_i)|_{W(m)}$ are the fibers of the line bundles $\BL(D_i)$ and $\widetilde{\BL}(E_i)$  over $ [W(m)] \in \FM_{X,m,d}$.  For this reason, the line bundles $\BL(D_i)$ and $\widetilde{\BL}(E_i)$ can be viewed as generalizations of the line bundles associated to the cotangent spaces of markings on curves. 

From the above identification, the fiber of  the calibration bundle $ \BM _{\widetilde{\FM}_{X,m,d}}$ over $[W(m)] \in \FM_{X,m,d}$ is naturally isomorphic to
\begin{equation} \label{theta}
(\BL(D_1)|_{W(m)}\otimes \BL(E_1)|_{W(m)}) \otimes \hdots \otimes  (\BL(D_k)|_{W(m)}\otimes \BL(E_k)|_{W(m)}).
 \end{equation}
 Note that the twist introduced in (\ref{twist}) is necessary for this identification, because to smooth out end components of weight $d_0$, ruled components of weight $0$ attached to them need to be smoothed out as well. 
\subsection{Torus fixed components} \label{secfixed}
There is a natural $\BC^*$-action on  $M\Mbar^{\epsilon_0}_{g,n,m}(X,\beta,\underline{\eta})$ given by  
\begin{equation} \label{actionmaster}
t \cdot (C,  W(m),\underline{p},f, e, \CL, v_1,v_2)= (C,  W(m),\underline{p},f, e, \CL, t \cdot v_1,v_2), \quad t \in \BC^*.
\end{equation}
By \cite[Lemma 6.5.1]{YZ}, an $\epsilon_0$-unramified map with calibrated end components is $\BC^*$-fixed, if and only if it is one of the following: 
\begin{itemize}
	\item[(1)] $v_1=0$ and $(C,  W(m),\underline{p},f)$ is  $\epsilon_+$-unramified,
	\item[(2)] $v_2=0$ and $(C,  W(m),\underline{p},f)$ is $\epsilon_-$-unramified,
	\item[(3)] $v_1\neq 0$ and $v_2\neq0 $, and all weight $d_0$ end components  of $(C,  W(m),\underline{p},f)$ are entangled strictly constant end components. 
	\end{itemize}
The $\BC^*$-fixed locus decomposes according to the types of maps above, 
\begin{equation} \label{fixed}
M\Mbar^{\epsilon_0}_{g,n,m}(X,\beta,\underline{\eta})^{\BC^*}=F_- \sqcup F_+ \sqcup \coprod_{\vec g} F_{\vec{g}},
\end{equation} 
the exact meaning of the components in the above decomposition is explained below. 

\subsection{Fixed components $F_-$ and $F_+$}

The fixed component $F_-$ is the simplest to describe. By \cite[Lemma 6.5.1]{YZ}, it admits the following description, 
\[
F_-= \Mbar^{\epsilon_-}_{g,n,m}(X,\beta,\underline{\eta}), \quad
N_{F_-}=\BM^{\vee}_-,\]
where $\BM^{\vee}_-$ is the dual of the calibration bundle $ \BM_{\widetilde{\FM}_{X,m,d}}$ on $\Mbar^{\epsilon_-}_{g,n,m}(X,\beta,\underline{\eta})$ with a trivial $\BC^*$-action of weight $-1$.

To describe the fixed component $F_+$, consider
\[\widetilde{M}^{\epsilon_+}_{g,n,m}(X,\beta,\underline{\eta}):=  \Mbar^{\epsilon_+}_{g,n,m}(X,\beta,\underline{\eta}) \times_{\FM_{X,m,d}} \widetilde{\FM}_{X,m,d},\] 
then
\[
F_+=\widetilde{M}^{\epsilon_+}_{g,n,m}(X,\beta, \underline{\eta}),\quad
N_{F_+}=\BM_+,
\]
where $\BM_+$ is the calibration bundle  $\BM_{\widetilde{\FM}_{X,m,d}}$ on $\widetilde{M}^{\epsilon_+}_{g,n,m}(X,\beta,\underline{\eta})$ with a trivial $\BC^*$-action of weight $1$. The above fiber product is required because FM degenerations in  $\widetilde{M}^{\epsilon_+}_{g,n,m}(X,\beta,\underline{\eta})$ have end components of weight $d_0$, which are calibrated in the master space by definition. 

\subsection{Virtual fundamental classes of $F_-$ and $F_+$} \label{normalplusminus}

Recall that, in Section \ref{secfixed}, we define the virtual fundamental classes of  master spaces $M\Mbar^{\epsilon_0}_{g,n,m}(X,\beta,\underline{\eta})$ by equipping targets and sources with auxiliary logarithmic structures. Hence, for the torus localisation formulas, a description of the $\BC^*$-fixed locus  $M\Mbar^{\epsilon_0}_{g,n,m}(X,\beta,\underline{\eta})^{\log}$ is needed. Since calibrations of logarithmic FM degenerations are pulled back from the space 	$\FM_{X,m,d}$, a logarithmic map is $\BC^*$-fixed, if and only if the underlying classical map is $\BC^*$-fixed. Hence, the $\BC^*$-fixed locus of the master space admits the same description as in Section~\ref{fixed}, with the only difference being that the sources and the targets are equipped with auxiliary logarithmic structures. This immediately implies that the fixed components

\[
F^{\mathrm{log}}_-, F^{\mathrm{log}}_+ \subseteq M\Mbar^{\epsilon_0}_{g,n,m}(X,\beta,\underline{\eta})
\]
 corresponding to the first two types of  maps in Section \ref{secfixed},  are simply logarithmic versions of their classical counterparts, i.e., they parametrise maps from $F_-$ and $F_+$ with auxiliary logarithmic structures, defined as in Section \ref{auxi}.  Moreover, the normal bundles of $F^{\mathrm{log}}_-$ and $F^{\mathrm{log}}_+$ are pulled back from $F_-$ and $F_+$, and the obstruction theories also match. Hence, we have 
\[
[F^{\mathrm{log}}_-]^{\mathrm{vir}}=[ \Mbar^{\epsilon_-}_{g,n,m}(X,\beta,\underline{\eta})^{\mathrm{log}}]^{\mathrm{vir}}
\]
and
\[
[F^{\mathrm{log}}_+]^{\mathrm{vir}}=[ \widetilde{M}^{\epsilon_+}_{g,n,m}(X,\beta,\underline{\eta})^{\mathrm{log}}]^{\mathrm{vir}},
\]
such that 
\begin{equation} \label{tildenottilde}
\pi_*[\widetilde{M}^{\epsilon_+}_{g,n,m}(X,\beta,\underline{\eta})]^{\mathrm{vir}}=[\Mbar^{\epsilon_+}_{g,n,m}(X,\beta,\underline{\eta})]^{\mathrm{vir}},
\end{equation}
where \[\pi \colon \widetilde{M}^{\epsilon_+}_{g,n,m}(X,\beta,\underline{\eta}) \rightarrow  \Mbar^{\epsilon_+}_{g,n,m}(X,\beta,\underline{\eta}) \] 
is the natural projection.

\subsection{Fixed components $F_{\vec{g}}$} \label{componentswall} The components $F_{\vec{g}}$ are responsible for the difference between the $\epsilon$-unramified invariants associated with $\epsilon_-$ and $\epsilon_+$. The closed points of $F_{\vec{g}}$ correspond to $\epsilon_0$-semi-unramified maps with entangled, strictly constant end components; that is, to the fixed maps of the third type described in Section~\ref{secfixed}.
 The vector $\vec{g}$ specifies the discrete data of the maps over these end components in terms of the genus, the number of markings, and the ramifications over the divisor at infinity. 

\begin{defn} \label{PartitionWeight}  We define  a \textit{stable partition} of $(g, \{1,\dots, n\})$  of weight $d_0$ and length $k$ as 
\[\vec{g}:=((g_0,N_0),(g_1,N_1, \eta^{m+1}),\dots,(g_k,N_k,\eta^{m+k})),\] 
where $g_i \in \BZ$,  $N_i \subseteq \{1,\dots,  n\}$, and $\eta^{m+i}=(\eta_1^{m+i},\dots, \eta_{\ell(\eta^{m+i})}^{m+i})$ are some  partitions. The sets $N_i$ are allowed to be empty, and the data above is subject to the following conditions:
\begin{itemize}
	\item $N_i \cap N_j=\emptyset \text{ if }i\neq j \text{, and }  \bigcup_i N_i =\{1,\dots,  n\},$
	\item $2g_i-2+ |\eta^{m+i}|+\ell(\eta^{m+i})+|N_i|/2=d_0, \text{ for all }i\neq 0,$
 \item $\sum^{k}_{i=1}(g_i+ \ell(\eta^{m+i}))+ g_0-k=g.$
\end{itemize}
 We say that an end component $P_i$ is of \textit{refined weight} $(g_i,N_i, \eta^{m+i})$  if maps to $P_i$ are of genus $g_i$ with  markings  given by the set $N_i$ and with the ramification profile $\eta^{m+i}$ over the divisor at infinity. 
\end{defn}
Using the above terminology, by  \cite[Lemma 6.5.1]{YZ}, we have 

\begin{equation*}	
F_{\vec{g}}:=	\left\{\xi \  \Bigl\rvert \ \arraycolsep=0.1pt\def\arraystretch{1} \begin{array}{c} \xi \text{ has exactly }k \text{ entangled end components, which are } \\[.001cm] \text{ strictly constant of refined weight } (g_i,N_i, \eta^{m+i}) \end{array} \right\}. 
\end{equation*}
We provide a more explicit description of $F_{\vec{g}}$ below. 

\subsection{Explicit description of $F_{\vec{g}}$}
Let 
\[ 
\mathfrak{D}_i \hookrightarrow	\mathfrak{P}_i \rightarrow  F_{\vec{g}}, \quad i=1, \dots, k,
\]
be the universal $i$-th entangled end component with the universal  divisor at infinity $\mathfrak{D}_i $, and
 let 
\begin{equation} \label{gluingmap}
\widetilde{\mathrm{gl}}_k \colon \widetilde{\FM}_{X, m+k,d-kd_0}\times_{X^k} \mathfrak{B}^k \rightarrow \widetilde{\FM}_{X, m, d}
 \end{equation}
be the gluing map, obtained by gluing an end component to the exceptional divisor of a blown-up FM degeneration. To make this gluing map compatible with the notion of a family from Definition \ref{family}, we can equivalently define $\mathfrak{B}$ via families of end components in the boundary divisors of $\widetilde{\FM}_{X, m, d}$.  See \cite[Section 3.3]{CGK} for the construction of gluing maps for the Fulton--MacPherson compactifications. 

Let  $\widetilde{\mathrm{gl}}_k^*F_{ \vec{g}, \mathrm{ord}}$ be the pullback of $F_{ \vec{g}}$ with respect to $\widetilde{\mathrm{gl}}_k$ together with the standard order on the nodes above the singularities of strictly constant end components, cf.\ Definition \ref{def}. Consider the line bundles $\BL(D_i)$ and $\BL(E_{m+i})$ associated to the divisor at infinity $\mathfrak{D}_i$ on $\mathfrak{P}_i$ and to the universal exceptional divisor $\mathfrak{E}_{m+i}$ on $\widetilde{\mathrm{gl}}_k^*F_{\vec{g},\mathrm{ord}}$  to which $\mathfrak{D}_i$ is glued, respectively. By the stability, sections $v_1$ and $v_2$ are both non-zero everywhere on $\widetilde{\mathrm{gl}}_k^*F_{\vec{g}, \mathrm{ord}}$, hence the calibration bundle $\BM_{d}$ of $\widetilde{\FM}_{X, m, d}$ trivializes over $\widetilde{\mathrm{gl}}_k^*F_{\vec{g}, \mathrm{ord}}$, 
\begin{equation} \label{triv}
\BM_{d} \overset{v_1/v_2}{\cong} \CO_{\widetilde{\mathrm{gl}}_k^*F_{\vec{g}, \mathrm{ord}}}.
\end{equation}
By (\ref{theta}), there is a canonical isomorphism on $\widetilde{\mathrm{gl}}_k^*F_{\vec{g}, \mathrm{ord}}$, 
\begin{equation} \label{bundles}
	\BM_{d-kd_0} \otimes \BL(E_{m+1}) \otimes \BL(D_1) \otimes \hdots \otimes  \BL(E_{m+k}) \otimes \BL(D_k) \cong \BM_{d} ,
	\end{equation}
where $\BM_{d-kd_0}$ is the calibration bundle of $ \widetilde{\FM}_{X, m+k,d-kd_0}$. 
 Since $\mathfrak{P}_i$ are entangled, \cite[Lemma 2.5.5]{YZ} gives  a canonical identification 
\begin{equation} \label {iso} \BL(E_{m+1})\otimes \BL(D_1) \cong   \dots \cong \BL(E_{m+k}) \otimes \BL(D_k).
	\end{equation}
From (\ref{triv}), (\ref{bundles}), and (\ref{iso}), we obtain that 
\begin{equation} \label{kpower}
(\BL(E_{m+1}) \otimes \BL(D_1))^k \cong \BM_{d-kd_0}^\vee.
 \end{equation}
Let 
\begin{equation} \label{spaceY}
	Y \rightarrow \widetilde{M}^{\epsilon_+}_{g_0,N_0,m+k}(X,\beta, \underline{\eta}, \underline{\eta}'):=\widetilde{M}^{\epsilon_+}_{g_0,N_0,m+k}(X,\beta, \underline{\eta}, \eta^{m+1}, \dots, \eta^{m+k})
\end{equation}	
	 be the stack of $k$-roots of $\BM_{d-kd_0}$ with the universal root $L$,
\[
L^k=\BM_{d-kd_0}.
\]
 Then, by (\ref{bundles}), we obtain an induced morphism, 
\begin{equation} \label{first}
 \widetilde{\mathrm{gl}}_k^*F_{\vec{g}, \mathrm{ord}} \rightarrow Y. 
 \end{equation}
We will now construct a space over which the family of end components $\mathfrak{P}_i$ can be canonically trivialized.  Namely, let 
\[
V^*_i \rightarrow Y
\] be the total space of 
\[
L\otimes \BL(E_{m+i})
\] without the zero section. We define
\[Y':= V^*_1 \times_Y \dots \times_Y V^*_k.\]
Consider now the base change 
\begin{equation} \label{basechangeY}
\widetilde{\mathrm{gl}}_k^*F_{\vec{g}, \mathrm{ord}} \times_Y Y'.
\end{equation}
By the definition of $L$ and $V_i^*$, and the identification (\ref{kpower}), a $B$-valued point of $\widetilde{\mathrm{gl}}_k^*F_{\vec{g}, \mathrm{ord}} \times_Y Y'$ is given  a $B$-valued point of $\widetilde{\mathrm{gl}}_k^*F_{\vec{g}, \mathrm{ord}}$ together with nowhere-vanishing sections 
\[
s_i \in H^0(B,\BL(D_i)^\vee_B).
\] 
Let $\mathfrak{P}'_i$ be the pullback of the universal entangled end components $\mathfrak{P}_i$  from   $\widetilde{\mathrm{gl}}_k^*F_{\vec{g}, \mathrm{ord}}$ to  $\widetilde{\mathrm{gl}}_k^*F_{\vec{g}, \mathrm{ord}} \times_Y Y'$. We choose a non-zero vector 
\[
v \in H^0(\BP(TX), \CN_{\BP(TX)} \otimes \CO_{\BP(TX)}(-1) ),
\]
then by the moduli interpretation of $\widetilde{\mathrm{gl}}_k^*F_{\vec{g}, \mathrm{ord}} \times_Y Y'$, we obtain the following morphisms over $X$,
\begin{equation} \label{induced}
 \mathfrak{P}_i' \rightarrow \BP(TX \oplus \BC), 
 \end{equation}
by sending the unique base point to $0$, the divisor at infinity $\mathfrak{D}_i$  to $\p(TX)$ and $s_i$ to $v$. Note that, in a smooth chart provided by  $X[n]$,  divisors at infinity of end components  can be canonically identified with $\BP(T_xX)$. This identification descends from a chart, as the permitted automorphisms of end components fix the divisor at infinity pointwise.    By composing $\epsilon_0$-semi-unramified maps with (\ref{induced}), we obtain 
\begin{equation} \label{second1} \widetilde{\mathrm{gl}}_k^*F_{\vec{g}, \mathrm{ord}}\times_Y Y' \rightarrow V_{ g_i, N_i,\eta^{m+i}}. 
\end{equation}
Since the end components are strictly constant, the morphism (\ref{second1}) lands into the $\BC^*$-fixed locus  
\[
V_{ g_i, N_i,\eta^{m+i}}^{\BC^*} \subset V_{ g_i, N_i,\eta^{m+i}}. 
\] 
Hence, this morphism is invariant under the $(\mathbb{C}^*)^k$-action given by scaling the sections $s_i$ (see (\ref{weighttwist}) for the relation between the scaling actions on both sides of (\ref{second1})). The quotient of this action is $\widetilde{\mathrm{gl}}_k^*F_{\vec{g}, \mathrm{ord}}$,
\[
\widetilde{\mathrm{gl}}_k^*F_{\vec{g}, \mathrm{ord}}\times_Y Y'/(\BC^*)^k=\widetilde{\mathrm{gl}}_k^*F_{\vec{g}, \mathrm{ord}}. 
\]
 The morphism (\ref{second1}) therefore descends to a morphism 
\begin{equation}\label{second2}
\widetilde{\mathrm{gl}}_k^*F_{\vec{g}, \mathrm{ord}} \rightarrow V^{\BC^*}_{ g_i, N_i,\eta^{m+i}}.
\end{equation}
Combining (\ref{first}) with (\ref{second2}), we obtain a morphism 
\begin{equation}\label{themap}
\widetilde{\mathrm{gl}}_k^*F_{\vec{g}, \mathrm{ord}}  \rightarrow Y \times_{\CI\Sym_X \p(TX)^k} \prod_{i=1}^{k} V^{\BC^*}_{ g_i, N_i,\eta^{m+i}},
\end{equation}
which is in fact an isomorphism. 
\begin{prop} \label{isomophism0} The morphism (\ref{themap}) is an isomorphism,
	\[ 
	\widetilde{\mathrm{gl}}_k^*F_{\vec{g}, \mathrm{ord}}  \cong Y \times_{\CI\Sym_X \p(TX)^k} \prod_{i=1}^{k} V^{\BC^*}_{ g_i, N_i,\eta^{m+i}}. 
	\]
\end{prop}	
\textit{Proof.}
The claim follows from exhibiting an explicit inverse morphism, the construction of which amounts to equipping the end components $\mathfrak{P}_i$ over
\[
Y \times_{\CI\Sym_X \p(TX)^k} \prod_{i=1}^{k} V^{\BC^*}_{ g_i, N_i,\eta^{m+i}},
\] corresponding to the factors $V^{\BC^*}_{ g_i, N_i,\eta^{m+i}}$, with calibrations.  The argument is identical to the one presented in \cite[Lemma, 6.5.5]{YZ}. 
\qed 
\\

\subsection{Virtual fundamental class and normal complex of $F_{\vec{g}}$}

Consider now a fixed locus of the third type with auxiliary logarithmic structures,
\[
F_{\vec{g}}^{\mathrm{log}} \subseteq (M\Mbar^{\epsilon_0}_{g,n,m}(X,\beta,\underline{\eta})^{\log})^{\BC^*}.
\]
Let 
\[
 b_i:=\mathrm{lcm}( \eta^{m+i}_1,\hdots, \eta^{m+i}_{\ell(\eta^{m+i})}) \in \BN,
 \]
 where $\eta^{m+i}$'s are the ramification profiles  over the end components, cf.\ Proposition \ref{isomophism0}. Let $Y^{\mathrm{log}}$ be defined by putting auxiliary logarithmic structures on maps belonging to the space $Y$ from (\ref{spaceY}), as in Section \ref{auxi}.
The following proposition provides a description of $\widetilde{\mathrm{gl}}_k^*F_{\vec{g},\mathrm{ord}}^{\mathrm{log}}$ similar to the one of $\widetilde{\mathrm{gl}}_k^*F_{\vec{g},\mathrm{ord}}$ from Proposition \ref{isomophism0}. Note that $\widetilde{\mathrm{gl}}_k$ is defined as in (\ref{gluingmap}), i.e., it corresponds to ordering the end components.  

\begin{prop} \label{isomophism1} There exists a morphism 
	\[
q\colon 	\widetilde{\mathrm{gl}}_k^*F^{\mathrm{log}}_{\vec{g}, \mathrm{ord}}  \rightarrow Y^{\mathrm{log}} \times_{\CI\Sym_X \p(TX)^k} \prod_{i=1}^{k} V^{\BC^*}_{ g_i, N_i,\eta^{m+i}},
	\] 
	which is \'etale  of degree 
	\[ 
	\frac{\prod_{i,j} \eta_j^{m+i}}{ \prod_i b_i}. 
	\]
\end{prop}	

\textit{Proof.} 
This is the same morphism as the one constructed in (\ref{themap}). Since calibrations of logarithmic FM degenerations are pulled back from the space 	$\FM_{X,m,d}$, a logarithmic map is $\BC^*$-fixed, if and only if the underlying classical map is $\BC^*$-fixed.  However,  in this case,  the morphism $q$ is not an isomorphism because the targets and sources in $\widetilde{\mathrm{gl}}_k^*F_{\vec{g}, \mathrm{ord}}^{\mathrm{log}} $ have auxiliary logarithmic structures associated to the singularities at which constant end components are attached, while in 
\[
Y^{\mathrm{log}} \times_{\CI\Sym_X \p(TX)^k} \prod_{i=1}^{k} V^{\BC^*}_{ g_i, N_i,\eta^{m+i}},
\] 
they are not present. As explained in \cite[Remark 6.3.1]{KiL}, these logarithmic structures are determined by the ramification profiles together with a choice of  $\eta^{m+i}_j$-th roots of unity. Roots of unity gives rise to isomorphic logarithmic structures if they are related by the action of  $b_i$-th roots of unity. This leads to the  degree of $q$ being equal to 
\[
\frac{\prod_{i,j} \eta_j^{m+i}}{ \prod_i b_i}. 
\]
We refer to \cite[Section 5.3]{MR} for more details on the splitting of logarithmic maps in the context of localisation formulas. See also \cite[Lemma 7.9.1]{ChenQ} and \cite[Lemma 9.2 (4)]{KLR} for discussions of the same phenomenon in the context of degeneration formulas in the logarithmic setting, and \cite{AF} for the orbifold setting.  \qed
\\

 For brevity, let
  \[	
 \CI_{g,n,\eta}(z):= e_{\BC^*}(N^{\mathrm{vir}}_{V_{g,n,\eta}(X)^{\BC^*}/ V_{g,n,\eta}(X)})^{-1} \in H^*(V_{g,n,\eta}(X)^{\BC^*})[z^\pm].
 \]
 In the next proposition, we finally determine the fundamental class and the virtual normal complex of $F_{\vec{g}}^{\mathrm{log}}$.

\begin{prop} \label{isomophism}  With respect to the morphism from Proposition \ref{isomophism1}, we have 
	\begin{align*} q_*([\widetilde{\mathrm{gl}}_k^*F^{\mathrm{log}}_{\vec{g}, \mathrm{ord}}]^{\mathrm{vir}}) =& \frac{\prod_{i,j} \eta_j^{m+i}}{ \prod_i b_i} \cdot  [Y^{\mathrm{log}}]^{\mathrm{vir}} \times_{\CI \Sym_X \p(TX)^k} \prod_{i=1}^{k} [V^{\BC^*}_{ g_i, N_i,\eta^{m+i}}]^{\mathrm{vir}},  \\	e_{\BC^*}(\widetilde{\mathrm{gl}}_k^* N_{F^{\mathrm{log}}_{\vec{g}, \mathrm{ord}}}^{\mathrm{vir}})^{-1}=&  \prod_i b_i \cdot q^* \left(   \frac{ \prod^k_{i=1}(\sum_j \mathrm{c}_{\dim-j}(X)(z/k+\Psi(D_i))^j) }{-z/k-\Psi(D_1)-\Psi_{m+1}-\sum^{\infty}_{i=k}\mathfrak{Y}_i} \right)  \\
	&	 \hspace{3.5cm} \cdot q^*\left(\prod^k_{i=1} \CI_{g_i,N_i,\eta^{m+i}}(z/k+\Psi(D_i)) \right).
	\end{align*}
\end{prop}
\textit{Proof.}  The claim follows from the analysis presented in \cite[Lemma 6.5.6]{YZ} with the main difference that the group of automorphisms of end components which fix the divisor at infinity and acts trivially on $\BL(D_j)$ is given by translations by $TX$. For this reason, Chern classes of $X$ appear in the numerator of the first factor in the expression of the Euler class of the normal bundle. Let us now summarise the analysis of \cite[Lemma 6.5.6]{YZ} applied to our setting. 

First, the space, introduced in (\ref{basechangeY}),
 \[ 
\widetilde{\mathrm{gl}}_k^*F_{\vec{g}, \mathrm{ord}}\times_Y Y' 
 \] admits a $\BC^* \times (\BC^*)^k$-action, such that the first factor $\BC^*$ comes form the scaling action on the master space (\ref{actionmaster}), while $(\BC^*)^k$ scales the sections $s_i$. On the other hand, the space
 \[
 Y \times_{\CI\Sym_X \p(TX)^k} \prod_{i=1}^{k} V_{ g_i, N_i,\eta^{m+i}}
 \] 
 admits a $(\BC^*)^k$-action given by scaling the factors $V_{ g_i, N_i,\eta^{m+i}}$ by (\ref{vertexaction}). With respect to the natural morphism, 
 \[
 \widetilde{\mathrm{gl}}_k^*F_{\vec{g}, \mathrm{ord}}\times_Y Y' \rightarrow  Y \times_{\CI\Sym_X \p(TX)^k} \prod_{i=1}^{k} V_{ g_i, N_i,\eta^{m+i}}
 \]
 induced by (\ref{first}) and (\ref{second1}), the two actions are related by 
\begin{align}\label{torusaction}
	\begin{split}
		\BC^*\times (\BC^*)^k &\rightarrow (\BC^*)^k \\
		(t, t_1, \hdots, t_k) &\mapsto (t^{1/k}t_1^{-1},\hdots, t^{1/k}t_k^{-1}).
	\end{split}
\end{align}
By the definition of $Y'$, after taking the quotient by $(\BC^*)^k$, this implies that the weight $z$ of $\BC^*$ acting on $\widetilde{\mathrm{gl}}_k^*F_{\vec{g}, \mathrm{ord}}$  is therefore identified with 
\begin{equation} \label{weighttwist}
	z/k+\Psi(D_i)
\end{equation}
on each factor $V_{ g_i, N_i,\eta^{m+i}}^{\BC^*}$ with respect to the isomorphism (\ref{themap}). The same conclusion holds for the logarithmic versions of these spaces.

 Next, we separate the absolute obstruction theory $\BL_{M\Mbar^{\epsilon_0}}^{\mathrm{vir}}$ of $M\Mbar^{\epsilon_0}_{g,n,m}(X,\beta,\underline{\eta})^{\mathrm{log}}$ into the  cotangent complex $\BL_{M\widetilde{\FM}^{\mathrm{log}}_{X,m,d}}$ of $M\widetilde{\FM}_{X,m,d}^{\mathrm{log}}$ and the $M\widetilde{\FM}^{\mathrm{log}}_{X,m,d}$-relative  obstruction theory $\BL^{\mathrm{vir}}_{\mathrm{rel}}$ of maps  to extended log twisted FM degenerations, 
\begin{equation} \label{relativeabs}
\BL_{M\widetilde{\FM}^{\log}_{X,m,d}} \rightarrow \BL_{M\Mbar^{\epsilon_0}}^{\mathrm{vir}} \rightarrow \BL^{\mathrm{vir}}_{\mathrm{rel}} \rightarrow. 
\end{equation}
Similarly to the one-dimensional case, arguments from \cite[Lemma 2.5.5]{YZ} and   \cite[Lemma 6.5.6]{YZ} imply that  the fixed and moving parts of $\BL_{M\widetilde{\FM}^{\mathrm{log}}_{X,m,d}}|_{\widetilde{\mathrm{gl}}_k^*F^{\mathrm{log}}_{\vec{g}, \mathrm{ord}}}$ admit the following expressions, 
\begin{align*} (\BL_{M\widetilde{\FM}^{\mathrm{log}}_{X,m,d}}|_{{\widetilde{\mathrm{gl}}_k^*F^{\mathrm{log}}_{\vec{g}, \mathrm{ord}}}})^f &\cong \BL_{\widetilde{\FM}^{\mathrm{log}}_{X,d-kd_0,m+k}} |_{\widetilde{\mathrm{gl}}_k^*F^{\mathrm{log}}_{\vec{g}, \mathrm{ord}}} , \\
	e_{\BC^*}\left( (\BL^\vee_{M\widetilde{\FM}^{\mathrm{log}}_{X,m,d}}|_{\widetilde{\mathrm{gl}}_k^*F^{\mathrm{log}}_{\vec{g}, \mathrm{ord}}})^m \right)^{-1}&=\prod_i b_i \cdot \frac{  \prod^k_{i=1}(\sum_j \mathrm{c}_{\dim-j}(X)(z/k+\Psi(D_i))^j)}{  -z/k-\Psi(D_1)-\Psi_{m+1}-\sum^{\infty}_{i=k}\mathfrak{Y}_i}.
\end{align*}
The denominator  of the expression for the Euler class  is the normal bundle of the locus of FM degenerations with exactly $k$ entangled end components by \cite[Lemma 2.5.5]{YZ}. 
By the expression of the formal versal deformation space from \cite[Section 6.2.2]{KiL}, the space of smoothings  $\Theta_i$ of end components of a logarithmically twisted FM degeneration is related to the space of smoothings $\Theta'_i$ of end components of the underling classical FM degeneration by the identity
\[
(\Theta_i)^{b_i} \cong \Theta'_i,
\]
such that $b_i$ is determined the ramification profiles over end components, as explained in \cite[Remark 6.3.1]{KiL}; see also \cite[Corollary 5, Corollary 7]{MR}. This results in the prefactor  $(b_i)^{-1}$ in  \cite[Lemma 2.5.5]{YZ} for each end component, and therefore in the prefactor $\prod b_i$ above. On the other hand, the numerator 
\[
\sum_j \mathrm{c}_{\dim-j}(X)z^j
\]
is the translation automorphisms of end components, which are identified with   $TX$; these are the automorphisms of end components that fix the unique base point together with the divisor at infinity, acting trivially on the tangent spaces of the divisor.  The substitution of $z$ by $z/k+\Psi(D_i)$ is explained in the preceding paragraph. 

Let us now consider the fixed and moving parts of  $\BL^{\mathrm{vir}}_{\mathrm{rel}}$ from (\ref{relativeabs}). Using the space 
\[
\widetilde{\mathrm{gl}}_k^*F^{\mathrm{log}}_{\vec{g}, \mathrm{ord}}\times_{Y^{\mathrm{log}}} Y'
\] and the $(\BC^*)^k$-scaling action on it, the fixed and moving parts of $ \BL^{\mathrm{vir}}_{\mathrm{rel}}$  admit the following expressions in terms of the obstruction theory of Vertex spaces, 
\begin{align*}
 (\BL^{\mathrm{vir}}_{\mathrm{rel}}|_{\widetilde{\mathrm{gl}}_k^*F^{\mathrm{log}}_{\vec{g}, \mathrm{ord}}})^f&\cong \BL^{\mathrm{vir}}_{Y^{\mathrm{log}}/\widetilde{\FM}^{\mathrm{log}}_{X,d-kd_0,m+k} } \boxplus \bigoplus^{k}_{i=1} (\BL^{\mathrm{vir}}_{ V^{\mathrm{log}}_{ g_i, N_i,\eta^{m+i}}}|_{\big(V^{{\mathrm{log}}}_{ g_i, N_i,\eta^{m+i}}\big)^{\BC^*}})^{f} , \\
	e_{\BC^*}\left( (\BL^{\mathrm{vir}}_{\mathrm{rel}}|_{\widetilde{\mathrm{gl}}_k^*F^{\mathrm{log}}_{\vec{g}, \mathrm{ord}}})^m \right)^{-1}&= \prod^k_{i=1} \CI_{g_i,N_i,\eta^{m+i}}(z/k+\Psi(D_i)).
\end{align*}
The difference with  \cite[Lemma 6.5.6]{YZ} is that Vertex spaces (or graph spaces in the quasimap terminology) are now maps to $\p(TX\oplus \BC)$ instead of quasimaps to $\CX \times \p^1$; see also \cite{MR} for a discussion of localisation formulas in the logarithmic setting. With the expression of the degree of $q$ from Proposition \ref{isomophism}, this finishes the proof. 
\qed
\\ 
\begin{rmk} \label{discon} The reason we choose our moduli spaces to parametrize maps with disconnected sources is because the decomposition of Proposition \ref{isomophism} is not compatible with connected curves. More specifically, curves in the first factor may be disconnected, however, they might become connected after attaching them to curves in the second factor. Hence we inevitably would need to consider disconnected curves. 
	\end{rmk}
Using the analysis of \cite[Section 7]{YZ}, we obtain the wall-crossing formula. 
	\begin{thm} \label{wallcrossing} We have
		\begin{multline*} 
			\langle \gamma_1 \psi^{k_1}_1 \cdot \!\cdot \!\cdot \gamma_n \psi^{k_n}_n \mid \gamma'_1 \Psi^{k'_1}_1\cdot \!\cdot \!\cdot \gamma'_m \Psi^{k'_m}_m   \rangle^{\epsilon_{-} }_{g,\beta}	-	\langle \gamma_1 \psi^{k_1}_1\cdot \!\cdot \!\cdot \gamma_n \psi^{k_n}_n \mid \gamma'_1 \Psi^{k'_1}_1\cdot \!\cdot \!\cdot \gamma'_m \Psi^{k'_m}_m   \rangle^{\epsilon_{+} }_{g,\beta}\\
			=\sum_{\vec{g}} \Big\langle \prod_{j\in N_0} \gamma_j \psi^{k_j}_j \ \Big| \ \gamma'_1 \Psi^{k'_1}_1\cdot \!\cdot \!\cdot \gamma'_m \Psi^{k'_m}_m \cdot  \prod^{k}_{i=1} \mu_{g_i,N_i, \eta^{m+i}}\Big( -\Psi_{m+i}, \prod_{j\in N_i} \gamma_j \psi^{k_j}_j \Big)  \Big\rangle^{\epsilon_+}_{g_0,\beta} \big/ k!,
		\end{multline*}
		where we sum over stable partitions 
		\[
		\vec{g}=((g_0,N_0),(g_1,N_1,\eta^{m+1}),\dots,(g_k,N_k,\eta^{m+k}))
		\] of $(g,\{1,\dots,n\})$ of  weight $d_0$ and length $k$ for all $k\geq 1$.
	\end{thm}
\textit{Proof.} Let us first introduce a compact notation for the total insertion with respect to markings on the source curve and exceptional divisors on the target, 
\[A:= \prod^{n}_{i=1}\ev^*_i(\gamma_i) \psi^{k_i}_i, \quad B:=\prod^{m}_{j=1}\widetilde{\ev}^*_j(\gamma'_j)\Psi^{k'_j}_j, \]
and, more generally, for a possibly  empty subset $N_i \subseteq \{1,\dots , n\}$, we denote 
\[ A_{N_i}:= \prod_{j\in N_i}\ev^*_j(\gamma_j) \psi^{k_j}_j. \]
  By Theorem \ref{masterprop} and Theorem \ref{obstmaster}, the master space $M\Mbar^{\epsilon_0}_{g,n,m}(X,\beta,\underline{\eta})^{\mathrm{log}}$ is proper. Hence by the virtual localization formula \cite{GP}, the sum of the $z=0$ residues of the localised classes associated to its $\BC^*$-fixed components (\ref{fixed}) is zero, 

\begin{multline} \mathrm{Res}_z\left( \frac{[F^{\mathrm{log}}_-]^{\mathrm{vir}}}{e_{\BC^*}(N^{\mathrm{vir}}_{F^{\mathrm{log}}_-})}+\frac{[F^{\mathrm{log}}_+]^{\mathrm{vir}}}{e_{\BC^*}(N^{\mathrm{vir}}_{F^{\mathrm{log}}_+})}+\sum_{\vec{g}}\frac{[F^{\mathrm{log}}_{\vec{g}}]^{\mathrm{vir}}}{e_{\BC^*}(N^{\mathrm{vir}}_{F^{\mathrm{log}}_{\vec{g}}})} \right) \\
	=\mathrm{Res}_z([M\Mbar^{\epsilon_0}_{g,n,m}(X,\beta,\underline{\eta})^{\mathrm{log}}]^{\mathrm{vir}})=0,
	\end{multline}
where $\mathrm{Res}_z$ is the $z=0$ residue, i.e.,\ the coefficient of the term $1/z$ of a Laurent polynomial.  We interpret the identity above as a relation of cycles on a space to which all fixed loci admit a projection, like a point. 
By Section \ref{normalplusminus}, the virtual tangent normal complexes of $F^{\mathrm{log}}_-$ and $F^{\mathrm{log}}_+$ are line bundles of weight $-1$ and $1$, respectively. Hence,
\[ \mathrm{Res}_z\left( \frac{[F^{\mathrm{log}}_{-}]^{\mathrm{vir}}}{e_{\BC^*}(N^{\mathrm{vir}}_{F^{\mathrm{log}}_{-}})}\right)=-[F^{\mathrm{log}}_{-}]^{\mathrm{vir}}, \quad \mathrm{Res}_z\left( \frac{[F^{\mathrm{log}}_{+}]^{\mathrm{vir}}}{e_{\BC^*}(N^{\mathrm{vir}}_{F^{\mathrm{log}}_{+}})}\right)=[F^{\mathrm{log}}_{+}]^{\mathrm{vir}}.\]
Using the description of the fixed loci from Section \ref{normalplusminus},   inserting the class $A\cdot B$ and pushing forward  to the point the residue relation, we therefore obtain 
\begin{equation} \label{theformula}
\int_{[ \Mbar^{\epsilon_-}_{g,n,m}(X,\beta,\underline{\eta})^{\mathrm{log}}]^{\mathrm{vir}}}A\cdot B- \int_{[ \Mbar^{\epsilon_+}_{g,n,m}(X,\beta,\underline{\eta})^{\mathrm{log}}]^{\mathrm{vir}}} A\cdot B=\sum_{\vec{g}}\mathrm{Res}_z \bigg(\int_{\frac{[F^{\mathrm{log}}_{\vec{g}}]^{\mathrm{vir}}}{e_{\BC^*}(N^{\mathrm{vir}}_{F^{\mathrm{log}}_{\vec{g}}})}}A\cdot B \bigg).
\end{equation}

We have to analyze the residue associated with the wall-crossing components $F_{\vec{g}}$. Using Proposition~\ref{isomophism}, together with the facts that the gluing map $\widetilde{\mathrm{gl}}_k$ has degree $k!$ and that the standard ordering of the ramification points has degree $|\Aut(\eta)|$, we obtain the following expression,
\begin{multline} \label{f1}
\int_{[\widetilde{\mathrm{gl}}_k^*F^{\mathrm{log}}_{\vec{g}, \mathrm{ord}}]^{\mathrm{vir}} } \frac{A\cdot B }{k!\prod_i|\Aut(\eta^{m+i})|} \\ \cdot \mathrm{Res}_{z}\bigg(  \frac{ \prod_{i,j} \eta_j^{m+i}  \cdot \prod^k_{i=1}(\sum_j \mathrm{c}_{\dim-j}(X)(z/k+\Psi(D_i))^j) )}{-z/k-\Psi(D_1)-\Psi_{m+1}-\sum^{\infty}_{i=k}\mathfrak{Y}_i} 
\\ \cdot \prod^k_{i=1} \CI_{g_i,N_i,\eta^{m+i}}(z/k+\Psi(D_i))\bigg). 
\end{multline}
We wish to apply the projection formula with respect to the composition 
\begin{equation} \label{projection}
 \widetilde{\mathrm{gl}}_k^*F^{\mathrm{log}}_{\vec{g},\mathrm{ord}}\rightarrow Y^{\mathrm{log}} \rightarrow \widetilde{M}^{\epsilon_+}_{g_0,n_0,m+k}(X,\beta, \underline{\eta}, \eta^{m+1}, \dots, \eta^{m+k})^{\mathrm{log}}.
 \end{equation}
Before doing that, we organize the expression above. First, by \cite[Lemma 2.7.3]{YZ}, the pullback of $\mathfrak{Y}_i$ to  $\widetilde{\mathrm{gl}}_k^*F^{\mathrm{log}}_{\vec{g},\mathrm{ord}}$ is equal to the pullback of the same boundary divisor $\mathfrak{Y}'_{i-k}$ but  on $ \widetilde{\FM}_{X, m+k,d-kd_0}$. Secondly, if $N_i\neq \emptyset$ for $i\neq 0$, then the insertion $A_{N_i}$ comes from the component $V_{g_i,N_i,\eta^{m+i}}$. We therefore split $A$ according to the partition $(N_0, \dots, N_{k})$,
\[A= \prod^{k}_{i=0}A_{N_i}, \]
and move $A_{N_i}$ to the factors corresponding to $V_{g_i,N_i,\eta^{m+i}}$ in the formula. 
 Finally, we make the change of variables, which affects the residue by the factor of $k$,  
\[ z \rightarrow k(z- \Psi(D_1)-\Psi_{m+1})=\hdots=k(z- \Psi(D_k)-\Psi_{m+k}).\]
Overall, we obtain that (\ref{f1}) is equal to
\begin{multline}\frac{1}{(k-1)!}\int_{[\widetilde{\mathrm{gl}}_k^*F^{\mathrm{log}}_{\vec{g},\mathrm{ord}}]^{\mathrm{vir}} } \frac{\prod_{i,j} \eta_j^{m+i}\cdot A_{N_0}\cdot B}{\prod_i|\Aut(\eta^{m+i})|} \\ \cdot \mathrm{Res}_{z}\bigg(  \frac{\prod^k_{i=1}(\sum_j \mathrm{c}_{\dim-j}(X)(z-\Psi_{m+i})^j) )}{-z-\sum^{\infty}_{i=0}\mathfrak{Y}'_i} 
	\cdot \prod^k_{i=1} A_{N_i}\cdot \CI_{g_i,N_i,\eta^{m+i}}(z-\Psi_{m+i})\bigg).
	\end{multline}

We now apply the projection formula with respect to (\ref{projection}), remembering that the second map in the composition is of degree $1/k$, 
\begin{multline} \label{contributionwall}
\frac{1}{k!} \int_{[\widetilde{M}^{\epsilon_+}_{g_0,N_0,m+k}(X,\beta, \underline{\eta}, \underline{\eta}')^{\mathrm{log}}]^{\mathrm{vir}} } \frac{A_{N_0}\cdot B}{\prod_i |\Aut(\eta^{m+i})|}  \\
\cdot  \mathrm{Res}_z\left( \frac{\prod^{k}_{i=1}\widetilde{\ev}^*_{m+i}	(I_{g_i,N_i,\eta^{m+i}}(z-\Psi_{m+i},A_{N_i}))}{-z-\sum^\infty_{i=0} \mathfrak{Y}'_i} \right),
 \end{multline}
where $\underline{\eta}':=(\eta^{m+1}, \dots, \eta^{m+k})$. 
Note that by our definition of $I$-functions in Definition \ref{Idefinition}, they absorb the factors $\prod^k_{i=1}(\sum_j \mathrm{c}_{\dim-j}(X)(z-\Psi_{m+i})^j) )$ and  $\prod_{i,j} \eta_j^{m+i}$. 

It remaines to determine the contribution of 

\[ \frac{1}{-z-\sum^{\infty}_{i=0}\mathfrak{Y}'_i}=-\frac{1}{z}+ \sum_{s\geq1}\sum_{r\geq 1}(-z)^{-s-1}\mathfrak{Y}'_{r-1}\cdot ( \sum_{i\geq 0} \mathfrak{Y}'_i)^{s-1}.\]
For a Laurent polynomial $f(z)=\sum a_iz^i$, we define $[f(z)]_i=a_i$, and denote \[
C:=\prod^{k}_{i=1}\widetilde{\ev}^*_{m+i}	(I_{g_i,N_i,\eta^{m+i}}(z-\Psi_{m+i},A_{N_i})).
\]By the same argument as in \cite[Lemma 7.3.1]{YZ}, we obtain 
\begin{multline}
	\int_{[\widetilde{M}^{\epsilon_+}_{g_0,N_0,m+k}(X,\beta, \underline{\eta}, \underline{\eta}')^{\mathrm{log}}]^{\mathrm{vir}} } A_{N_0}\cdot B\cdot C \cdot \mathfrak{Y}'_{r-1}\cdot ( \sum_{i\geq 0} \mathfrak{Y}'_i)^{s-1}\\
	=\sum_{\vec{g}_0}\sum_{\underline{j}} \int_{[\widetilde{M}^{\epsilon_+}_{g'_0,N'_0,m+k+r}(X,\beta, \underline{\eta}, \underline{\eta}',\underline{\eta}'')^{\mathrm{log}}]^{\mathrm{vir}}}\frac{(-1)^{s-r} A_{N_0'}\cdot B\cdot C}{r!} \\
	\cdot \prod_{a=1}^{a=r}\widetilde{\ev}^*_{m+k+a}[I_{g_{k+a},N_{k+a},\eta^{m+k+a}}(z-\Psi_{m+k+a},A_{N'_{a}} )]_{-j_a-1},
\end{multline} 
where  the sum is taken over stable partitions $\vec{g}_0$ of $(g_0, N_0)$ of weight $d_0$ and length $r$, and  $\underline{j}=(j_1,\dots, j_r)$ is an $r$-tuple of non-negative integers such that $j_1+\hdots+j_r=s-r$. 
The difference with \cite[Lemma 7.3.1]{YZ} is that the map-theoretic insertions  $A_{N_0}$ have to be taken into account, i.e., they split among the components that appear in the formula above. 

Plugging this expression into (\ref{contributionwall}), taking the residue, and relabelling the marked points, we obtain that the right hand side of (\ref{theformula}) is equal to 
\begin{multline*}
	\sum_{k\geq 1} \sum_{\vv{g}} \sum^{k-1}_{r=0} \sum_{\underline{b}} \frac{(-1)^r}{r!(k-r)!} \int_{[\widetilde{M}^{\epsilon_+}_{g_0,N_0,m+k}(X,\beta, \underline{\eta}, \underline{\eta}')^{\mathrm{log}} ]^{\mathrm{vir}}} \frac{A_{N_0} \cdot B}{\prod_i |\Aut(\eta^{m+i})|} \\ \cdot \prod_{i=1}^{k}\widetilde{\ev}^*_{m+i}[	I_{g_i,N_i,\eta^{m+i}}(z-\Psi_{m+i},A_{N_i})]_{b_\ell}, 
	\end{multline*}
where we sum over integers $k$, stable partitions $\vec{g}$ of $(g, \{1,\hdots, n\})$, and integers $r$; $\underline{b}$ runs through $k$-tuples of integers such that $b_1 + \hdots + b_k=0$, and $b_{k-r+1},\dots, b_k<0$. By noticing that the terms with $b_i<0$ cancel out due to the sign (as explained in the proof of \cite[Theorem 7.3.3]{YZ}), we obtain that the right-hand side of (\ref{theformula}) simplifies to
\begin{multline} \sum_{k\geq 1}\sum_{\vec{g}}  \frac{1}{k!} \int_{[\widetilde{M}^{\epsilon_+}_{g_0,N_0,m+k}(X,\beta, \underline{\eta}, \underline{\eta}')^{\mathrm{log}}]^{\mathrm{vir}}} \frac{A_{N_0} \cdot B}{\prod_i |\Aut(\eta^{m+i})|}\\
	\cdot  \prod_{i=1}^{k}\widetilde{\ev}^*_{m+i}[	I_{g_i,N_i,\eta^{m+i}}(z-\Psi_{m+i},A_{N_i})]_{0}, 
\end{multline}
such that 
\[ 
	[I_{g_i,N_i,\eta^{m+i}}(z-\Psi_{m+i},A_{N_i})]_0=\mu_{g_i,N_i,\eta^{m+i}}(-\Psi_{m+i}, A_{N_i}).
\]
By~(\ref{tildenottilde}), we can instead integrate the above classes over
\[
\Mbar^{\epsilon_+}_{g_0,N_0,m+k}(X,\beta,\underline{\eta},\underline{\eta}')^{\mathrm{log}}.
\]Overall, we derived the following expression for the wall-crossing formula (\ref{theformula}),
\begin{multline}
\int_{[ \Mbar^{\epsilon_-}_{g,n,m}(X,\beta,\underline{\eta})^{\mathrm{log}}]^{\mathrm{vir}}}A\cdot B- \int_{[ \Mbar^{\epsilon_+}_{g,n,m}(X,\beta,\underline{\eta})^{\mathrm{log}}]^{\mathrm{vir}}} A\cdot B \\
=\sum_{\vec{g}} \frac{1}{k!}\int_{[\Mbar^{\epsilon_+}_{g_0,N_0,m+k}(X,\beta, \underline{\eta}, \underline{\eta}')^{\mathrm{log}}]^{\mathrm{vir}} } \frac{A_{N_0} \cdot B}{\prod_i |\Aut(\eta^{m+i})|} \\
 \cdot  \prod^{k}_{i=1}\widetilde{\ev}_{m+i}^*	(\mu_{g_i,N_i,\eta^{m+i}}(-\Psi_{m+i}, A_{N_i})),
\end{multline}
  completing the proof. \qed 



\section{Divisor and dilaton equations} \label{sectiondil}

\subsection{Divisor equations} 
In this section, we will prove dilaton and divisor equations for relative unramified invariants, they will be crucial for evaluating the wall-crossing formula for threefolds in Section \ref{sectionthreefolds}. One can use the wall-crossing to reduce them to the standard dilaton and divisor equations for some divisor classes. However, in the unramified theory, there is an extra class one can insert -- the hyperplane class -- which cannot be compared to anything on the side of the standard theory. Hence, because of the hyperplane class, we will prove the divisor equation in a more direct way by using the forgetful morphisms. 

Like in the standard Gromov--Witten theory, there exists a forgetful morphism between $\Mbar^+_{g,n,m}(X,\beta)^{\log}$ and $\Mbar^+_{g,n-1,m}(X,\beta)^{\log}$. However, unlike in the standard Gromov--Witten theory, the forgetful morphism cannot be identified with the universal curve. The reason is that when two marked points come together (or even when one marked point approaches a point in the fiber of another marked point), the target sprouts an end component with lines attached to the points in the fiber of the marked points. These lines have positive-dimensional moduli whenever there are at least two points in the fiber, including the marked points themselves. Nevertheless, the forgetful morphism still admits a compatible obstruction theory; see Proposition~\ref{forgobst}.

\begin{prop} \label{forg}There exists a forgetful morphism,
	\[\mathrm{fg}_{n} \colon \Mbar^+_{g,n,m}(X, \beta)^{\log} \rightarrow \Mbar^+_{g,n-1,m}(X, \beta)^{\log}.\] 
\end{prop}
\textit{Proof, Geometric points.} 
We start with exhibiting the forgetful morphism on geometric points. For a geometric point $(\Spec(K), N)$, the map $\mathrm{fg}_{n}$ is  constructed as follows.  Let
  \[
  \begin{tikzcd} [row sep=small, column sep = small]
  	(C, \underline{p}, M_{C}) \arrow{dr} \arrow{rr}{f} && (W(m), M_{W(m)}) \arrow{dl} \\
  	& (\Spec(K),N) 
  \end{tikzcd}
  \]
  be a log unramified map.  We forget the $n$-th marking on the source curve $C$, obtaining a log curve $(C, \underline{p}', M_{C})$.  This might destabilize the resulting map  
  \[
  f \colon (C,\underline{p}', M_{C})\rightarrow  (W(m), M_{W(m)}).
  \] We therefore need to contract unstable end and ruled components of $(C, \underline{p}') $ and $W(m)$. On the target, such unstable components will be one of the following three types: 
 \begin{itemize}
\item  end components with only line images and images of exactly two markings from $C$, including the $n$-th marking, 
\item  chains of ruled components attached to unstable end components  with  fiber and non-fiber line images (see Section \ref{FMdeg} for the definition of fiber and non-fiber lines), and without images of markings from $C$,
\item (marked and non-marked) ruled components  with only fiber line images and the image of  exactly one  marking from $C$, namely the  $n$-th marking. 

\end{itemize} 

 We will refer to such components as $n$-\textit{unstable components}.  A pictorial representation of $n$-unstable components is given in Figure \ref{fig:unstable}. An unstable end component together with an unstable ruled component are shown on the left, while an unstable marked ruled component is on the right.

Before forgetting the $n$-th marking, such components can be characterized solely in terms of weights and markings:
 \begin{itemize}
 \item unstable end components are end components with exactly two markings from $C$, including the $n$-th marking, whose weight with respect to the class $\big[\Omega^\mathrm{log}_C\big(\sum_{j,i} p_{j,i}\big) \big] -\big[f^*\Omega^{\mathrm{log}}_{W(m)}\big(\log\big(\sum_{j} E_j \big)\big)\big]/\dim$ is 0,\footnote{More precisely, the degree of  $\big[\Omega^\mathrm{log}_C\big(\sum_{j,i} p_{j,i}\big) \big] -\big[f^*\Omega^{\mathrm{log}}_{W(m)}\big(\log\big(\sum_{j} E_j \big)\big)\big]/\dim$  over such a component is 0; this follows from Lemma \ref{weightP} and Lemma \ref{less}.}
 
 \item unstable ruled components of the first type are ruled components without images of markings from $C$, attached to an unstable end component via an unstable chain, whose  weight with respect to the class $\big[\Omega^\mathrm{log}_C\big(\sum_{j,i} p_{j,i}\big) \big] -\big[f^*\Omega^{\mathrm{log}}_{W(m)}\big(\log\big(\sum_{j} E_j \big)\big)\big]/\dim$ is 0, 
 \item unstable (marked and non-marked) ruled components of the second type  are   ruled components with the image of one marking from $C$, namely, the $n$-th marking, whose weights with respect to the classes $\big[\Omega^\mathrm{log}_C\big(\sum_{j,i} p_{j,i}\big) \big] -\big[f^*\Omega^{\mathrm{log}}_{W(m)}\big(\log\big(\sum_{j} E_j \big)\big)\big]/\dim$ and  $-\big[f^*\Omega^{\mathrm{log}}_{W(m)}\big(\log\big(\sum_{j} E_j \big)\big)\big]$ are 0.\footnote{This follows from the fact that for a marked (or non-marked) ruled component $P_\mathrm{rul}$, the class $\mathrm{c}_1\big(\Omega^{\mathrm{log}}_{W(m)}\big(\log\big(\sum_{j} E_j \big)\big) \big|_{P_\mathrm{rul}}\big)= \dim \cdot ([E_j]-[D])$ pairs trivially with a curve, if and only if it is a multiple of the fiber line class.  Note that it is not zero on a union of fiber and non-fiber lines (which do not become destabilised after forgetting a marking); this is why we have to use an extra class to distinguish such components.}
  \end{itemize}

\begin{figure}[!ht]
	\centering
	\vspace{0.15cm}
	\begin{tikzpicture}		
		\node (0) at (1,0.5) {\large$\mathbb{P}^{\dim}$};
		\node (1) at (3.5,0.5) {\large$\mathrm{Bl}_0(\mathbb{P}^{\dim})$};
		\node (2) at (0.7,-0.55) {$\mathbb{P}^1$};
		\node (3) at (3,-3) {};

		\draw[black, thick] (0,0)--(2,0);
		\draw[black, thick] (0,-3)--(2,-3);
		\draw[black, thick] (2,-3)--(2,0);
		\draw [black,thick] plot [smooth,tension=1.5] coordinates {(0,0)  (-1,-1.5) (0,-3)};
		
		\draw[black, thick, opacity = 0.3] (2,-0.75)--(0,-1);
		\draw[black, thick, opacity=0.3] (2,-1.5)--(0,-1.75);
		\draw[black, thick, opacity=0.3] (4.5,-2.25)--(2.5,-2.5);
		
		\node (5) at (0.5,-1.7) {};
		\node (6) at (1.25,-0.85) {};
		\filldraw[thick, fill = black] (0.5,-1.7) circle (.05cm) node at (5) {};
		\filldraw[thick, fill = black] (1.25,-0.85) circle (.05cm) node at (6) {};
		
		\draw[black, thick, opacity=0.3] (2,-0.75)--(4.5,-0.75); 
		\draw[black, thick, opacity=0.3] (2,-1.5)--(4.5,-1.5);

		\draw[black, thick] (2,0)--(4.5,0);
		\draw[black, thick] (2,-3)--(4.5,-3);
		\draw[black, thick] (4.5,0)--(4.5,-3);
		
		\draw[black, thick, dashed] (4.5,0)--(5,0);
		\draw[black, thick, dashed] (4.5,-3)--(5,-3);
		\draw[black, thick, opacity=0.3, dashed] (4.5,-0.75)--(5,-0.75); \draw[black, thick, opacity=0.3, dashed] (4.5,-1.5) --(5,-1.5);  
		\draw[black, thick, opacity=0.3, dashed] (4.5,-2.25) --(5,-2.25);

		\draw[black, thick] (8,0)--(11,0);
		\draw[black, line width=1mm] (8,0)--(8,-3);
		\draw[black, thick] (8,-3)--(11,-3);
		\draw[black, thick] (11,0)--(11,-3);
		
		\draw[black, thick, opacity=0.3] (8,-1.5) --(11,-1.5);  
		
		\node (7) at (9.5,-0.75) {};
		\filldraw[thick, fill = black] ((9.5,-1.5) circle (.05cm) node at (7) {};
		
		\draw[black, thick, dashed] (11,0)--(11.5,0);
		\draw[black, thick, dashed] (11,-3)--(11.5,-3);
		\draw[black, thick, opacity=0.3, dashed] (11,-1.5) --(11.5,-1.5);  
		
		\node (9) at (9.5,0.5) {\large$\mathrm{Bl}_0(\mathbb{P}^{\dim})$};
		\node (9) at (7.7,0.45) {\large$E_j$};
		\node (10) at (9,-1.15) {$\mathbb{P}^1$};

	\end{tikzpicture}
	
	\caption{Depiction of $n$-unstable components.}
	\label{fig:unstable}

\end{figure}

Let 
\begin{equation*} \label{contr0}
	\pi_1 \colon (C,\underline{p}',M_C) \rightarrow (C',\underline{p}',M_{C'}), \quad  \pi_2 \colon  (W(m),M_{W(m)})\rightarrow  (W(m)',M_{W(m)'})
\end{equation*}
 be the contractions of $n$-unstable components  together with the unstable components in the source after forgetting the $n$-th marking, such that log structures on stabilisations are defined as follows,
 \[ M_{C'}:=\pi_{1*}(M_C), \quad M_{W(m)'}:=\pi_{2*}(M_{W(m)}),\]  which indeed are log structures by \cite[Lemma B.3]{AMW}.  By  \cite[Theorem B.6]{AMW}, the map $f$ descends to a map between $(C',\underline{p}',M_{C'})$ and $(W(m)', M_{W(m)'})$ over $(\Spec(K),N)$. However, the minimality of log structures is violated, as  $N$ will contain  copies of $\BN$  that correspond to simple normal crossing singularities of $W(m)$ that disappeared after the contraction. One can replace $N$ by the associated minimal log structure $N'$, which is defined by removing copies of $\BN$ as follows.  By \cite[Section 4.1]{KiL}, we have an identification 
 \[N\cong \BN^{k} \oplus K^*,\]
 where $k$ is the number of connected components of  $W(m)^{\mathrm{sing}}$.\footnote{For unramified maps, there are no non-distinguished nodes in Kim's terminology, as nodes are not allowed to map to the smooth locus of the target by Definition \ref{stab} (2).} Now let 
 $N'\subseteq N$ 
 to be the submonoid that corresponds to the connected components of $W(m)'^{\mathrm{sing}}$, such that for a (non-marked) ruled component of the second type, we consider a diagonal embedding $\BN \rightarrow \BN^2$, where $\BN^2$ is associated to two neighbouring components. By \cite[Section 4.1]{KiL}, this is independent of the identification above, as there is a bijection between irreducible elements of $\BN^k$ and connected components of $W(m)^{\mathrm{sing}}$.
   By construction of $N'$ and the minimality of $M_C$, the resulting log curve 
   \[
   (C',\underline{p}',M_{C'}) \rightarrow (\Spec(K),N') 
   \] is minimal. Other required properties of log structures from   \cite[Section 5.2]{KiL} continue to hold, since they are local and the contraction left the remaining components untouched.   We thereby obtain a log unramfied map 
 \[
 \begin{tikzcd}[row sep=small, column sep = small]
 	(C',\underline{p}',M_{C'}) \arrow{dr} \arrow{rr}{f'} && (W(m)', M_{W(m)'}) \arrow{dl} \\
 	& (\Spec(K),N') 
 \end{tikzcd}
 \]
See \cite[Section 5.2.3]{KiL} for more details on how log unramified maps look like in local charts.
We now show how to extend this construction to families.
\vspace{0.3cm}

\textit{Proof, Families.} The construction extends to families by the existence of the forgetful morphisms for Fulton--MacPherson compactifications established in \cite[Section~2.6]{KKO}. We begin by constructing the contraction of $n$-unstable components for the moduli space of \emph{biweighted} extended log twisted FM degenerations $\FM^{\log}_{X,m,n,\underline{d}}$, defined below.  As we shall see, this  is possible because $n$-unstable components are characterized solely by their weights and markings. Two weights are needed because $n$-unstable ruled components of the second type are characterized by an additional weight. The smoothness of $\FM^{\log}_{X,m,n,\underline{d}}$ will allow us to descend  the construction of the contraction morphism from a smooth cover.  

Let 
\[
\FM^{\log}_{X,m,n,\underline{d}}  \colon (Sch/\BC)^\circ \rightarrow Grpd
\]
be the moduli space of  FM degenerations with $n$ (non-blown-up) markings and $m$ markings corresponding to exceptional divisors (equivalently, to blown-up markings), such that they are equipped with two weights assigned to each irreducible component,
\[
\underline{d}=(\underline{d}_1, \underline{d}_2).
\] Analogously to Definition~\ref{moduliwieght}, a FM degeneration of $\FM^{\log}_{X,m,n,\underline{d}}$ and its smoothings are required to satisfy the following conditions:
\begin{itemize}
	 \item[(1)]  the number of irreducible components and their weights are bounded,
 \item[(2)] weights of end components and ruled components of a FM degeneration  are non-negative, 
 \item[(3)] if the weight $\underline{d}_1$ of an end component is zero, then it must contain at least two markings.
 \item[(4)] the weight $\underline{d}_2$ of an end component is positive,
 \item[(5)] if the weight $\underline{d}_2$ of a ruled component is zero, then the weight $\underline{d}_1$  is positive or it contains at least one marking,
 \item[(6)] if both $\underline{d}_1$ and $\underline{d}_2$ are zero on a ruled component, this ruled component cannot degenerate further into a non-ruled one. 


 \end{itemize}
Note that since we require these conditions to be satisfied for all smoothings of a FM degeneration,  they are open in the space of biweighted FM degenerations. 
See (\ref{biweighted}) for what the weights $(\underline{d}_1,\underline{d}_2)$ are in the case of maps; the conditions above are tailored to have the properties of the weights given by the degrees of log cotangent bundles.  

Choose a FM degeneration  with $n$ markings in $\FM^{\log}_{X,m,n,\underline{d}}$, 
\[ (W(m), \underline{x}):=(W(m), x_1,\dots, x_{n}).\] 
Next, choose two sets of auxiliary markings in the smooth locus of $W(m)$,
 \[
 \underline{\sigma}_u \quad \text{and} \quad \underline{\sigma}_s,
 \] 
 disjoint from $\underline{x}$, such that $\underline{\sigma}_u$ and the $n$-marking stabilise\footnote{This means that the automorphism group of a component fixing the markings on it is trivial, while without the markings $\underline{\sigma}_u$ and the $n$-th marking, it is not.} $n$-unstable components of $(W(m), \underline{x})$.  While the markings  $\underline{\sigma}_s$ and the remaining markings from $\underline{x}$ stabilise other end and ruled components. Since the markings are contained in the smooth locus of $W(m)$, we can spread 
 these markings to the universal family  $\mathfrak{W}(m)$ over a smooth schematic chart
 \[
 B_i \rightarrow  \FM^{\log}_{X,m,n,\underline{d}}
 \] around $(W(m), \underline{x})$, which satisfies the following conditions: 
\begin{itemize}
\item{} $B_i$ contains only smoothings\footnote{In other words, we do not permit further degenerations of $(W(m), \underline{x})$.} of the singularities of $(W(m), \underline{x})$, and markings $\underline{\sigma}_u$ and $\underline{\sigma}_s$ are disjoint from $\underline{x}$,
\item{} for all FM degenerations in $B_i$, the markings $\underline{\sigma}_u$  stabilise \textit{only} $n$-unstable  components (which are $n$-unstable end and ruled components of the first type with respect to the first weight, and $n$-unstable ruled components of the second type with respect to the second weight). 
\end{itemize}

By the definition of $\FM^{\log}_{X,m,n,\underline{d}}$, the above conditions indeed define a smooth chart, i.e., they define an open subset of some Fulton--MacPherson chart $X[n_i]$ of $\FM^{\log}_{X,m,n,\underline{d}}$. As we explain below, to see this, it is enough to inductively smooth chains of ruled components together with end components and verify, using the properties of the weights established in Lemma \ref{lemmaweights}, that the resulting FM degenerations satisfy the same conditions.

First, we may assume that $(W(m),\underline{x}) \in B_i$ already contains an $n$-unstable component, because, by the defining conditions, $n$-unstable components cannot be created as a result of smoothing a FM degeneration without $n$-unstable components; for non-marked ruled components of the second type, this is prohibited by the open condition (6). Next, if a singularity between an  $n$-unstable component and a non-$n$-unstable component is smoothed, then the resulting component is either an end or a ruled component of positive weight, or it contains other markings from $\underline{x}$ or blown-up markings, or it is a non-end and non-ruled component. In the latter case, if we further smooth the trees of components attached to it, then by the positivity assumption on the weights in $\FM^{\log}_{X,m,n,\underline{d}}$, it either results in an end component or a ruled component of positive weight, or it contains other markings from $\underline{x}$ or blown-up markings, making it non-$n$-unstable.  This implies that the markings $\underline{\sigma}_s$ cannot stabilise a $n$-unstable component after smoothings (and vice versa),  which in turn implies that the above conditions are indeed open.   Moreover, by the choice of  $\underline{\sigma}_u$ and $\underline{\sigma}_s$, the first defining condition of $B_i$  above implies the second. 

We  identify $\mathfrak{W}(m)|_{B_i}$ with the restriction of the (blow-up of) the universal family $\FW[n_i]|_{B_i}$ for some Fulton--MacPherson chart $X[n_i]$, using the auxiliary markings and the markings $\underline{x}$. We then forget the $n$-th marking  and the auxiliary markings $\underline{\sigma}_u$, and apply the forgetful morphisms from  \cite[Section 2.6]{KKO} to the family $\mathfrak{W}[n_i]|_{B_i}$.  Note that by construction of $B_i$, the result of the forgetful morphism applied to $\mathfrak{W}[n_i]|_{B_i}$ will be the contraction of exactly $n$-unstable components. We define log structures  by taking pushforwards of log structures with respect to contractions, while on the base $B_i$ we take the minimal log structure by removing copies of $\BN$ in $N_{B_i}$ corresponding to the contracted singularities like in the case of the construction presented for geometric points; this is well-defined because by \cite[Section 4.1]{KiL}, in a logarithmic chart, there is a bijection between irreducible elements of the logarithmic structure and  components  of the singular locus. 
 We therefore obtain
\[
\mathrm{fg}_{n}|_{B_i} \colon B_i \rightarrow  \FM^{\log}_{X,m,n-1,\underline{d}}.
 \]
 We can cover $\FM^{\log}_{X,m,n,\underline{d}}$ by such smooth charts, 
 \[\coprod_i B_i \rightarrow \FM^{\log}_{X,m,n,\underline{d}}.\]  
 
  The resulting forgetful morphisms on each $B_i$ do not depend on the choice of markings. This can be seen by using the space $X[n_1,n_2]$ from \cite[Section 2.7]{KKO} and \cite[Proposition 2.7.2]{KKO}. The space $X[n_1,n_2]$ parametrises stable FM degenerations with two sets of possibly overlapping markings.  Assume  that over a chart $B_i$, we have two different sets of markings,
  \[\underline{\sigma}_{u,1}\quad \text{and} \quad \underline{\sigma}_{u,2},\] that stabilize components we want to contract in the way it is described above. We define two maps, 
  \[
  g_j \colon B_i \rightarrow X[n'_{1},n'_{2}], \quad j=1,2,
  \]
  by applying the forgetful morphisms to $\underline{\sigma}_{u,1}$ and then to $\underline{\sigma}_{u,2}$, together with the $n$-th marking,  and vice versa. 
Since the moduli spaces of FM degenerations are smooth, and in particular, reduced,  a chart $B_i$ is reduced as well.  It is clear that $g_j$ agree on geometric points, hence they must be equal by the equality of the associated graphs. This implies that the resulting families given by contractions associated to markings $\underline{\sigma}_{u,1}$ and $\underline{\sigma}_{u,2}$ are 
 isomorphic; the same applies to the log structures, since they are given by the pushforwards of the log structures via the contraction maps on the total family, and by the removal of copies of $\BN$ on the base.   In other words, forgetful morphisms are isomorphic on double intersections  
  \[
 B_i \times _{\FM^{\log}_{X,m,n,\underline{d}}} B_j.
 \]
 
  By using a similarly defined space $X[n_{1},n_{2},n_3]$, we also conclude that the  forgetful  morphisms  $\mathrm{fg}_{n}|_{B_i}$ satisfy the cocycle condition on triple intersections 
  \[
  B_i \times _{\FM^{\log}_{X,m,n,\underline{d}}} B_j\times _{\FM^{\log}_{X,m,n,\underline{d}}} B_k.
  \] 

Hence, the forgetful morphisms on $B_i$ glue, or, in other words, descend from the smooth cover $\coprod_i B_i$.  We  therefore obtain the forgetful morphism for the moduli space of log FM degenerations in families, 
\begin{equation} \label{theformap}
\mathrm{fg}_{n} \colon \FM^{\log}_{X,m,n,\underline{d}} \rightarrow  \FM^{\log}_{X,m,n-1,\underline{d}}.
\end{equation}

We now proceed to the construction of the forgetful morphism for 	$\Mbar^+_{g,n,m}(X, \beta)^{\log}$.  By associating to an unramified map the weights given by 
 \begin{align} \label{biweighted}
 	\begin{split} &d_{1,\ell}=\deg\big(\Omega^\mathrm{log}_C\big(\sum_{j,i} p_{j,i}\big)|_{W(m)_\ell} \big) -\deg\big(f^*\Omega^{\mathrm{log}}_{W(m)}\big(\log\big(\sum_{j} E_j \big)\big)|_{W(m)_\ell}\big)/\dim, \\
& d_{2,\ell}= - \deg\big(f^*\Omega^{\mathrm{log}}_{W(m)}\big(\log\big(\sum_{j} E_j \big) \big) \big|_{W(m)_\ell}\big),
\end{split}
\end{align}
we obtain a morphism

\begin{equation*} 
	\Mbar^+_{g,n,m}(X, \beta)^{\log} \rightarrow \FM^{\log}_{X,m,n,\underline{d}},
\end{equation*}
where the additional $n$ markings correspond to the images of the markings from curves.  Note that we do not include markings in our weights unlike in (\ref{weights}), but remember them in our moduli spaces instead. The fact that the weights above satisfy conditions (1)--(5) of $\FM^{\log}_{X,m,n,\underline{d}}$ follows from Lemma \ref{weightP} and the arguments from Lemma \ref{weights}.  Condition (6) is satisfied because maps covering fiber-lines in a ruled component, with no ramifications in the smooth locus and containing only one marking, cannot degenerate further, since their moduli space is a point.

We now apply the forgetful morphism (\ref{theformap}) to FM degenerations  in the space $\Mbar^+_{g,n,m}(X, \beta)^{\log}$, stabilizing the composition of the contraction with the universal map. Such a stabilisation is more standard and is achieved by using the linear system associated to the pullback of a relative ample line bundle $L$ from FM degenerations twisted by the log canonical line bundle of the source, $f^*L^{\otimes3}\otimes  \Omega^{\mathrm{log}}_{(C,\underline{p})}$. The line bundle $L$ exists in smooth charts, while the stabilization is independent of its choice and therefore descends from a smooth cover. We use \cite[Appendix~B]{AMW} to obtain a natural morphism between the log structures of the contracted targets and sources. By the same arguments as in the case of geometric points, these contractions preserve the defining properties of log unramified maps in the sense of \cite[Section~5.2]{KiL}. Overall, we obtain the desired forgetful morphism, 
 \[\mathrm{fg}_{n} \colon \Mbar^+_{g,n,m}(X, \beta)^{\log} \rightarrow \Mbar^+_{g,n-1,m}(X, \beta)^{\log},\]
 this completes the proof. 
  \qed
 
 \begin{prop} \label{forgobst} The forgetful morphism	\[\mathrm{fg}_{n} \colon \Mbar^+_{g,n,m}(X, \beta)^{\log} \rightarrow \Mbar^+_{g,n-1,m}(X, \beta)^{\log}\] 
 	admits a compatible relative perfect obstruction theory in the sense of \cite[Defintion 4.5]{Man2}.
 	\end{prop}
 
\textit{Proof.} For brevity, we denote 
\[M_1=\Mbar^+_{g,n,m}(X, \beta)^{\mathrm{log}}, \quad  M_2=\Mbar^+_{g,n-1,m}(X, \beta)^{\mathrm{log}}.\] 
A compatible relative perfect obstruction theory is a perfect 2-term complex $\BL^{\mathrm{vir}}_{M_1/M_2}$ with a morphism to the cotangent complex (or its truncation), 
\[ \BL^{\mathrm{vir}}_{M_1/M_2}  \rightarrow \BL_{M_1/M_2} ,\]
which is isomorphism in degree 0 and surjection in degree $-1$, and which fits in a triangle
\[\mathrm{fg}_{n}^*\BL^{\mathrm{vir}}_{M_2}  \rightarrow \BL^{\mathrm{vir}}_{M_1}  \rightarrow \BL^{\mathrm{vir}}_{M_1/M_2},\]
respecting the obstruction-theory morphisms. By \cite[Construction 3.13]{Man2}, in order to construct a compatible relative perfect obstruction theory, it is enough to construct a morphism  
\begin{equation} \label{map0}
	\mathrm{fg}_{n}^*\BL^{\mathrm{vir}}_{M_2} \rightarrow  \BL^{\mathrm{vir}}_{M_1}, 
	\end{equation}
	that respects the obstruction-theory morphisms, whose cone is a 2-term complex.

 Since we contract components both on the target and on the source,  we have to work with the absolute obstruction theory of maps. 
Let $\FM^{\log}_{g,n+n'}$ be the moduli space of minimal log prestable curves. Extra markings  correspond to intersection points of curves with the exceptional divisors.  Consider the stacks
\begin{align*}
\mathfrak{S}_1&=\FM^{\log}_{X,m,n,\underline{d}} \times_{\log} \FM^{\log}_{g,n+n'}, \\
 \mathfrak{S}_2&=\FM^{\log}_{X,m,n-1,\underline{d}} \times_{\log} \FM^{\log}_{g,n+n'-1},
\end{align*}
 defined for example in \cite[Section 6.2.3]{KiL} (the same construction applies to the weighted versions of the moduli spaces). We have natural projections
\[ \tau_i\colon  M_i \rightarrow\mathfrak{S}_i, \quad i=1,2,  \]
such extra markings on the target are given by the images of markings on the source. Using the same argument as before, we can define the forgetful morphism  also between  moduli spaces $\mathfrak{S}_1$ and $\mathfrak{S}_2$, 
\[ \mathrm{fg}_{n} \colon \mathfrak{S}_1 \rightarrow \mathfrak{S}_2, \]
by forgetting the $n$-th marking and  contracting the $n$-unstable components both on the target and the source\footnote{By taking the log product in the definition of $\mathfrak{S}_i$, we know which components of the curves lie over which components of the FM degenerations. This information is encoded in the monoid $N$ of the base $\operatorname{Spec}(B)$; namely, we know how the intersection graphs of the curves map to the intersection graphs of the FM degenerations. We can therefore contract $n$-unstable components in the target together with the corresponding unstable components in the source simultaneously, using the same argument as in Proposition~\ref{forg}.}.  We therefore obtain the following commutative diagram,    

\begin{equation}\label{diagram}
	\begin{tikzcd}[row sep=scriptsize, column sep = scriptsize]
		& M_1 \arrow[d,"\tau_1"]  \arrow[r,"\mathrm{fg}_{n}"] &  M_2 \arrow[d,"\tau_2"]\\
		&  \mathfrak{S}_1 \arrow[r,"\mathrm{fg}_{n} "] &   \mathfrak{S}_2.
	\end{tikzcd}
\end{equation}


By Theorem \ref{obst}, the $\tau_i$-relative obstruction theories are given by the complexes
\[\BL^{\mathrm{vir}}_{\tau_i}=(R\pi_*F^*T^{\mathrm{log}}_{\pi, \mathfrak{E}})^\vee , \quad i=1,2  \]
where the notation is as in Theorem~\ref{obst}.  Since the stack $\mathfrak{S}_i$ is smooth, using the construction from  \cite[Proposition 3]{KKP}, the absolute obstruction theory of $M_i$ is therefore given by the cone, 
\begin{equation} \label{cone}
\BL_{M_i}^\mathrm{vir}=\mathrm{cone}(\BL^{\mathrm{vir}}_{\tau_i}[-1] \rightarrow \tau_i^*\BL_{\mathfrak{S}_i}), \quad i=1,2.
\end{equation}
 Hence to construct a morphism between absolute obstruction theories, it is enough to  construct a morphism between relative obstruction theories
\begin{equation} \label{map2}
\mathrm{fg}_{n}^*\BL^{\mathrm{vir}}_{\tau_2}\rightarrow \BL^{\mathrm{vir}}_{\tau_1},
\end{equation}
 such that the following diagram commutes 
\begin{equation}\label{diagram2}
\begin{tikzcd}[row sep=scriptsize, column sep = scriptsize]
	&\mathrm{fg}_{n}^*\BL^{\mathrm{vir}}_{\tau_2}[-1] \arrow[d]  \arrow[r] & \mathrm{fg}_{n}^*( \tau_2^*\BL_{\mathfrak{S}_2})\arrow[d]\\
	&\BL^{\mathrm{vir}}_{\tau_2}[-1] \arrow[r] &  \tau_1^*\BL_{\mathfrak{S}_1},
\end{tikzcd}
\end{equation}
where the second vertical arrow is given by the pullback of $\mathrm{fg}_{n}^*\BL_{\mathfrak{S}_2}\rightarrow \BL_{\mathfrak{S}_1}$ and the commutativity of the diagram (\ref{diagram}).

For simplicity, we exhibit the construction on geometric points; it easily extends to families. Consider the image  
\[
[f' \colon (C',\underline{p}', M_{C'}) \rightarrow  (W(m)',M_{W(m)'})]
\]
 of  an unramified map 
 \[
[ f \colon (C,\underline{p},M_{C}) \rightarrow  (W(m),M_{W(m)})]
 \] with respect to the forgetful morphism 
$\mathrm{fg}_{n}$, and consider the contraction morphisms from (\ref{contr0}), 
\[
\begin{tikzcd}[row sep=scriptsize, column sep = scriptsize]
	& C \arrow[d,"\pi_1"]  \arrow[r,"f"] &  W(m)  \arrow[d,"\pi_2"]\\
	& C' \arrow[r,"f'"] &  W(m) '.
\end{tikzcd}
\]
We have a natural morphism between log tangent bundles
\[
	T^\mathrm{log}_{W(m)} \rightarrow  \pi_2^* T^\mathrm{log}_{W(m)'}.
\]
By pulling back this morphism and applying the commutativity of the diagram, we obtain 
\begin{equation} \label{map1}
f^*T^\mathrm{log}_{W(m)} \rightarrow \pi_1^*(f')^*T^\mathrm{log}_{W(m)'}.
\end{equation} 
Note that since $\pi_1$ involves contractions of end and ruled components (rational tails and bridges), we have 	
\[R\pi_{1*}\CO_C= \CO_{C'}.\] 
Hence, by the projection formula, we have a natural identification of cohomologies,
\[H^*(C,\pi_1^*(f')^*T^\mathrm{log}_{W(m)'})= H^*(C',(f')^*T^\mathrm{log}_{W(m)'}).\]
Overall, by (\ref{map1}), we obtain 
\begin{equation} \label{map3}
	H^*(C,f^*T^\mathrm{log}_{W(m)}) \rightarrow H^*(C',(f')^*T^\mathrm{log}_{W(m)'}).
\end{equation}
Dualizing and  extending to families, we obtain the desired morphism, 
\[\mathrm{fg}_{n}^*\BL^{\mathrm{vir}}_{\tau_2}\rightarrow \BL^{\mathrm{vir}}_{\tau_1}.  \]

  The fact that the obstruction-theory morphisms commute with the morphism above follows from the construction of the obstruction-theory morphisms from \cite[Section 7]{KiL} (the obstruction-theory morphisms are induced by the canonical morphisms between cotangent bundles of the target and source curves, which commute with (\ref{map1})). The diagram (\ref{diagram2}) commutes because  the obstruction-theory morphisms commute and the same canonical diagram for cotangent complexes commutes.  Overall, we obtain a morphism of absolute obstruction theories
  \[ \mathrm{fg}_{n}^*\BL^{\mathrm{vir}}_{M_2} \rightarrow  \BL^{\mathrm{vir}}_{M_1}.\] 
It remains to show that the cone of this morphism is 2-term. It is enough to show it over closed points. Using  (\ref{cone}), and the fact that $\mathfrak{S}_i$ is smooth, we therefore have  to show that the morphism
(\ref{map3}) is surjective.
By construction, it is given by applying the cohomology functor to the morphism
\[f^*T^{\mathrm{log}}_{ W(m)} \rightarrow  \pi_1^*(f')^* T^{\mathrm{log}}_{ W(m)'} .\] 
Let $Q$ be the cokernel of this morphism. The sheaf $Q$ is zero-dimensional on the component of $C$ corresponding to $C'$, while on contracted rational tails it is a trivial locally free sheaf. Hence its first cohomology is zero, 
\[H^1(C,Q)=0,\]
this implies the surjectivity, concluding the proof of the proposition. \qed 

\begin{prop} \label{degreegenus} Let $\p(f) \colon C \rightarrow \p(TX)$ be the map associated to an unramified map $f \colon C \rightarrow  W$ and $H \in H^2(\p(TX))$ be the relative hyperplane class, then
	\[  \deg(\p(f)^*H)=2g(C)-2. \]
\end{prop}
\textit{Proof.}   We have to analyze the degree of $H$ on each component of $C$. Away from the nodes, by the construction of $\p(f)$, the pullback of the relative tautological line bundle $\CO_{\p(TX)}(1)$ is isomorphic to the cotangent line bundle $\Omega_C$. After passing to the normalisation, the pullback of $\CO_{\p(TX)}(1)$ acquires a twist over the preimages of nodes, depending on whether a preimage maps to an exceptional divisor or a divisor at infinity. Determining this twist is a local problem.  Hence, we can assume that $ W$ is just a one-step FM degeneration of $X$, i.e.,\ $ W$ has two components, one is a blow-up of $X$ at a point $\mathrm{Bl}_x(X)$, another one is $\p(T_xX \oplus \BC)$.   Let 
\[ \p(f_{\mathrm{main}}) \colon C_{\mathrm{main}} \rightarrow \p(TX)\]
be the projectivised  derivative associated to the restriction of $f$ over $\mathrm{Bl}_x(X)$, constructed in  \cite[Lemma 3.2.4]{KKO}.  After projecting to $X$, we may view $ \p(f_{\mathrm{main}})$ as the map associated to 
\[ f_{\mathrm{main}} \colon C_{\mathrm{main}} \rightarrow X,\]
which might be ramified only over the point $x \in X$. Let $\eta=(\eta_1,\dots, \eta_\ell)$ be the ramification profile of $f_{\mathrm{main}}$ over  $x$, and $\{p_1,\dots, p_\ell \}$ be the fiber of $f$ over $x$, then the morphism
\[ f_{\mathrm{main}}^* \Omega_X \rightarrow  \Omega_{C_{\mathrm{main}}}(- \sum_i(\eta_i-1) p_i) \]
is surjective. Hence, by the universal property of $\p(TX)$, the map  $\p(f_{\mathrm{main}})$ is associated to the injective morphism of vector bundles
\[  T_{C_{\mathrm{main}}}(\sum_i( \eta_i-1) p_i) \rightarrow f_{\mathrm{main}}^*T_X. \]
By the universal property of the tautological line bundle $\CO_{\p(TX)}(-1)$ on $\p(TX)$, we obtain that 
\[\p(f_{\mathrm{main}})^*\CO_{\p(TX)}(-1)=T_{C_{\mathrm{main}}}(\sum_i( \eta_i-1) p_i),\]
hence 
\begin{equation} \label{main}
	\p(f_{\mathrm{main}})^*H=2g(C_{\mathrm{main}})-2- \sum_i(\eta_i-1).
\end{equation}

Now  let 
\[ \p(f_{\mathrm{end}}) \colon C_{\mathrm{end}} \rightarrow \p(T_xX)\]
be the map associated to the restriction $\p(f_{\mathrm{end}})  \colon C_{\mathrm{end}} \rightarrow \p(T_xX \oplus \BC)$ of $f$ over the end component $\p(T_xX \oplus \BC)$.  Let $U:=T_xX$ to be the complement  of $\p(T_xX) \subset \p(T_xX \oplus \BC)$,  then over $U$ the map  $\p(f_{\mathrm{end}})$ is associated to the injective map of vector bundles, 
\begin{equation} \label{U}
	T_{C_{\mathrm{end}}}|_{f^{-1}(U)} \rightarrow T_xX \otimes \CO_{f^{-1}(U)}.
\end{equation}
We have to analyse what happens over $\p(T_xX) \subset \p(T_xX \oplus \BC)$. By the admissibility of $f$, the ramification profile of $f_{\mathrm{end}}$ over $\p(T_xX)$ is also given by $\eta=(\eta_1,\dots, \eta_{\ell(\eta)})$. Since it is a local question, we may restrict to a formal neighbourhood of a point $p_i \in C_{\mathrm{end}}$ mapping to $\p(T_xX)$. We also choose an identification $T_xX \cong \BC^{\dim}$. In the formal neighbourhood of the point $p_i \in C_{\mathrm{end}}$, the map $f_{\mathrm{end}}$ is of the form 
\[ t \mapsto \left( \frac{1}{t^{\eta_i}},\frac{z_2(t)}{t^{\eta_i}}, \dots, \frac{z_{\dim}(t)}{t^{\eta_i}}  \right).\]
The derivative of this map is 
\[ \left( -\frac{\eta_i}{t^{\eta_i+1}},\frac{z'_2(t)}{t^{\eta_i}}-\frac{\eta_i z_2(t)}{t^{\eta_i+1}}, \dots, \frac{z'_{\dim}(t)}{t^{\eta_i}}-\frac{\eta_i z_{\dim}(t)}{t^{\eta_i+1}} \right).\] 
The leading order of the poles is $\eta_i+1$, hence, we conclude that the morphism (\ref{U}) extends to the following injective morphism of vector bundles, 
\[ T_{C_{\mathrm{end}}}(-\sum_i(\eta_i+1)p_i) \rightarrow T_xX \otimes \CO_{C_{\mathrm{end}}},\]
which therefore induces the map $\p(f_{\mathrm{end}})\colon C_{\mathrm{end}} \rightarrow \p(T_xX)$, such that
\begin{equation*}
	\p(f_{\mathrm{end}})^*\CO_{\p(T_xX)}(-1)= T_{C_{\mathrm{end}}}(-\sum_i(\eta_i+1)p_i),
\end{equation*}
which implies that 
\begin{equation}  \label{end}
	\p(f_{\mathrm{end}})^*H=2g(C_{\mathrm{end}})-2+ \sum_i(\eta_i+1).
\end{equation}
Adding the degrees from the main and end components, we obtain 
\[\p(f)^*H=(\ref{main})+(\ref{end})=2g(C)-2,\]
which is the desired equality. With the same analysis, the proof extends to a FM degeneration with an arbitrary number of irreducible components. 
\qed
\\

We are now ready to prove divisor equations for relative unramified Gromov--Witten theory. 
\begin{prop} \label{dilatondivisor} Relative unramified Gromov--Witten invariants satisfy the following divisor equations:
	\begin{itemize}
		\item if $D \in H^2(X) \subset H^2(\p(TX))$, then 
		\[\langle  \gamma_1 \cdot \! \cdot \! \cdot  \gamma_n \cdot D\mid \gamma_1' \cdot \! \cdot \! \cdot  \gamma_m'  \rangle^+_{g,\beta}=(D \cdot \beta) \langle  \gamma_1 \cdot \! \cdot \! \cdot  \gamma_n  \mid \gamma_1' \cdot \! \cdot \! \cdot  \gamma_m'   \rangle^+_{g,\beta},\] 
		\item	if $H \in H^2(\p(TX))$ is the relative hyperplane class,  then 
		\[\langle  \gamma_1 \cdot \! \cdot \! \cdot  \gamma_n \cdot  H \mid \gamma'_1 \cdot \! \cdot \! \cdot  \gamma'_m  \rangle^+_{g,\beta}=(2g-2) \langle  \gamma_1 \cdot \! \cdot \! \cdot  \gamma_n   \mid \gamma_1' \cdot \! \cdot \! \cdot  \gamma_m'    \rangle^+_{g,\beta}.\] 
	\end{itemize}
\end{prop}
\textit{Proof.}  
  Let $F$ be the fiber of the forgetful morphism  from Proposition \ref{forg},
\[\mathrm{fg}_{n+1}\colon\Mbar^+_{g,n+1,m}(X,\beta)^{\mathrm{log}} \rightarrow \Mbar^+_{g,n,m}(X,\beta)^{\mathrm{log}},\]
over $
[f\colon (C, \underline{p}, M_C) \rightarrow  (W(m),M_{W(m)})] \in \Mbar^+_{g,n,m}(X,\beta)^{\mathrm{log}}$.
By Proposition \ref{forgobst}, arguments from \cite[Corollary 4.10]{Man} and the projection formula, the fiber $F$ carries a perfect obstruction theory, such that 
\begin{multline} \label{proj}
\langle  \gamma_1 \cdot \! \cdot \! \cdot \gamma_n \cdot  \gamma_{n+1}\mid \gamma'_1 \cdot \! \cdot \! \cdot  \gamma'_m \rangle^{+}_{g,\beta}\\
=([F]^{\mathrm{vir}}\cdot \p(\ev_{n+1})^*\gamma_{n+1})\langle  \gamma_1 \cdot \! \cdot \! \cdot \gamma_n \mid \gamma'_1 \cdot \! \cdot \! \cdot  \gamma'_m \rangle^{+}_{g,\beta}. 
\end{multline}
We therefore have to compute the quantity $[F]^{\mathrm{vir}}\cdot \p(\ev_{n+1})^*\gamma_{n+1}$. For  the unramified map $f$,  we define a finite subset $Z_f \subset  C$, 
\begin{equation} \label{Zf}
Z_f:= \{p \in C \mid f(p)=f(p_i) \text{ for some }i, \text{ or } f(p) \in E_j \text{ for some }j \}.
\end{equation}
Away from the subset $Z_f\subset C$, the fiber $F$ over $f$ is isomorphic to $C$ and has a trivial obstruction theory. While at each point in $Z_f$, we attach certain blown-up projective spaces to $C$, which parametrise lines on end components of FM degenerations.\footnote{Since markings cannot lie in the same fiber of a map, the target sprouts end components over which the source acquires rational tails. These rational tails are mapped isomorphically onto lines in the end components and therefore can  have positive-dimensional moduli, which admit an iterated blow-up description, similar to the one discussed in \cite[Section~4.6]{10A}.} In fact, the exact form of these spaces is not necessary, as we will see below.  We obtain that 
\[F\cong C \cup_{Z_f} Y, \] 
for some space $Y$, such that \[[F]^{\mathrm{vir}}=[C]+A,\] where $A$ is a class of homological degree 2 supported on $Y$. The evaluation morphism $\p(\ev_{n+1})$ restricted to the fiber $F$  is equal to \[\p(f)\colon C \rightarrow \p(TX)\] on the component $C$ and is constant on the component $Y$. Since $\p(\ev_{n+1})$ is constant on $Y$, the pullback of a class by $\p(\ev_{n+1})$ will pair trivially with the class $A$. Hence only the component corresponding to $C$ contributes,
\[[F]^{\mathrm{vir}}\cdot \p(\ev_{n+1})^*\gamma_{n+1}=[C]\cdot \p(f)^*\gamma_{n+1}. \]
  Using Proposition \ref{degreegenus}, we therefore obtain  
\begin{align*}
[F']^{\mathrm{vir}}\cdot\p(\ev_{n+1})^*D& =\deg(\p(f)^*D) =D \cdot \beta, \\
[F']^{\mathrm{vir}}\cdot \p(\ev_{n+1})^*H&= \deg(\p(f)^*H) =2g(C)-2. 
\end{align*}
The claim then follows from (\ref{proj}). 
\qed 
\subsection{The  chamber without markings} \label{secret} If $\epsilon > 2$, then the $\epsilon$-unramification forbids marked points. This implies that $\epsilon_+$-unramified invariants vanish, while $\epsilon_-$-unramified invariants are equal to wall-crossing invariants.   In fact, the wall-crossing formula that compares $\epsilon_-$ and $\epsilon_+$, such that 
\begin{align} \label{sec}
	\begin{split}
		1< &\epsilon_-\leq2, \\
		2<&\epsilon_+ < \infty,
	\end{split}
\end{align}
will allow us to relate relative invariants to absolute ones, which will in turn be used in the proof of the dilaton equation, Proposition \ref{dilaton}. We can also achieve this using a single marking by assigning weight $0$ to all markings except the chosen one. More specifically, if the markings are labeled by the set $\{1, \dots, n+1\}$, we assign weight $1/2$ to the $(n+1)$-th marking and weight $0$ to all the others. The wall-crossing formula then takes the following form.
\begin{thm} \label{relpsi} Assume $\dim(X)=3$ and $k_{n+1}=0, 1$. Then we have the following relation between relative and absolute invariants, 
	
	\begin{multline*} 
		\langle \gamma_1 \psi^{k_1}_1 \cdot \!\cdot \!\cdot \gamma_{n} \psi^{k_{n}}_{n} \cdot  \gamma_{n+1} \psi_{n+1}^{k_{n+1}}\mid \gamma'_1 \Psi^{k'_1}_1\cdot \!\cdot \!\cdot \gamma'_m \Psi^{k'_m}_m   \rangle^{+}_{g,\beta} \\
		= \langle \gamma_1 \psi^{k_1}_1 \cdot \!\cdot \!\cdot \gamma_{n} \psi^{k_{n}}_{n} \mid  \gamma'_1 \Psi^{k'_1}_1\cdot \!\cdot \!\cdot \gamma'_m \Psi^{k'_m}_m \cdot \gamma_{n+1} (\Psi_{m+1}+H)^{k_{n+1}}  \rangle^{+ }_{g,\beta},
	\end{multline*}
where, following Definition \ref{inv}, the last relative insertion on the right is associated to the partition $\eta^{m+1}=(1)$.
\end{thm}

\textit{Proof.} Let $\epsilon_-$ and $\epsilon_+$ be as in (\ref{sec}). The claim is a specialization of Theorem~\ref{wallcrossing} for these values of $\epsilon$, provided we assign the weights to the markings as above. In fact, the situation is even simpler,  as  FM degenerations in the maser space will have at most one strictly constant end component -- the one which contains the $(n+1)$-th marked point. Hence entanglement is not required, and just the calibration suffices. 

 In this case, the $\epsilon_+$-invariants vanish because the corresponding moduli space is empty, as explained in Section \ref{sectionepsilonun}. The only Vertex space contributing to this wall-crossing is the one of weight $1/2$ in the sense of Definition~\ref{PartitionWeight}; namely,
\[
V_{1-\ell,1,(1^\ell)}(X),
\]
where $(1^\ell):=(1,\dots,1)$.
It parametrises $\ell$ copies of $\p^1$, each mapping with degree~$1$, one of which contains the $(n+1)$-th marking.  If $\dim(X)=3$, then by the dimension constraint, the truncated $I$-functions corresponding to $V_{1-\ell,1,(1^\ell)}(X)$ vanish unless $\ell=1$, 
 \[
 \mu_{1-\ell,1,(1^\ell)}(z, \gamma_{n+1} \psi_{n+1}^{k_{n+1}} )=\gamma_{n+1} \cdot \mu_{1-\ell,1,(1^\ell)}(z,  \psi_{n+1}^{k_{n+1}} )=0 \quad \text{ if }  \ell>1.
 \]
Indeed, the codimension of the class $\mu_{1-\ell,1,(1^\ell)}(z,  \psi_{n+1}^{k_{n+1}} )$ is equal to 
\[
2(1-\ell)+k_{n+1}.
\] Hence, if $ k_{n+1}=0,1$, it is negative for $\ell \geq 2$, and the class therefore must vanish; cf. Lemma~\ref{vanishing} and Remark~\ref{remarkmu}.  Note that we used the projection formula for the first equality above.

 Let $\ell=1$, then the $\BC^*$-locus of $V_{0,1}(X)$ is naturally isomorphic to $\p(TX)$, 
\[ V_{0,1}(X)^{\BC^*}\cong \p(TX). \]
The virtual normal bundle of $V_{0,1}(X)^{\BC^*}$  inside $V_{0,1}(X)$ is $TX\otimes \mathbf{z}$, and, with respect to the identification above,   the $\psi$-class on $V_{0,1}(X)^{\BC^*}$ is equal to 
\[ -z+H,\]
where $z$ is equivariant parameter, and $H$ is the hyperplane class. Hence, we have that 
\[
\mu_{0,1}(z, \gamma_{n+1} \psi_{n+1}^{k_{n+1}} )= \gamma_{n+1}  (-z+H)^{k_{n+1}} \in H^*(\BP (TX))[z].  
 \]
Plugging this expression  into the wall-crossing formula from Theorem \ref{wallcrossing}, we obtain the claim. 
\qed

\begin{rmk} \label{remarkmu}
	 In general, the codimension of the class $ \mu_{1-\ell,1,(1^\ell)}(z, \psi^k )$  on the space $\CI_{(1^\ell)} \Sym_X \p(TX)$ is equal to $(\dim-1)(1-\ell)+k$. In particular, for higher powers $k$, there can be contributions from Vertex spaces associated to higher $\ell$.  For $\dim=1,2$, there can be such contributions even for $k=1$; if $\dim=1$, then $\ell$ must be equal to the degree of maps. 
\end{rmk}

\subsection{Dilaton equation}
\begin{prop} \label{dilatonGW} Relative Gromov--Witten invariants satisfy the following dilaton equation, 
		\begin{multline*}
		\langle \gamma_1 \psi^{k_1}_1 \cdot \!\cdot \!\cdot \gamma_n \psi^{k_n}_n \cdot \psi_{n+1}\mid \gamma'_1 \Psi^{k'_1}_1\cdot \!\cdot \! \cdot \gamma'_m \Psi^{k'_m}_m   \rangle^{-}_{g,\beta} \\
		=\Big(2g-2+n+\sum_j \ell(\eta^j)\Big) \langle \gamma_1 \psi^{k_1}_1 \cdot \!\cdot \!\cdot \gamma_n \psi^{k_n}_n \mid \gamma'_1 \Psi^{k'_1}_1\cdot \!\cdot \! \cdot \gamma'_m \Psi^{k'_m}_m  \rangle^{-}_{g,\beta},		
	\end{multline*}
where $\eta^j$ is the ramification profile over the $j$-th exceptional divisor. 
\end{prop}

\textit{Proof.} Exactly the same as in the non-relative Gromov--Witten theory. The morphism forgetting a marking   can be identified with the universal curve, such that we treat the points which map to the exceptional divisors as marked points. However, stabilisation of the target might be required,  see the proof of Proposition \ref{forg} for how this is done.  We refer to \cite[Section 1.5.4]{MP06} for the case of a fixed relative divisor.
\qed 	
\\

With the help of Theorem \ref{relpsi}  and Lemma \ref{vanishing}, we can greatly simplify the proof of the dilaton equation for unramified Gromov--Witten invariants in the case of threefolds for low powers of $\psi$-classes by reducing it to Proposition \ref{dilatonGW}.  However, we expect that the following result can be proved in any dimension and for arbitrary powers of the $\psi$-classes using the same arguments. In general, this would require an analysis of a more complicated wall-crossing formula.

\begin{prop} \label{dilaton} Assume $\dim(X)=3$, $k_{i}=0, 1$ for all $i$, and $\eta^j=(1)$ for all $j$. Then relative unramified Gromov--Witten  invariants satisfy the dilaton equation, 
	\begin{multline*}
	\langle \gamma_1 \psi_1^{k_1} \cdot \!\cdot \!\cdot \gamma_n \psi_n^{k_n} \cdot \psi_{n+1}\mid \gamma'_1 \Psi^{k'_1}_1\cdot \!\cdot \! \cdot \gamma'_m \Psi^{k'_m}_m  \rangle^{+}_{g,\beta} \\
	= (2g-2+n+m) \langle \gamma_1 \psi^{k_1}_1 \cdot \!\cdot \!\cdot \gamma_n \psi^{k_n}_n \mid \gamma'_1 \Psi^{k'_1}_1\cdot \!\cdot \! \cdot \gamma'_m \Psi^{k'_m}_m  \rangle^{+}_{g,\beta}.
	\end{multline*}
	\end{prop}
\textit{Proof.} 
Proving this directly is difficult because the forgetful morphism cannot be identified with the universal curve for unramified maps. We would therefore need an explicit description of the fibers of $\mathrm{fg}_{n+1}$ and of the restrictions of the $\psi$-classes to them; cf.\ the proof of Proposition~\ref{dilatondivisor}. Luckily, we can use the wall-crossing formula, Theorem \ref{wallcrossing},  to reduce the dilaton equation for unramified Gromov--Witten invariants to the dilaton equation for Gromov--Witten invariants, Proposition \ref{dilatonGW}. First, using Theorem \ref{relpsi}, we can move all absolute insertions to the relative ones except the last insertion, 
\begin{multline*}
	\langle \gamma_1 \psi_1^{k_1} \cdot \!\cdot \!\cdot \gamma_n \psi_n^{k_n} \cdot \psi_{n+1} \mid \gamma'_1 \Psi^{k'_1}_1\cdot \!\cdot \!\cdot \gamma'_m \Psi^{k'_m}_m   \rangle^{+}_{g,\beta} \\
	=\langle  \psi_{n+1} \mid \gamma'_1 \Psi^{k'_1}_1\cdot \!\cdot \!\cdot \gamma'_m \Psi^{k'_m}_m \cdot \gamma_1 (\Psi_{m+1}+H)^{k_1}\cdot \!\cdot \!\cdot \gamma_n (\Psi_{m+n}+H)^{k_n}   \rangle^{+}_{g,\beta}. 
	\end{multline*}
We do this because the wall-crossing formula, Theorem \ref{wallcrossing},   does not affect the relative insertions. After relabelling the insertions, it is therefore sufficient to prove the dilaton equation for invariants of the form
\[\langle  \psi_{1} \mid \gamma'_1 \Psi^{k'_1}_1\cdot \!\cdot \!\cdot \gamma'_m \Psi^{k'_m}_m   \rangle^{+}_{g,\beta}.\]
 By Lemma \ref{vanishing},  there will be two types of Vertex spaces involved in the wall-crossing formula, 
\begin{equation*}
 V_{g',0}(X) \quad \text{and} \quad V_{g',1}(X). 
 \end{equation*}
If $g'\neq 0$, the corresponding $\mu$-functions are related by the dilaton equation for moduli spaces of maps to $\p(TX\oplus \BC)$ relative to the divisor $\p(TX)$, 
\begin{equation} \label{dilatoncurves}
 \mu_{g',1}(z,  \psi_1 )= (2g'-2+1) \mu_{g',0}(z,  \emptyset ). 
 \end{equation}
If $g'=0$, then  $\mu_{g',1}(z,  \psi )$ is responsible for the wall-crossing to the chamber with $\epsilon>2$, which we already considered in Section \ref{secret}.

In what follows, we denote 	$\mu_{g,0}(z,\emptyset)$ by	$\mu_{g,0}(z)$. Hence, using (\ref{dilatoncurves}) and Theorem \ref{wallcrossing}, we get the following relation,
\begin{multline} \label{wallcross}
	\langle  \psi_{1} \mid \gamma'_1 \Psi^{k'_1}_1\cdot \!\cdot \!\cdot \gamma'_m \Psi^{k'_m}_m   \rangle^{-}_{g,\beta}-\langle  \psi_{1} \mid \gamma'_1 \Psi^{k'_1}_1\cdot \!\cdot \!\cdot \gamma'_m \Psi^{k'_m}_m   \rangle^{+}_{g,\beta} \\
	=\sum_{\underline{g}}\bigg(\langle  \psi_{1} \mid \gamma'_1 \Psi^{k'_1}_1\cdot \!\cdot \!\cdot \gamma'_m \Psi^{k'_m}_m \cdot \prod^{k}_{i=1}\mu_{g_{i},0}(-\Psi_{m+i})  \rangle^{+}_{g_0,\beta}/k! \\
	+ \Big(\sum^{k}_{i=1} (2g_i-2+1)\Big) \cdot \langle \emptyset \mid \gamma'_1 \Psi^{k'_1}_1\cdot \!\cdot \!\cdot \gamma'_m \Psi^{k'_m}_m \cdot \prod^{k}_{i=1}\mu_{g_{i},0}(-\Psi_{m+i})  \rangle^{+}_{g_0,\beta}/k! \bigg),
\end{multline}
where we sum over ordered partitions 
\[
\underline{g}=(g_0, g_1 ,\dots, g_k)
\] of $g$, such that $g_0<g$, and $g_i>0$ if $i>0$. We now apply induction on the genus. The base case is the minimal value of the genus of a curve in the class $\beta$, which we denote by $g_\beta$. In this case, the wall-crossing must be trivial, otherwise it would contradict the minimality of $g_\beta$. We therefore obtain
\[ \langle  \psi_{1} \mid \gamma'_1 \Psi^{k'_1}_1\cdot \!\cdot \!\cdot \gamma'_m \Psi^{k'_m}_m   \rangle^{-}_{g_\beta,\beta}=\langle  \psi_{1} \mid \gamma'_1 \Psi^{k'_1}_1\cdot \!\cdot \!\cdot \gamma'_m \Psi^{k'_m}_m   \rangle^{+}_{g_\beta,\beta}. \] 
The same holds for  invariants without absolute $\psi$-insertions. Hence, if the dilaton equation holds for Gromov--Witten invariants of genus $g_\beta$, then it holds for unramified Gromov--Witten invariants of genus  $g_\beta$. The former holds by Proposition \ref{dilatonGW}. Consider now an arbitrary genus $g$, and assume that the dilaton equation holds for unramified Gromov--Witten invariants of genus $g'$, such that $g'<g$.  By (\ref{wallcross}), we therefore obtain 

	\begin{multline} 
		\langle  \psi_{1} \mid \gamma'_1 \Psi^{k'_1}_1\cdot \!\cdot \!\cdot \gamma'_m \Psi^{k'_m}_m   \rangle^{-}_{g,\beta}-\langle  \psi_{1} \mid \gamma'_1 \Psi^{k'_1}_1\cdot \!\cdot \!\cdot \gamma'_m \Psi^{k'_m}_m \rangle^{+}_{g,\beta}   \\
		= \sum_{\underline{g}} \bigg((2g_0-2+m+k)\cdot\langle \emptyset \mid \gamma'_1 \Psi^{k'_1}_1\cdot \!\cdot \!\cdot \gamma'_m \Psi^{k'_m}_m \cdot \prod^{k}_{i=1}\mu_{g_{i},0}(-\Psi_{m+i})  \rangle^{+}_{g_0,\beta}/k! \\
		+\Big(\sum^{k}_{i=1} (2g_i-2+1) \Big)\cdot\langle \emptyset \mid \gamma'_1 \Psi^{k'_1}_1\cdot \!\cdot \!\cdot \gamma'_m \Psi^{k'_m}_m \cdot \prod^{k}_{i=1}\mu_{g_{i},0}(-\Psi_{m+i})  \rangle^{+}_{g_0,\beta} /k!\bigg)\\
		 =(2g-2+m)\cdot \sum_{\underline{g}} \langle \emptyset \mid  \gamma'_1 \Psi^{k'_1}_1\cdot \!\cdot \!\cdot \gamma'_m \Psi^{k'_m}_m \cdot \prod^{k}_{i=1}\mu_{g_{i},0}(-\Psi_{m+i})  \rangle^{+}_{g_0,\beta}/k!.
	\end{multline}
By Theorem \ref{wallcrossing}, a similar relation holds for the invariants without absolute $\psi$-insertions, 
\begin{multline}
\langle  \emptyset \mid \gamma'_1 \Psi^{k'_1}_1\cdot \!\cdot \!\cdot \gamma'_m \Psi^{k'_m}_m   \rangle^{-}_{g,\beta}-\langle  \emptyset \mid \gamma'_1 \Psi^{k'_1}_1\cdot \!\cdot \!\cdot \gamma'_m \Psi^{k'_m}_m \rangle^{+}_{g,\beta}   \\
 = \sum_{\underline{g}} \langle \emptyset \mid \gamma'_1 \Psi^{k'_1}_1\cdot \!\cdot \!\cdot \gamma'_m \Psi^{k'_m}_m \cdot \prod^{k}_{i=1}\mu_{g_{i},0}(-\Psi_{m+i})  \rangle^{+}_{g_0,\beta}/k!.
\end{multline}
 We conclude that if the dilaton equation holds for Gromov--Witten invariants of arbitrary genus and for unramified Gromov--Witten invariants of genus $g' < g$, then it also holds for unramified Gromov--Witten invariants of genus $g$. The former follows from Proposition \ref{dilatonGW}. Hence, by induction, the dilaton equation holds for unramified Gromov--Witten invariants of arbitrary genus.
  \qed

\section{Threefolds} \label{sectionthreefolds}
\subsection{Vertex for threefolds}

\begin{lemma} \label{vanishing}
	Assume $\dim(X)=3$. 	If $(n,|\eta|)\neq (0,1), (1,1)$, then
	\[ \mu_{g,n,\eta}(z, \prod^{n}_{i=1} \gamma_i )=0.\] 
	If 	$(n,|\eta|)=(0,1), (1,1)$, then
	\begin{align*}
	\mu_{g,0}(z,\emptyset)&=A_1\cdot z +A_2, \\
		\mu_{g,1}(z,\gamma_1)&=B \cdot \gamma_1, 
	\end{align*}
	where $A_1	, B \in  H^0( \p(TX))$, and  $A_2\in  H^2( \p(TX))$.
\end{lemma}
\textit{Proof.}
The result follows from a simple calculation of virtual dimensions. Let $\ell(\eta)$ be the length of a partition, then
\[\mathrm{vdim}(V_{g,n,\eta}(X))=4|\eta|-(|\eta|-\ell(\eta))+n+3,\]
while 
\[\dim(\CI_{\eta}\Sym_X\p(TX))=2\ell(\eta )+3,\]
hence 
\begin{align*} \dim(\CI_{\eta}\Sym_X\p(TX))- \mathrm{vdim}(V_{g,n,\eta}(X))&= -3|\eta|+\ell(\eta)-n,\\
&\geq -2|\eta|-n.
\end{align*}

 Since the classes $\gamma_i$ come from  $\CI_\eta\Sym_X\p(TX))$, by the projection formula, it is enough to consider the case when $\gamma_i=\mathbb{1}$ for all $i$.  Since $\mu_{g,n,\eta}(z, \mathbb 1)$ is a truncation of  $I_{g,n,\eta}(z, \mathbb 1)$ by non-negative powers, the bound above implies that it is non-zero only if 
 \[
 3-2|\eta|-n\geq 0,
 \]
 where $3$ is the Chow degree of the prefactor $\sum_j \mathrm{c}_{\dim-j}(X)z^{j}$. This can happen only if $(n,|\eta|)=(0,1)$ or $(1,1)$. By the virtual dimension calculation, we also get the second claim.  
\qed
\\

By Lemma \ref{vanishing}, if $\dim(X)=3$, the only contributing invariants are those corresponding to partitions $\eta$ with $|\eta|=1$.  In particular, the source curves of these invariants are connected, and can be attached via only one marking. Hence, the following result also holds for invariants associated to maps with connected sources. In what follows, we denote 	$\mu_{g,0}(z,\emptyset)$ by	$\mu_{g,0}(z)$, and $\mu_{g,1} \in \BQ$ is defined by  $\mu_{g,1}(z,\mathbb{1})=\mu_{g,1} \cdot \mathbb{1}$. 

	\begin{thm} \label{wallthreefold} If $\dim(X)=3$, then
	\begin{multline*}
		\langle \gamma_1 \cdot \! \cdot \! \cdot \gamma_n  \rangle^{-}_{g,\beta} - \langle \gamma_1\cdot \! \cdot \! \cdot\gamma_n  \rangle^{ +}_{g,\beta}= \\
		\sum_{\substack{g_0, \underline{g}_1,   \underline{g}_2 \\ \ell(\underline{g}_1) \leq n}} \frac{1}{k_2!}\binom{n}{k_1} \binom{k_1}{m_1,\dots,m_h} \cdot  \bigg( \prod^{k_1}_{i=1}  \mu_{g_{1,i},1} \cdot  \langle \gamma_1 \cdot \! \cdot \! \cdot \gamma_n \cdot 
		\prod^{k_2}_{i=1} \mu_{g_{2,i},0}(-\psi_{n+i}+H) \rangle^{+}_{g_0,\beta} \bigg),	
	\end{multline*}
	where we sum over unordered partitions $\underline{g}_1$ of length $k_1$, ordered partitions $\underline{g}_2$ of length $k_2$, and integers $g_0$, such that $k_1\leq n$, and 
	\[
	|\underline{g}_1|+|\underline{g}_2|+g_0=g.
	\]  
	The set of integers $\{ m_1,\dots, m_h \}$ are multiplicities of the parts in the partition $\underline{g}_1$. The same holds for invariants associated to maps with connected sources. 
	
\end{thm}
\textit{Proof.} We use Theorem \ref{wallcrossing}, crossing all walls between  $\epsilon=-$ and $\epsilon=+$.   By Lemma \ref{vanishing}, the only non-zero wall-crossing insertions are those for which $(n,|\eta|)=(0,1),(1,1)$. Moreover, again by Lemma \ref{vanishing}, 
\[ \mu_{g,1}(z,\gamma)= \mu_{g,1} \cdot \gamma\] 
for some rational number $\mu_{g,1} \in \BQ$. Let $k_1$ be the number of wall-crossing insertions of the type $\mu_{g,1}$, and $k_2$ be the number of wall-crossing insertions of the type  $\mu_{g,0}$. We move all the wall-crossing insertions of the type $\mu_{g,1}$ to the left, and we also order them according to the order of the corresponding marked points. We denote the marked point corresponding to a wall-crossing insertion of the type $\mu_{g,1}$ by $n_i$, such that $n_{i+1}>n_i$.    There are  
\[ k_1!\cdot\binom{k}{k_1}=\frac{k!}{k_2!}\] ways to make an arbitrary stable partition out of a stable partition which respects the order of marked points and keeps the parts with $N_i\neq \emptyset$ to the left. The wall-crossing terms in Theorem \ref{wallcrossing}  therefore take the following form, 

\begin{equation} 
	\sum_{\vec{g}} \prod^{k_1}_{i=1} \mu_{g_{1,i},1} \Big\langle \prod_{j\in N_0} \gamma_j \ \Big|   \ \prod^{k_1}_{i=1}\gamma_{n_i} \cdot \prod^{k_2}_{i=1} \mu_{g_{2,i},0}\Big( -\Psi_{i} \Big)  \Big\rangle^{+}_{g_0,\beta}\Big/ k_2!,
\end{equation}
where we sum over stable ordered partitions which respect the order of marked points and keeps parts with $N_i\neq \emptyset$ to the left. 

 By Theorem \ref{relpsi}, we can turn all relative insertions for $\eta=(1)$ into absolute insertions. We then move all insertions $\prod^{k_1}_{i=1}\gamma_{n_i}$ back to their initial positions.  Wall-crossing insertions of the type $\mu_{g,1}$ involve the following choices:
\begin{itemize}
	\item a choice of $k_1$ markings out of $n$ ordered markings,
	\item putting an order on an unordered partition $\underline g_1=(g_{1,1}, \dots, g_{1,k_1})$. 
\end{itemize}

These choices contribute in terms of $\binom{n}{k_1}$  and $\binom{k_1}{m_1,\dots,m_h}$ to the formula, where $m_i$ are the multiplicities of the genera in the partition $(g_{1,1}, \dots, g_{1,k_1})$. This leads us to the formula in the statement of the theorem. 
\qed

\subsection{ Hodge integrals} \label{Threefolds} We will now compute $\mu_{g,n}(z)$ in the case when $\dim(X)=3$ and $n=0, 1$, which by Lemma \ref{vanishing} are the only non-vanishing truncated vertex functions for primary insertions.  Let us first look at $V_{g,0}(X)$. Its $\BC^*$-fixed locus admits a simple expression, 
\[ V_{g,0}(X)^{\BC^*}\cong \Mbar_{g,1}\times \p(TX),\]
where $\p(TX)$ parametrises lines $L$ going through $[0,1] \in \p(TX\oplus \BC)$, while the space $\Mbar_{g,1}$ parametrises curves with one marked point $(C,p)$ which are attached to a line $L$  via the marked point $p \in C$ and $[0,1] \in L$. 
Let $\ell_i$ be the Chern roots of $TX$.  With respect to the identification above, we can write $I_{g,0}(z,\emptyset)$  as follows,
\begin{align*}
I_{g,0}(z,\emptyset)& =\int_{\Mbar_{g,1}} \frac{e(\BE^\vee \otimes TX\otimes \mathbf z)}{e(T_pC\otimes T_0L \otimes \mathbf z)}\\
	&=\int_{\Mbar_{g,1}}\frac{\Lambda(\ell_1+z)\cdot \Lambda(\ell_2+z) \cdot  \Lambda(\ell_3+z)}{(z-\psi_1-H)}, 
\end{align*}
where $\BE$ is the \textit{Hodge bundle}, i.e., the vector bundle whose fiber over a curve $(C,p)$ is given by $H^0(C, \Omega^{\mathrm{log}}_C)$, and 
\[\Lambda(\ell_i+z):=\sum^{g}_{j=0}(-1)^{g-j}\lambda_{g-j}(\ell_i+z)^j,\]
such that $\lambda_{j}$ are the Chern classes of $\BE$. We treat $\ell_i$ and $z$ as commuting variables with a relation $\ell_i^3=0$. This quantity is known as a \textit{triple Hodge integral} \cite{PandDeg ,FP}. 

Let us now look at the case $n=1$. By the same  arguments, it admits a similar expression with the difference that the integral takes place on $\Mbar_{g,2}$, 
\[I_{g,1}(z,\mathbb{1})=\int_{\Mbar_{g,2}}\frac{\Lambda(\ell_1+z)\cdot \Lambda(\ell_2+z) \cdot  \Lambda(\ell_3+z)}{(z-\psi_1-H)}.  \]
We can reduce this integral to the one on $\Mbar_{g,1}$ by the string equation.  Let 
\[
\pi\colon \Mbar_{g,2} \rightarrow \Mbar_{g,1}
\] be the natural projection.  The Hodge bundle pullbacks  to the Hodge bundle via $\pi$, hence 
\[ \pi^*\Lambda(\ell_i+z)=\Lambda(\ell_i+z),\]
while $\psi$-classes pullback to $\psi$-classes with a difference given by the boundary divisor parametrizing curves with rational tails that contain both marked points,  denoted by $D_{12} \subset \Mbar_{g,2}$, 
\[  \psi_1^k=\pi^* \psi_1^k+(-D_{12})^{k-1}\cdot D_{12},\]
such that $\pi_*((-D_{12})^{k-1}\cdot D_{12})=\psi_1^{k-1}$. Hence by the projection formula, we get that 
\[ I_{g,1}(z,\mathbb{1})= \frac{I_{g,0}(z,\emptyset)}{z}.\] 
It therefore remains to compute $I_{g,0}(z,\emptyset)$.  This particular Hodge integral can be explicitly evaluated. 
	\begin{prop}[Maximilian Schimpf] \label{Max} We have
	\begin{multline*} \sum_{g>0} u^{2g} \int_{\Mbar_{g,1}}\frac{\Lambda(\ell_1+z)\cdot  \Lambda(\ell_2+z) \cdot  \Lambda(\ell_3+z)}{(z-\psi_1-H)} \\
	=	\left( \frac{\mathrm{sin}(u/2)}{u/2}\right)^2 \left((z-H)+ \log\left( \frac{\mathrm{sin}(u/2)}{u/2}\right) \left(3H+ \mathrm{c}_1 (X) \right)\right) -z+H  +O(1/z).
	\end{multline*}
\end{prop}
\textit{Proof.} The evaluation of this integral is implicitly presented in \cite{FP}.  We do a variable substitution, 
\[z'=z-H,\]
then by the arguments from \cite{N22}, we get that 
\begin{equation*}
	z'+\sum_{g>0} u^{2g} \int_{\Mbar_{g,1}}\frac{\prod^{3}_{i=1}\Lambda(\ell_i+H+z')}{(z'-\psi_1)} =z'\cdot \left(\frac{\sin(u/2)}{u/2}\right)^{\frac{2z'+3H+\mathrm{c}_1(X)}{z'}},
	\end{equation*}
and the claim follows by expanding the expression in powers of $z'$. 
 \qed 
\\

 After taking truncations by non-negative powers of $z$, we obtain expressions for $\mu_{g,0}(z)$ and $\mu_{g,1}$. 
\begin{cor} We have 
	\begin{align*}
		&\sum_{g>0} u^{2g}\mu_{g,0}(z)= \left( \frac{\mathrm{sin}(u/2)}{u/2}\right)^2 \left((z-H)+ \log\left( \frac{\mathrm{sin}(u/2)}{u/2}\right) \left(3H+ \mathrm{c}_1 (X) \right)\right) -z+H, \\ 
		&\sum_{g>0}u^{2g}\mu_{g,1}= \left( \left( \frac{\mathrm{sin}(u/2)}{u/2}\right)^2 -1\right).
		\end{align*}
	\end{cor}
\noindent We are now ready to evaluate the wall-crossing formula for threefolds. 
\begin{thm} \label{pand} If $\dim(X)=3$, then there is an explicit wall-crossing formula relating  Gromov--Witten and unramified Gromov--Witten invariants of $X$ in a class $\beta \in H_2(X)$, 
	\[ \sum_g	\langle \gamma_1 \cdot \! \cdot \! \cdot\gamma_n  \rangle^{-}_{g,\beta} u^{2g-2}=  \sum_g	\langle \gamma_1  \cdot \! \cdot \! \cdot\gamma_n  \rangle^{+}_{g,\beta} u^{2g-2} \left(\frac{\mathrm{sin}(u/2)}{u/2}\right)^{2g-2+\beta\cdot \mathrm{c}_1(X)}  . \] 
	In particular, if $\beta\cdot \mathrm{c}_1(X)>0$, or  $\beta\cdot \mathrm{c}_1(X)=0$ and $\beta$ is primitive, we obtain that unramified Gromov--Witten invariants and Gopakumar--Vafa invariants are equal. 
\end{thm}
\textit{Proof.} For simplicity, we define 
\begin{align*} 
t_1&:=  \left(\frac{\mathrm{sin}(u/2)}{u/2}\right)^2-1  ,   &t_2:= \left(\frac{\mathrm{sin}(u/2)}{u/2}\right)^2  \log\left(\frac{\mathrm{sin}(u/2)}{u/2}\right), \\
A_1&:=-\psi, &t'_2:=\log\left(\frac{\mathrm{sin}(u/2)}{u/2}\right), \\
 & & A_2:=3H+ \mathrm{c}_1(X). 
\end{align*}

Let us first disregard the wall-crossing insertions of the type $\mu_{g,1}$, as by Theorem \ref{wallthreefold} they appear as  coefficients in front of the  brackets. Hence, by Proposition \ref{Max} and the wall-crossing formula from Theorem \ref{wallthreefold}, for each $g$ we have to evaluate the following series:

\begin{multline} \label{eq}
\sum_{k \geq 0} \langle \gamma_1\cdot \! \cdot \! \cdot \gamma_n  \cdot \underbrace{( t_1A_1+ t_2A_2) \cdot \! \cdot \! \cdot (t_1A_1+ t_2A_2)}_{k\text{ times}} \rangle^+_{g,\beta}/ k! \\
=\sum_{a_1,a_2\geq 0} \frac{1}{a_1!a_2!}\langle \gamma_1\cdot \! \cdot \! \cdot \gamma_n \cdot \underbrace{t_1A_1 \cdot \! \cdot \! \cdot t_1A_1}_{a_1 \text{ times}} \cdot \underbrace{ t_2A_2 \cdot \! \cdot \! \cdot  t_2A_2}_{a_2 \text{ times}} \rangle^+_{g,\beta},
\end{multline}
where we used that the bracket is symmetric in its entries. Denote 
\[
\langle \gamma_1\cdot \! \cdot \! \cdot \gamma_n \cdot  (t_1A_1)^{a_1} \cdot (t_2A_2)^{a_2} \rangle^+_{g,\beta} :=\langle \gamma_1\cdot \! \cdot \! \cdot \gamma_n \cdot \underbrace{t_1A_1 \cdot \! \cdot \! \cdot t_1A_1}_{a_1 \text{ times}} \cdot \underbrace{ t_2A_2 \cdot \! \cdot \! \cdot  t_2A_2}_{a_2 \text{ times}} \rangle^+_{g,\beta}.
\]
In order to evaluate these invariants,  we apply the divisor and dilaton equations from Proposition \ref{dilatondivisor} and Proposition \ref{dilaton}. More precisely, we first apply the dilaton equation and then the divisor one, because the dilaton equation behaves consistently in the presence of other $\psi$-classes.  We obtain that
\begin{multline}
\frac{1}{a_1!a_2!} \langle \gamma_1\cdot \! \cdot \! \cdot \gamma_n \cdot  (t_1A_1)^{a_1} \cdot (t_2A_2)^{a_2} \rangle^+_{g,\beta} \\
=\frac{(m_2)^{a_2} (m_1-(a_1+a_2))!}{a_1! a_2!(m_1-a_2)!}t_1^{a_1}t_2^{a_2}\langle \gamma_1\cdot \! \cdot \! \cdot \gamma_n \rangle^+_{g,\beta}, 
\end{multline}
where 
\begin{align*}
m_1&:=-(2g-2+n), \\
 m_2&:=\beta \cdot A_2=3(2g-2)+\mathrm{c}_1(S)\cdot \beta .
\end{align*}
We then sum over $a_1$, keeping $a_2$ fixed, to obtain  
\begin{multline*} (\ref{eq})=\sum_{a_2} \sum_{a_1}\frac{(m_2)^{a_2} (m_1-(a_1+a_2)!}{a_1! a_2!(m_1-a_2)!}t_1^{a_1}t_2^{a_2} \langle \gamma_1\cdot \! \cdot \! \cdot \gamma_n \rangle^+_{g,\beta}  \\ =\sum_{a_2} \frac{(m_2t_2)^{a_2}}{a_2!}(1+t_1)^{m_1-a_2} \langle \gamma_1\cdot \! \cdot \! \cdot \gamma_n \rangle^+_{g,\beta}.
	\end{multline*}
Using that $t_2=(1+t_1)\cdot t_2'$, we have
 \begin{multline*} (\ref{eq})=\sum_{a_2} \frac{(m_2t_2)^{a_2}}{a_2!}(1+t_1)^{m_1-a_2}\langle \gamma_1\cdot \! \cdot \! \cdot \gamma_n \rangle^+_{g,\beta} \\
 	=\left(\sum_{a_2}\frac{(m_2t')^{a_2}}{a_2!}\right)  (1+t_1)^{m_1}\langle \gamma_1\cdot \! \cdot \! \cdot \gamma_n \rangle^+_{g,\beta}  =
 	\exp(m_2t_2')(1+t_1)^{m_1} \langle \gamma_1\cdot \! \cdot \! \cdot \gamma_n \rangle^+_{g,\beta} .
 	\end{multline*}
Finally, we substitute into the above formula the values of $m_i$ and $t_i$ to obtain 
\begin{multline*}(\ref{eq})= \left( \frac{\mathrm{sin}(u/2)}{u/2}\right)^{3(2g-2)+ \beta \cdot \mathrm{c}_1(X) } \left( \frac{\mathrm{sin}(u/2)}{u/2}\right)^{-2(2g-2)-2n} \langle \gamma_1\cdot \! \cdot \! \cdot \gamma_n \rangle^+_{g,\beta} \\ =\left( \frac{\mathrm{sin}(u/2)}{u/2}\right)^{2g-2-2n+ \beta\cdot \mathrm{c}_1(X)}\langle \gamma_1\cdot \! \cdot \! \cdot \gamma_n \rangle^+_{g,\beta}.
	\end{multline*}
Let us now incorporate the contributions from insertions of the type $\mu_{g,1}$. By Theorem \ref{wallthreefold} and Proposition \ref{Max}, the full wall-crossing contribution of genus $g$ is therefore equal to 
\[\sum^{n}_{i=0} \binom{n}{i} \left(\left(\frac{\mathrm{sin}(u/2)}{u/2}\right)^2-1 \right)^{i}\cdot \left( \frac{\mathrm{sin}(u/2)}{u/2}\right)^{2g-2-2n+ \beta \cdot \mathrm{c}_1(X)}\langle \gamma_1\cdot \! \cdot \! \cdot \gamma_n \rangle^+_{g,\beta},\]
which simplifies to 
\begin{multline}
	\left(\frac{\mathrm{sin}(u/2)}{u/2}\right)^{2n} \cdot   \left( \frac{\mathrm{sin}(u/2)}{u/2}\right)^{2g-2-2n+ \beta \cdot \mathrm{c}_1(X)}\langle \gamma_1\cdot \! \cdot \! \cdot \gamma_n \rangle^+_{g,\beta}\\
	 =\left( \frac{\mathrm{sin}(u/2)}{u/2}\right)^{2g-2+\beta \cdot \mathrm{c}_1(X)}\langle \gamma_1\cdot \! \cdot \! \cdot \gamma_n \rangle^+_{g,\beta}.
	\end{multline}
We thus arrive at the desired formula.
\qed 

\bibliographystyle{amsalpha}
\bibliography{uGWJAMS}
	\end{document}